# Free boundary regularity for harmonic measures and Poisson kernels

By Carlos E. Kenig and Tatiana Toro*

## 1. Introduction

One of the basic aims of this paper is to study the relationship between the geometry of "hypersurface like" subsets of Euclidean space and the properties of the measures they support. In this context we show that certain doubling properties of a measure determine the geometry of its support. A Radon measure is said to be doubling with constant $C$ if $C$ times the measure of the ball of radius $r$ centered on the support is greater than the measure of the ball of radius $2r$ and the same center. We prove that if the doubling constant of a measure on $\mathbf{R}^{n+1}$ is close to the doubling constant of the $n$-dimensional Lebesgue measure then its support is well approximated by $n$-dimensional affine spaces, provided that the support is relatively flat to start with. Primarily we consider sets which are boundaries of domains in $\mathbf{R}^{n+1}$. The $n$-dimensional Hausdorff measure may not be defined on the boundary of a domain in $\mathbf{R}^{n+1}$. Thus we turn our attention to the harmonic measure which is well behaved under minor assumptions (see Section 3). We obtain a new characterization of *locally flat domains* in terms of the doubling properties of their harmonic measure (see Section 3). Along these lines we investigate how the "weak" regularity of the Poisson kernel of a domain determines the geometry of its boundary. Sections 5 and 6 pursue this goal, as in Alt and Caffarelli's work (see [AC], [C1], [C2]), and also Jerison's [J]. In both cases the goal is to prove that, under the appropriate technical conditions at "flat points" of the boundary, the oscillation of the Poisson kernel controls the oscillation of the unit normal vector. The difference between our work and the work in [AC] is that we measure the oscillation in an integral sense (BMO estimates) while they do so in a pointwise sense (Hölder estimates).

In [KT1] we studied a boundary regularity problem. Namely, we proved that if the boundary of a domain $\Omega \subset \mathbf{R}^{n+1}$ can be well approximated by

*The first author was partially supported by the NSF. The second author was partially funded by an NSF Postdoctoral Fellowship and an Alfred P. Sloan Research Fellowship.

$n$-dimensional affine spaces then the harmonic measure behaves like the $n$-dimensional Lebesgue measure from a doubling point of view. We also showed that if the unit normal vector to $\partial\Omega$ has small mean oscillation then so does the logarithm of the Poisson kernel. The results discussed in this paper provide answers to the corresponding free boundary regularity problems (i.e. the converse problems). In our context the doubling character of a measure or the mean oscillation of the logarithm of its density, and the flatness of a set or the mean oscillation of its unit normal vector replace stronger notions of regularity.

We now introduce formally the notions of doubling and flatness and sketch briefly the contents of each section of the paper.

*Definition* 1.1. Let $\Sigma \subset \mathbf{R}^{n+1}$ be a locally compact set, and let $\delta > 0$. We say that $\Sigma$ is *$\delta$-Reifenberg flat* if for each compact set $K \subset \mathbf{R}^{n+1}$, there exists $R_K > 0$ such that for every $Q \in K \cap \Sigma$ and every $r \in (0, R_K]$ there exists an $n$-dimensional plane $L(r,Q)$ containing $Q$ such that

$$\frac{1}{r} D[\Sigma \cap B(r,Q), L(r,Q) \cap B(r,Q)] \leq \delta. \tag{1.1}$$

Here $B(r,Q)$ denotes the $(n+1)$-dimensional ball of radius $r$ and center $Q$, and $D$ denotes the Hausdorff distance.

Recall that for $A, B \subset \mathbf{R}^{n+1}$,

$$D[A,B] = \max\{\sup\{d(a,B) : a \in A\},\ \sup\{d(b,A) : b \in B\}\}.$$

Thus $D[A,B] \leq \delta$ is equivalent to $A \subset (B;\delta)$ and $B \subset (A;\delta)$, where for a set $A \subset \mathbf{R}^{n+1}$, $(A;\delta)$ denotes the $\delta$ neighborhood of $A$. Note that the previous definition is only significant for $\delta > 0$ small. Thus we assume that $\delta \in (0, \frac{1}{4\sqrt{2}})$ whenever we talk about $\delta$-Reifenberg flat sets. The relevance of the constant $\frac{1}{4\sqrt{2}}$ will become clear in the proof of Theorem 2.1.

We let

$$\theta(r,Q) = \inf_L \left\{ \frac{1}{r} D[\Sigma \cap B(r,Q), L \cap B(r,Q)] \right\}, \tag{1.2}$$

where the infimum is taken over all $n$-planes containing $Q$. Our work requires uniform control of several quantities on compact sets; thus for each compact set $K \subset \mathbf{R}^{n+1}$ we define

$$\theta_K(r) = \sup_{Q \in \Sigma \cap K} \theta(r,Q). \tag{1.3}$$

This notion was initially introduced by Reifenberg who proved the following remarkable theorem:

THEOREM ([Mo], [R]). *There exists $\delta > 0$ depending only on $n$ so that if $\Sigma$ is $\delta$-Reifenberg flat then locally $\Sigma$ is a topological disc.*

From his proof it is easy to see that the embedding is tame; that is, the local homeomorphism constructed extends to a local homeomorphism of a tubular neighborhood of $\Sigma$ (see [T]). Moreover the proof also shows that this homeomorphism yields a local Hölder parametrization for $\Sigma$. Thus there exists $\beta \in (0,1)$ depending on $\delta$ such that $\Sigma$ is a $C^{0,\beta}$ $n$-dimensional submanifold. Furthermore $\beta$ tends to 1 as $\delta$ approaches 0.

*Definition* 1.2. When $\Sigma \subset \mathbf{R}^{n+1}$, we say that $\Sigma$ is *Reifenberg flat with vanishing constant* if it is $\delta$-Reifenberg flat for some $\delta \in (0, \frac{1}{4\sqrt{2}})$ and for each compact set $K \subset \mathbf{R}^{n+1}$

$$\lim_{r\to 0} \theta_K(r) = 0.$$

*Remark.* Consider $\varphi : \mathbf{R}^n \to \mathbf{R}$, $\varphi \in \lambda_*$, the little-o Zygmund class (see [St] for a definition). Then $\operatorname{graph}\varphi = \{(x,t) \in \mathbf{R}^{n+1} : t = \varphi(x)\}$ is Reifenberg flat with vanishing constant. In fact $\varphi$ is well approximated by affine functions, and therefore its graph is well approximated by $n$-planes. On the other hand there are examples (see [Z]) of such functions $\varphi$ which are almost nowhere differentiable (note, for instance that $\varphi : \mathbf{R} \to \mathbf{R}$ defined by $\varphi(x) = \sum_{k\geq 1} \frac{\cos(2^k x)}{2^k\sqrt{k}}$ belongs to $\lambda_*$, is continuous and almost nowhere differentiable). Thus $\operatorname{graph}\varphi$ does not admit a tangent plane almost everywhere, and hence $\operatorname{graph}\varphi$ is not rectifiable (see [Si] for a definition).

*Definition* 1.3. Let $\mu$ be a Radon measure on $\mathbf{R}^{n+1}$. We say that $\mu$ is a *doubling measure* if for each compact set $K \subset \mathbf{R}^{n+1}$, there exist constants $C = C_K > 1$ and $R_K > 0$ so that for every $Q \in \operatorname{spt}\mu \cap K$ and every $r \in (0, R_K]$

$$\mu(B(2r,Q)) \leq C\mu(B(r,Q)).$$

Here $\operatorname{spt}\mu$ denotes the support of the measure $\mu$. We now quantify further the doubling character of a Radon measure.

*Definition* 1.4. Let $\mu$ be a doubling Radon measure in $\mathbf{R}^{n+1}$, and let $\varepsilon > 0$. We say that $\mu$ is $\varepsilon$-*approximately optimally doubling* if for each compact set $K \subset \mathbf{R}^{n+1}$, and for each $\tau \in (0,1)$ there exists $R = R(\tau, K) > 0$ so that for $r \in (0, R]$, and $Q \in K \cap \operatorname{spt}\mu$

$$\frac{\tau^n}{1+\varepsilon} \leq \frac{\mu(B(r\tau,Q))}{\mu(B(r,Q))} \leq (1+\varepsilon)\tau^n.$$

*Notation.* For $a, b \in (0,\infty)$, and $\varepsilon > 0$ we write

$$a \underset{\varepsilon}{\sim} b \quad \text{if and only if} \quad \frac{1}{1+\varepsilon} \leq \frac{a}{b} \leq 1+\varepsilon.$$

*Definition* 1.5. Let $\mu$ be a doubling Radon measure in $\mathbf{R}^{n+1}$. We say that $\mu$ is *asymptotically optimally doubling* if for each compact set $K \subset \mathbf{R}^{n+1}$ such that $K \cap \operatorname{spt}\mu \neq \emptyset$, and for each $\tau \in (0,1]$

$$\lim_{r\to 0} \inf_{Q\in K\cap \operatorname{spt}\mu} \frac{\mu(B(r\tau,Q))}{\mu(B(r,Q))} = \lim_{r\to 0} \sup_{Q\in K\cap \operatorname{spt}\mu} \frac{\mu(B(r\tau,Q))}{\mu(B(r,Q))} = \tau^n.$$

*Definition* 1.6. A measure $\mu$ in $\mathbf{R}^{n+1}$ is said to be Ahlfors regular if there exists $C > 1$ such that for $Q \in \operatorname{spt}\mu$ and $r \in (0, \operatorname{diam}\operatorname{spt}\mu)$

$$C^{-1}r^n \leq \mu(B(r,Q)) \leq Cr^n.$$

*Definition* 1.7. A nonzero Radon measure $\nu$ in $\mathbf{R}^{n+1}$ is called a *uniform measure* if there exists $C > 0$ such that for every $Q \in \operatorname{spt}\nu$, and every $r \in (0,\infty)$

$$\nu(B(r,Q)) = Cr^n.$$

In Section 2 we show that if $\mu$ is an $\varepsilon$-approximately optimally doubling measure in $\mathbf{R}^{n+1}$ (with $\varepsilon$ small) whose support $\Sigma$ is $\delta$-Reifenberg flat then $\Sigma$ is a $C^{0,\beta}$ $n$-dimensional submanifold for some $\beta > 0$ that only depends on $n$ and $\varepsilon$. In Section 2 we also prove the following theorem (see Theorem 2.1):

THEOREM. *Let $\mu$ be an asymptotically optimally doubling measure. Then its support $\Sigma \subset \mathbf{R}^{n+1}$ is also Reifenberg flat with vanishing constant if $n = 1, 2$ and if $n \geq 3$ and $\Sigma$ is $\delta$-Reifenberg flat, then $\Sigma$ is also Reifenberg flat with vanishing constant.*

For $n \geq 3$ the extra assumption of $\Sigma$ being $\delta$-Reifenberg flat is in fact needed; see [KP] and Proposition 3.1. To prove these results we introduce the notion of pseudo-tangent measure (see Definition 2.2). We then prove that the pseudo-tangent measures of an asymptotically optimally doubling measure are uniform (see Theorem 2.2). Kowalski and Preiss (see [KP]) showed that the support of a uniform measure is either a plane, or a cone (only if $n \geq 3$). The assumption that $\Sigma$ is $\delta$-Reifenberg flat allows us to rule out the cone. This is enough to conclude that $\Sigma$ is Reifenberg flat with vanishing constant.

*Remark.* David and Semmes (see [DS, Part III, Ch. 5]) looked at the relationship between the doubling properties of an Ahlfors regular measure and the regularity of its support. Their doubling condition involves a comparison with the Lebesgue measure in a large set of points. Under their assumptions one concludes that the support is uniformly rectifiable. Note that since Reifenberg flat sets with vanishing constant are not rectifiable in general, their results do not apply in this setting.

In Section 3 we study the special case when $\mu$ is the harmonic measure of a domain $\Omega$ and $\operatorname{spt}\mu = \partial\Omega$. We show that if $\Omega \subset \mathbf{R}^{n+1}$ is a domain with the separation property (see Definition 1.9 below), whose boundary is $\delta$-Reifenberg flat (this condition is only necessary if $n \geq 3$) then $\partial\Omega$ is Reifenberg flat with vanishing constant if and only if its harmonic measure is asymptotically optimally doubling (here we use the results in [KT1] for the only if direction). We also give examples of domains in $\mathbf{R}^{n+1}$, for $n \geq 3$, which are not $\delta$-Reifenberg flat but whose harmonic measure is uniform, and hence asymptotically optimally doubling.

In order to state the main results in Sections 4, 5 and 6 we need to define the space $\mathrm{VMO}(\partial\Omega)$ for $\Omega \subset \mathbf{R}^{n+1}$, a set of locally finite perimeter (see [EG] for a precise definition) whose boundary is Ahlfors regular (i.e. the surface measure $\sigma$ of $\partial\Omega$ is Ahlfors regular). Moreover we need to introduce the notion of the "separation property."

*Definition* 1.8. Let $\Omega \subset \mathbf{R}^{n+1}$ be a set of locally finite perimeter whose boundary is Ahlfors regular. For $f \in L^2_{\mathrm{loc}}(d\sigma)$ we say that $f \in \mathrm{VMO}(\partial\Omega)$ if for each compact set $K \subset R^{n+1}$

$$\lim_{r\to 0} \sup_{Q\in\partial\Omega\cap K} \|f\|_*(B(r,Q)) = 0,$$

where

$$\|f\|_*(B(r,Q)) = \Big(\sup_{0<s\leq r} \frac{1}{\sigma(\Delta(s,Q))} \int_{\Delta(s,Q)} |f - f_{s,Q}|^2 d\sigma\Big)^{\frac{1}{2}},$$

$\Delta(s,Q) = \partial\Omega \cap B(s,Q)$, and $\sigma = \mathcal{H}^n \llcorner \partial\Omega$ denotes the surface measure of the boundary.

*Definition* 1.9. Let $\Omega \subset \mathbf{R}^{n+1}$. We say that $\Omega$ has the *separation property* if for each compact set $K \subset \mathbf{R}^{n+1}$ there exists $R > 0$ such that for $Q \in \partial\Omega \cap K$ and $r \in (0, R]$ there exists an $n$-dimensional plane $\mathcal{L}(r,Q)$ containing $Q$ and a choice of unit normal vector to $\mathcal{L}(r,Q)$, $\overrightarrow{n}_{r,Q}$ satisfying

$$\mathcal{T}^+(r,Q) = \{X = (x,t) = x + t\overrightarrow{n}_{r,Q} \in B(r,Q) : x \in \mathcal{L}(r,Q),\ t > \frac{1}{4}r\} \subset \Omega,$$

and

$$\mathcal{T}^-(r,Q) = \{X = (x,t) = x + t\overrightarrow{n}_{r,Q} \in B(r,Q) : x \in \mathcal{L}(r,Q),\ t < -\frac{1}{4}r\} \subset \Omega^c.$$

Moreover if $\Omega$ is an unbounded domain we also require that $\mathbf{R}^{n+1}\backslash\partial\Omega$ divide $\mathbf{R}^{n+1}$ into two distinct connected components $\Omega$ and $\operatorname{int}\Omega^c \neq \emptyset$.

The notation $(x,t) = x + t\overrightarrow{n}_{r,Q}$ is used to denote a point in $\mathbf{R}^{n+1}$. The first component, $x$, of the pair belongs to an $n$-dimensional affine space whose unit normal vector is $\overrightarrow{n}_{r,Q}$. The second component $t$ belongs to $\mathbf{R}$. From the context it will always be clear what affine hyperplane $x$ belongs to, and what the orientation of the unit normal vector is.

*Remark* 1.1. The separation property required in this definition seems slightly stronger than the one used in [KT1]. Nevertheless a similar argument to the one used in the proof of Proposition 2.2 in [KT1] guarantees that the two definitions are equivalent provided that $\partial\Omega$ is a $\delta$-Reifenberg flat set, for $\delta$ small enough. Moreover it was noted in [KT1] that without loss of generality we could assume that $\mathcal{L}(r,Q) = L(r,Q)$, where $L(r,Q)$ is taken as in Definition 1.1. In particular we obtain that

$$\{X = (x,t) = x + t\overrightarrow{n}_{r,Q} \in B(r,Q) : x \in L(r,Q),\ t > 2\delta r\} \subset \Omega,$$

and

$$\{X = (x,t) = x + t\overrightarrow{n}_{r,Q} \in B(r,Q) : x \in L(r,Q),\ t < -2\delta r\} \subset \Omega^c.$$

It was also shown in [KT1] that there exists $\delta_n < \frac{1}{4\sqrt{2}}$ a small constant depending only on the dimension $n$, such that a domain with the separation property, whose boundary is a $\delta$-Reifenberg flat set, with $\delta \in (0,\delta_n)$ is an NTA domain.

*Definition* 1.10. Let $\delta \in (0,\delta_n)$, and $\Omega \subset \mathbf{R}^{n+1}$. We say that $\Omega$ is a *$\delta$-Reifenberg flat domain* or a *Reifenberg flat domain* if $\Omega$ has the separation property and $\partial\Omega$ is $\delta$- Reifenberg flat. Moreover if $\Omega$ is an unbounded domain we also require that

$$\sup_{r>0}\sup_{Q\in\partial\Omega} \theta(r,Q) < \delta_n.$$

When we consider $\delta$-Reifenberg flat domains in $\mathbf{R}^{n+1}$ we assume that $\delta \in (0,\delta_n)$, in order to insure that we are working on NTA domains. On the other hand when discussing $\delta$-Reifenberg flat sets it is enough to assume that $\delta \in (0, \frac{1}{4\sqrt{2}})$.

*Definition* 1.11. Let $\delta \in (0,\delta_n)$. A set of locally finite perimeter $\Omega$ (see [EG]) is said to be a *$\delta$-chord arc domain* or a *chord arc domain with small constant* if $\Omega$ is a $\delta$-Reifenberg flat domain, $\partial\Omega$ is Ahlfors regular and for each compact set $K \subset \mathbf{R}^{n+1}$ there exists $R > 0$ so that

$$\sup_{Q\in\partial\Omega\cap K} \|\overrightarrow{n}\|_*(B(R,Q)) < \delta.$$

Here $\overrightarrow{n}$ denotes the unit normal vector to the boundary.

We only use the notation $\delta$-Reifenberg flat domain or $\delta$-chord arc domain when we want to emphasize the dependence on $\delta$; otherwise we simply refer to them as Reifenberg flat domains or chord arc domains with small constant.

In Section 4 we provide some new characterizations of chord arc domains with small (enough) constant in terms of the properties of the surface measure of the boundary. In Sections 5 and 6 we provide a characterization of chord arc domains with small (enough) constant in terms of the doubling properties of their harmonic measure and the oscillation of the logarithm of the corresponding Poisson kernel. Let $\Omega$ be a chord arc domain with small (enough) constant, and $X \in \Omega$; then the harmonic measure with pole at $X$, $\omega^X$ (see Section 3 for a definition), and $\sigma$ are mutually absolutely continuous. The Radon-Nikodym theorem insures that the Poisson kernel $k_X(Q) = \frac{d\omega^X}{d\sigma}(Q) = \frac{\partial G(X,-)}{\partial n}(Q) \in L^1_{\mathrm{loc}}(d\sigma)$. Here $G(X,-)$ denotes the Green's function of $\Omega$ with pole at $X$ and $\frac{\partial}{\partial n}$ denotes the normal derivative at the boundary. We prove that if $\omega^X$ is asymptotically optimally doubling and $\log k_X$ has vanishing mean oscillation then the unit normal vector to $\partial\Omega$ also has vanishing mean oscillation (i.e. $\Omega$ is a chord arc domain with vanishing constant; see Definition 4.1). Jerison's paper [J] introduced this endpoint problem, but treated it under more restrictive assumptions, namely that the boundary is given locally as a Lipschitz graph, and that the normal derivative data is continuous as opposed to having vanishing mean oscillation. His paper is based on the work of Jerison-Kenig [JK2] and Alt-Caffarelli [AC]. There is an error in Lemma 4 of Jerison's paper. Nevertheless we still make considerable use of the ideas in [J].

It is interesting to compare our result with the following one by Alt and Caffarelli [AC]. In both cases the oscillation of the logarithm of the Poisson kernel controls the oscillation of the unit normal.

THEOREM ([AC]). *Assume that*:

1. $\Omega \subset \mathbf{R}^{n+1}$ *is a set of locally finite perimeter whose boundary is Ahlfors regular;*
2. $\Omega \subset \mathbf{R}^{n+1}$ *is a $\delta$-Reifenberg flat domain for some $\delta > 0$ small enough*;
3. $\log k_X \in C^{0,\beta}$ *for some* $\beta \in (0,1)$.

*Then $\Omega$ is a $C^{1,\alpha}$ domain for some $\alpha \in (0,1)$ which depends on $\beta$ and $n$.*

Jerison showed that $\alpha = \beta$ (see [J]). We would like to emphasize that the hypotheses 1 and 2 above are necessary. Keldysh and Lavrentiev constructed a domain in $\mathbf{R}^2$ whose boundary is rectifiable but not Ahlfors regular, whose Poisson kernel is identically equal to 1 and which is not $C^1$. Moreover there are examples of domains in $\mathbf{R}^2$ whose boundary is Reifenberg flat with vanishing constant, rectifiable but not Ahlfors regular, for which the logarithm of the Poisson kernel is Hölder continuous and which are not even $C^1$ domains (see [Du]). Furthermore if $n \geq 2$ there are examples of domains satisfying 1 and 3, whose boundaries are not $C^1$, which contain a neighborhood of the vertex that is a double cone (see [AC, Remark 3.2]).

In Section 5 we prove that if $\Omega$ is a chord arc domain with small (enough) constant such that $\omega^X$ is asymptotically optimally doubling and $\log k_X \in$

VMO($\partial\Omega$), then for each compact set $K$ there exists $r > 0$ so that for $Q \in \partial\Omega \cap K$, $\partial\Omega \cap B(r,Q) = \mathcal{A}(r,Q) \cup \mathcal{F}(r,Q)$, where $\mathcal{F}(r,Q)$ is a set of very small surface measure. On $\mathcal{A}(r,Q)$ the pointwise oscillation of the unit normal is controlled, and the Poisson kernel of $\Omega$ (with appropriate pole) can be compared to the Poisson kernel of a half space. In Section 6 using an integration by parts in Rellich's identity (see [J], [JK3] and [KT2]), in both $\partial\Omega$ and the boundary of the half space mentioned above, we gain control on the mean oscillation of the unit normal vector. This is used to prove a decay-type estimate which is the main ingredient in the proof of the following theorem (which also uses Theorem 4.2, the Main Theorem in Section 5 and Theorem 6.3).

THEOREM. *Assume that*:

1. $\Omega \subset \mathbf{R}^{n+1}$ *is a chord arc domain with small* (*enough*) *constant*;
2. $\omega^X$ *is asymptotically optimally doubling*;
3. $\log k_X \in \mathrm{VMO}(\partial\Omega)$.

*Then* $\Omega$ *is a chord arc domain with vanishing constant*; *i.e.*, $\overrightarrow{n} \in \mathrm{VMO}(\partial\Omega)$.

*Remark.* In [KT1] it is shown that if $\Omega$ is a chord arc domain with vanishing constant, 2 and 3 above hold. Thus these 2 theorems characterize chord arc domains with vanishing constant in terms of the behavior of their harmonic measure and their Poisson kernel.

These results should be also be compared with Pommerenke's theorem [P]:

THEOREM ([P]). *Let* $\Omega \subset \mathbf{R}^2$ *be a chord arc domain. Then* $\Omega$ *is a chord arc domain with vanishing constant if and only if* $\log k_X \in \mathrm{VMO}(\partial\Omega)$.

Recall that a domain $\Omega \subset \mathbf{R}^2$ is called a chord arc domain if the distance along the boundary and the Euclidean distance in $\mathbf{R}^2$ are equivalent, i.e. there exists $\kappa > 0$ such that for any $P, Q \in \partial\Omega$, $|P-Q| \le d(P,Q) \le (1+\kappa)|P-Q|$. Here $d$ denotes the distance along $\partial\Omega$. In the case when $\kappa$ is small this definition of chord arc domain is equivalent to the definition of chord arc domain with small constant given above (see Section 4).

CONJECTURE. *Assume that*:

1. $\Omega \subset \mathbf{R}^{n+1}$ *is a set of locally finite perimeter whose boundary is Ahlfors regular*;
2. $\Omega \subset \mathbf{R}^{n+1}$ *is a* $\delta$-*Reifenberg flat domain for some* $\delta > 0$ *small enough*;
3. $\log k_X \in \mathrm{VMO}(\partial\Omega)$.

*Then* $\Omega$ *is a chord arc domain with vanishing constant*; *i.e.*, $\overrightarrow{n} \in \mathrm{VMO}(\partial\Omega)$.

Recently we have succeeded in proving a weaker version of the conjecture ([KT2]). Namely, if conditions (1) and (2) above are replaced by the statement: $\Omega \subset \mathbf{R}^{n+1}$ is a chord arc domain with small (enough) constant, then the above conclusion holds

*Acknowledgment.* Part of this work was carried out while the first author was visiting Princeton University and the second author was visiting the Institute for Mathematics and its Applications. They wish to thank Princeton University and IMA for their hospitality.

## 2. Doubling and flatness

In this section we study the relationship between the doubling properties of a Radon measure and the flatness of its support. Roughly speaking, under the appropriate hypothesis a set of codimension 1 is Reifenberg flat with vanishing constant if it supports a Radon measure which is asymptotically optimally doubling. Recently we showed that the converse is also true. Namely, every Reifenberg flat set with vanishing constant supports an asymptotically optimally doubling measure (see [DKT]). In order to prove the main theorem in this section (Theorem 2.1) we first prove that if $\mu$ is an asymptotically optimally doubling measure in $\mathbf{R}^{n+1}$, then all its pseudo-tangent measures are uniform (see Definitions 1.7 and 2.2). Then we use Kowalski-Preiss' characterization for the support of uniform measures to determine what the support of the pseudo-tangent measure looks like. This, combined with the fact that the support is $\delta$-Reifenberg flat, allows us to conclude that all pseudo-tangent measures are multiples of Lebesgue measures on $n$-planes, and that the support of $\mu$ is Reifenberg flat with vanishing constant. Theorem 2.1 plays a key role in Sections 4, 5, and 6.

We first introduce some definitions and recall a classification result for uniform measures from [KP]. Let $\mu$ be a doubling Radon measure in $\mathbf{R}^{n+1}$, $r \in (0,1)$, and $Q \in \operatorname{spt}\mu$. Consider the new measure $\mu_{r,Q}$ defined by

$$\mu_{r,Q}(A) = \frac{\mu(rA+Q)}{\mu(B(r,Q))} \quad \text{for every Borel set } A \subset \mathbf{R}^{n+1}.$$

Note that $\mu_{r,Q}$ is a Radon measure. Since $\mu$ is a doubling measure an easy computation shows that for each compact set $K \subset \mathbf{R}^{n+1}$, and $Q \in \operatorname{spt}\mu \cap K$,

$$\sup_{0<r\leq 1} \mu_{r,Q}(K) \leq C,$$

where $C$ is a positive constant depending only on $n$, and $K$.

*Definition* 2.1. Let $\mu$ be a doubling Radon measure in $\mathbf{R}^{n+1}$. We say that $\nu$ is a *tangent measure* of $\mu$ at the point $Q \in \operatorname{spt}\mu$ if $\nu$ is a nonzero Radon measure in $\mathbf{R}^{n+1}$ and if there exists a sequence of positive numbers $\{r_i\}$ such that $r_i \downarrow 0$, and $\mu_{r_i,Q}$ converges weakly to $\nu$ in the sense of Radon measures, i.e. $\mu_{r_i,Q} \rightharpoonup \nu$. In particular for $\varphi \in C_c(\mathbf{R}^{n+1})$,

$$\lim_{i\to\infty} \frac{1}{\mu(B(r_i,Q))} \int \varphi(\frac{X-Q}{r_i}) d\,\mu(X) = \int \varphi d\,\nu.$$

For a detailed discussion of tangent measures for general Radon measures see [M] and [Pr]. Roughly speaking the tangent measures of a given doubling Radon measure $\mu$ provide pointwise information about the behavior of $\mu$. For our purposes we need some sort of uniform control on the behavior of $\mu$. Thus we introduce the notion of pseudo-tangent measure. This notion is analogous to the notion of pseudo-tangent map (resp. pseudo-tangent cone) introduced by Simon [Si1], [Si2] in order to study the regularity of the singular set of an energy minimizing harmonic map (resp. minimal surface).

Let $\mu$ be a doubling Radon measure in $\mathbf{R}^{n+1}$. A simple computation shows that for every compact set $K \subset \mathbf{R}^{n+1}$ there exists a constant $C(K) > 0$ depending only on $n$, and $K$ such that $Q \in \operatorname{spt}\mu \cap K$, $|Q - Q_i| \le 1$ and $r_i \le 1$,

$$\sup_{i\ge 0} \mu_{r_i,Q_i}(K) \le C(K).$$

*Definition* 2.2. Let $\mu$ be a doubling Radon measure in $\mathbf{R}^{n+1}$. We say that $\nu$ is a *pseudo-tangent measure* of $\mu$ at the point $Q \in \operatorname{spt}\mu$ if $\nu$ is a nonzero Radon measure in $\mathbf{R}^{n+1}$ and if there exists a sequence of points $Q_i \in \operatorname{spt}\mu$ such that $Q_i \to Q$, and a sequence of positive numbers $\{r_i\}$ such that $r_i \downarrow 0$, and $\mu_{r_i,Q_i} \rightharpoonup \nu$. In particular for $\varphi \in C_c(\mathbf{R}^{n+1})$,

$$\lim_{i\to\infty} \frac{1}{\mu(B(r_i,Q_i))} \int \varphi(\frac{X-Q_i}{r_i}) d\,\mu(X) = \int \varphi d\,\nu.$$

Note that a tangent measure of $\mu$ at the point $Q$ is a pseudo-tangent measure of $\mu$ at $Q$. We mentioned in the introduction that all the pseudo-tangent measures of an asymptotically optimally doubling measure are uniform. The following theorem gives a complete description of uniform measures.

THEOREM ([KP]). *Let $\nu$ be a nonzero Radon measure in $\mathbf{R}^{n+1}$ such that for every $Q \in \operatorname{spt}\nu$, and $r \in (0,\infty)$*

$$\nu(B(r,Q)) = \omega_n r^n,$$

*where $\omega_n$ denotes the volume of the unit ball in $\mathbf{R}^n$. Then after translation and rotation, either*

$$\nu = \mathcal{H}^n \llcorner \{(x_1,\ldots,x_{n+1}) \in \mathbf{R}^{n+1} : x_{n+1} = 0\},$$

*or*

$$n \ge 3 \quad and \quad \nu = \mathcal{H}^n \llcorner \{(x_1,\ldots,x_{n+1}) \in \mathbf{R}^{n+1} : x_4^2 = x_1^2 + x_2^2 + x_3^2\}.$$

THEOREM 2.1. *Let $\mu$ be a doubling Radon measure in $\mathbf{R}^{n+1}$ which is asymptotically optimally doubling. Then*:

1. *If $n = 1, 2$, $\operatorname{spt}\mu$ is Reifenberg flat with vanishing constant.*

2. *If $n \geq 3$ and $\operatorname{spt}\mu$ is $\delta$-Reifenberg flat, then $\operatorname{spt}\mu$ is Reifenberg flat with vanishing constant.*

THEOREM 2.2. *Let $\mu$ be a doubling Radon measure in $\mathbf{R}^{n+1}$ which is asymptotically optimally doubling. Then all pseudo-tangent measures of $\mu$ are uniform measures.*

*Remark* 2.1. The following fact is often used. If for each $i \geq 1$, $\Gamma_i \subset \mathbf{R}^{n+1}$ contains the origin then, given a compact set $K \subset \mathbf{R}^{n+1}$, there exists a subsequence such that $\Gamma_{i'} \cap K$ converges in the Hausdorff distance sense. Taking an exhaustion of $\mathbf{R}^{n+1}$ by compact sets, we can insure that there exists another subsequence $\{i_k\}$ such that $\Gamma_{i_k}$ converges to $\Gamma \ni 0$, in the Hausdorff distance sense, uniformly on compact sets. A similar result, known as Gromov's precompactness theorem [G, Prop. 5.2], holds in general metric spaces.

LEMMA 2.1. *Let $\mu$ be a doubling Radon measure in $\mathbf{R}^{n+1}$ which is asymptotically optimally doubling. Assume that there exist a sequence $Q_i \in \operatorname{spt}\mu$, so that $Q_i \to Q$, and a sequence $\lambda_i \downarrow 0$ such that $\mu_{\lambda_i,Q_i} \rightharpoonup \nu$. If $\Sigma = \operatorname{spt}\mu$ and $\eta_{\lambda_i,Q_i}(\Sigma) = \frac{1}{\lambda_i}(\Sigma - Q_i)$ then $X \in \operatorname{spt}\nu$ if and only if there exists a sequence $X_i \in \eta_{\lambda_i,Q_i}(\Sigma)$ such that $X_i \to X$.*

*Proof of Lemma* 2.1. Let $X_i \in \eta_{\lambda_i,Q_i}(\Sigma)$, and assume that $X_i \to X$. There exists a sequence $\{Z_i\}_i \subset \Sigma$ such that $X_i = \frac{1}{\lambda_i}(Z_i - Q_i)$. For $r \in (0,1)$, there exists $i_0 \geq 1$, such that for $i \geq i_0$, $|X - X_i| < \frac{r}{2}$, and $|Z_i - Q_i| = \lambda_i |X_i| \leq M\lambda_i$, where $M = |X| + 1$. Let $\tau = \frac{r}{2(M+1)}$. Since $\mu$ is asymptotically optimally doubling there exists $R > 0$ such that for $P \in B(1,Q) \cap \Sigma$, if $0 < \rho \leq R$ then

$$\tau^n \underset{1}{\sim} \frac{\mu(B(\tau\rho, P))}{\mu(B(\rho, P))}.$$

Since $Q_i \to Q$, $\lambda_i \downarrow 0$, and $|Z_i - Q_i| \leq M\lambda_i$, there exists $i_1 \geq i_0$ such that for $i \geq i_1$, $Z_i \in B(1,Q) \cap \Sigma$, $\frac{r}{2}\lambda_i \leq R$, and $(M+1)\lambda_i \leq R$. Under these assumptions we have that

$$\begin{aligned} \mu_{\lambda_i,Q_i}(B(r,X)) &= \frac{\mu(B(r\lambda_i, Q_i + \lambda_i X))}{\mu(B(\lambda_i, Q_i))} \\ &\geq \frac{\mu(B(r\lambda_i - \lambda_i|X - X_i|, Z_i))}{\mu(B(\lambda_i, Q_i))} \\ &\geq \frac{\mu(B(\frac{r}{2}\lambda_i, Z_i))}{\mu(B(\lambda_i, Q_i))} \end{aligned}$$

$$\geq \frac{\mu(B(\frac{r}{2}\lambda_i, Z_i))}{\mu(B(\lambda_i(M+1), Z_i))}$$
$$\geq \frac{1}{2}\Big(\frac{r}{2(M+1)}\Big)^n.$$

Recall that $\mu_{\lambda_i,Q_i}$ converges weakly to $\nu$. Therefore

$$\begin{aligned} \nu(B(2r,X)) &\geq \nu(\overline{B}(r,X)) \\ &\geq \limsup_{i\to\infty} \mu_{\lambda_i,Q_i}(\overline{B}(r,X)) \\ &\geq \frac{1}{2}\Big(\frac{r}{2(M+1)}\Big)^n. \end{aligned} \tag{2.1}$$

For $r > 0$, (2.1) guarantees that $\nu(B(r,X)) \geq C(n,|X|)r^n > 0$; thus $X \in \operatorname{spt}\nu$.

In order to prove that if $X \in \operatorname{spt}\nu$, there exists a sequence $X_i \in \eta_{\lambda_i,Q_i}(\Sigma)$ such that $X_i \to X$, assume not. Then $X \in \operatorname{spt}\nu$ and there exist $\varepsilon_0 > 0$ and a subsequence $\{i_k\}$ such that $d(X, \eta_{\lambda_{i_k},Q_{i_k}}(\Sigma)) \geq \varepsilon_0$. In particular $B(\frac{\varepsilon_0}{2}, X) \cap \eta_{\lambda_{i_k},Q_{i_k}}(\Sigma) = \emptyset$. When $\varphi \in C_c(B(\frac{\varepsilon_0}{2}, X))$,

$$\int \varphi d\nu = \lim_{i\to\infty} \frac{1}{\mu(B(\lambda_{i_k}, Q_{i_k}))} \int \varphi(\frac{Y - Q_{i_k}}{\lambda_{i_k}}) d\mu(Y).$$

If $Y \in \operatorname{spt}\mu = \Sigma$, then $|\frac{Y-Q_{i_k}}{\lambda_{i_k}} - X| \geq \frac{\varepsilon_0}{2}$, which implies that $\varphi(\frac{Y-Q_{i_k}}{\lambda_{i_k}}) = 0$. Thus $\int \varphi d\nu = 0$ for all $\varphi \in C_c(B(\frac{\varepsilon_0}{2}, X))$. This contradicts the fact that $X \in \operatorname{spt}\nu$. □

*Proof of Theorem* 2.2. Assume that $\nu$ is a pseudo-tangent measure of $\mu$ at the point $Q \in \operatorname{spt}\mu$. There exist a sequence $Q_i \in \operatorname{spt}\mu$, so that $Q_i \to Q$, and a sequence $\lambda_i \downarrow 0$, such that $\mu_{\lambda_i,Q_i} \rightharpoonup \nu$. When $X \in \operatorname{spt}\nu$, Lemma 2.1 guarantees that there exists a sequence $X_i \in \eta_{\lambda_i,Q_i}(\Sigma)$, such that $X_i \to X$. Moreover, there exists a sequence $\{Z_i\}_i \subset \Sigma$ such that $X_i = \frac{1}{\lambda_i}(Z_i - Q_i)$. Fix $r > 0$. Given $\varepsilon > 0$, there exists $i_0 \geq 1$, so that for $i \geq i_0$, $|X - X_i| < \min\{1, \varepsilon r\}$, and $|Z_i - Q_i| = \lambda_i |X_i| \leq M\lambda_i$, where $M = |X| + 1$. Under these assumptions,

$$\begin{aligned} \mu_{\lambda_i,Q_i}(B(r,X)) &\leq \frac{\mu(B(r\lambda_i + \lambda_i|X - X_i|, Z_i))}{\mu(B(\lambda_i, Q_i))} \\ &\leq \frac{\mu(B(r\lambda_i(1+\varepsilon), Z_i))}{\mu(B(\lambda_i, Q_i))} \\ &\leq \frac{\mu(B(r\lambda_i(1+\varepsilon), Z_i))}{\mu(B(r\lambda_i(1+\varepsilon), Q_i))} \cdot \frac{\mu(B(r\lambda_i(1+\varepsilon), Q_i))}{\mu(B(\lambda_i, Q_i))}. \end{aligned}$$

Let $\kappa$ be a large constant to be determined. Since $\mu$ is asymptotically optimally doubling, $Q_i \to Q$, $Z_i \to Q$, and $\lambda_i \downarrow 0$, there exists $i_1(\varepsilon, \kappa, M, r) \geq i_0$ so that for $i \geq i_1$,

$$\frac{\mu(B(r\lambda_i(1+\varepsilon), Q_i))}{\mu(B(\lambda_i, Q_i))} \leq (1+\varepsilon)[r(1+\varepsilon)]^n, \tag{2.2}$$

and

$$(2.3)\quad \begin{aligned} \frac{\mu(B(r\lambda_i(1+\varepsilon),Z_i))}{\mu(B(r\lambda_i(1+\varepsilon),Q_i))} &\leq (1+\varepsilon)^2\frac{\mu(B(\kappa r\lambda_i(1+\varepsilon),Z_i))}{\mu(B(\kappa r\lambda_i(1+\varepsilon),Q_i))} \\ &\leq (1+\varepsilon)^2\frac{\mu(B(\kappa\lambda_i r(1+\varepsilon+\frac{M}{r\kappa}),Q_i))}{\mu(B(\kappa\lambda_i r,Q_i))} \\ &\leq (1+\varepsilon)^3\Big(1+\varepsilon+\frac{M}{r\kappa}\Big)^n. \end{aligned}$$

Choosing $\kappa$ large enough so that $\frac{M}{\kappa r}<\varepsilon$, (2.2) and (2.3) yield that for $i\geq i_1$,

$$\mu_{\lambda_i,Q_i}(B(r,X))\leq(1+\varepsilon)^4(1+2\varepsilon)^n r^n.$$

Thus

$$\limsup_{i\to\infty}\mu_{\lambda_i,Q_i}(B(r,X))\leq r^n.$$

A similar argument shows that

$$\liminf_{i\to\infty}\mu_{\lambda_i,Q_i}(B(r,X))\geq r^n.$$

Therefore for $X\in\operatorname{spt}\nu$ and $r>0$,

$$\lim_{i\to\infty}\mu_{\lambda_i,Q_i}(B(r,X))=r^n.$$

Since $\mu_{\lambda_i,Q_i}$ converges weakly to $\nu$ in the sense of Radon measures, for $X\in\operatorname{spt}\nu$ and $r>0$,

$$\nu(B(r,X))\leq\liminf_{i\to\infty}\mu_{\lambda_i,Q_i}(B(r,X))=r^n,$$

and for any $\varepsilon>0$,

$$\begin{aligned}\nu(B(r,X))\geq\nu(\overline{B}(r(1-\varepsilon),X)) &\geq \limsup_{i\to\infty}\mu_{\lambda_i,Q_i}(\overline{B}(r(1-\varepsilon),X)) \\ &\geq \lim_{i\to\infty}\mu_{\lambda_i,Q_i}(B(r(1-\varepsilon),X))=(1-\varepsilon)^n r^n.\end{aligned}$$

We conclude that for $X\in\operatorname{spt}\nu$ and $r>0$, $\nu(B(r,X))=r^n$. □

*Proof of Theorem* 2.1. Let $\mu$ be a doubling Radon measure in $\mathbf{R}^{n+1}$ which is asymptotically optimally doubling. Let $\Sigma=\operatorname{spt}\mu$, and $K\subset\mathbf{R}^{n+1}$ be a compact set such that $K\cap\Sigma\neq\emptyset$. Let $\theta(r,Q)$ and $\theta_K(r)$ be as in (1.2) and (1.3) respectively, and let

$$\ell=\limsup_{r\downarrow 0}\theta_K(r).$$

We need to prove that $\ell=0$. By definition there exists a sequence $\lambda_i\downarrow 0$, such that $\theta_K(\lambda_i)\to\ell$. Also there exists a sequence of points $Q_i\in\Sigma\cap K$ such that $\theta(\lambda_i,Q_i)\to\ell$. Since $K$ is compact and $\Sigma$ is closed we may assume that $Q_i\to Q\in\Sigma\cap K$. Note that for every $i\geq 1$, $0\in\eta_{\lambda_i,Q_i}(\Sigma)$. Therefore there exists a subsequence $\{i_k\}$ such that $\eta_{\lambda_{i_k},Q_{i_k}}(\Sigma)$ converges to $\Sigma_\infty$ in the Hausdorff distance sense uniformly on compact sets (see Remark 2.1). Modulo passing to a further subsequence, which we relabel, Theorem 2.2 guarantees that

$$\eta_{\lambda_k,Q_k}(\Sigma) \to \Sigma_\infty \text{ in the Hausdorff distance sense, uniformly on compact sets,}$$

and

$$\mu_{\lambda_k,Q_k} \rightharpoonup \nu \text{ where } \nu \text{ is a uniform measure.}$$

Lemma 2.1 states that $\eta_{\lambda_k,Q_k}(\Sigma)$ converges to $\operatorname{spt}\nu$ in a pointwise sense. The fact that $\eta_{\lambda_k,Q_k}(\Sigma)$ converges to $\Sigma_\infty$ in the Hausdorff distance sense, uniformly on compact sets guarantees that

$$\Sigma_\infty = \operatorname{spt}\nu.$$

The Kowalski and Preiss theorem asserts that: i) if $n = 1, 2$, $\Sigma_\infty$ is an $n$-plane; ii) if $n \geq 3$, $\Sigma_\infty$ is either an $n$-plane or a cone. Our next step is to rule out the cone as a possibility for $\Sigma_\infty$ in the case where $n \geq 3$. It is here where the assumption about the Reifenberg flatness of $\Sigma = \operatorname{spt}\mu$ is used. In fact we prove that since $\Sigma$ is $\delta$-Reifenberg flat there exists an $n$-plane $L$ such that the Hausdorff distance between $L \cap B(1, 0)$ and $\Sigma_\infty \cap B(1, 0)$ is at most $\frac{1}{4\sqrt{2}}$. This condition is not satisfied by the Kowalski-Preiss cone.

When $X \in \Sigma_\infty$, there exists a sequence $\{Z_k\}_k \subset \Sigma \cap B(1, Q)$ such that $X_k = \frac{1}{\lambda_k}(Z_k - Q_k) \to X$ as $k \to \infty$. Without loss of generality we may assume that $|X - X_k| \leq 1/2$. Since $\Sigma$ is $\delta$-Reifenberg flat,

$$\sup_{0<r\leq R}\ \sup_{P\in\Sigma\cap B(1,Q)} \theta(r, P) \leq \delta \quad \text{for some } \delta \in (0, \frac{1}{4\sqrt{2}}) \text{ and } R > 0.$$

There exists $k_1 \geq 1$ such that for $k \geq k_1$, $s_k = \lambda_k(1 - |X - X_k|)$, $r_k = \lambda_k(1 + |X - X_k|) \leq R$. For each such $k$ there exist $n$-planes $L_k$ and $L'_k$ containing $Z_k$ and such that

$$D[\Sigma \cap B(r_k, Z_k), L_k \cap B(r_k, Z_k)] \leq \delta r_k, \tag{2.4}$$

and

$$D[\Sigma \cap B(s_k, Z_k), L'_k \cap B(s_k, Z_k)] \leq \delta s_k. \tag{2.5}$$

Note that for $Y \in L_k \cap B(s_k, Z_k)$, there exists $\overline{Y} \in L_k \cap B(s_k - \delta r_k, Z_k)$ such that $|Y - \overline{Y}| \leq \delta r_k$. (2.4) insures that there exists $P \in \Sigma \cap B(r_k, Z_k)$ such that $|\overline{Y} - P| \leq \delta r_k$. Note that $|P - Z_k| \leq s_k$ and $|P - Y| \leq 2\delta r_k$. Now (2.5) guarantees that there exists $Y' \in L'_k \cap B(s_k, Z_k)$ so that $|P - Y'| \leq \delta s_k$. Thus $|Y - Y'| \leq \delta s_k + 2\delta r_k$, and

$$L_k \cap B(s_k, Z_k) \subset \Big(L'_k \cap B(s_k, Z_k); \delta s_k + 2\delta r_k\Big).$$

Since $\eta_{\lambda_k,Q_k}(\Sigma) \to \Sigma_\infty$ in the Hausdorff distance sense, uniformly on compact sets, given $\varepsilon > 0$ there exists $k_0 \geq k_1$ such that for $k \geq k_0$

$$D[\Sigma_\infty \cap B(1,X), \eta_{\lambda_k,Q_k}(\Sigma) \cap B(1,X)] \leq \varepsilon. \tag{2.6}$$

Hence

$$\begin{aligned} \Sigma_\infty \cap B(1,X) &\subset \Big(\eta_{\lambda_k,Q_k}(\Sigma) \cap B(1,X); \varepsilon\Big) \\ &\subset \Big(\eta_{\lambda_k,Q_k}(\Sigma) \cap B(1+|X-X_k|, X_k); \varepsilon\Big), \end{aligned}$$

and

$$\begin{aligned} &\lambda_k\Big(\Sigma_\infty \cap B(1,X)\Big) + Q_k \\ &\quad\subset \Big(\Sigma \cap B(r_k, Z_k); \varepsilon\lambda_k\Big) \\ &\quad\subset \Big(L_k \cap B(r_k, Z_k); \varepsilon\lambda_k + \delta r_k\Big) \\ &\quad\subset \Big(L_k \cap B(s_k, Z_k); \varepsilon\lambda_k + \delta r_k + 2\lambda_k|X-X_k|\Big) \\ &\quad\subset \Big(L'_k \cap B(s_k, Z_k); \varepsilon\lambda_k + 3\delta r_k + \delta s_k + 2\lambda_k|X-X_k|\Big) \\ &\quad\subset \Big(L'_k \cap B(\lambda_k, Z_k); \varepsilon\lambda_k + 4\delta\lambda_k + 4\lambda_k|X-X_k|\Big). \end{aligned} \tag{2.7}$$

Moreover (2.6) also implies that

$$\begin{aligned} \eta_{\lambda_k,Q_k}(\Sigma) \cap B(1,X) &\subset \Big(\Sigma_\infty \cap B(1,X); \varepsilon\Big) \\ \eta_{\lambda_k,Q_k}(\Sigma) \cap B(1-|X-X_k|, X_k) &\subset \Big(\Sigma_\infty \cap B(1,X); \varepsilon\Big). \end{aligned} \tag{2.8}$$

Combining (2.5) and (2.8) we have that

$$\begin{aligned} L'_k \cap B(\lambda_k, Z_k) &\subset \Big(L'_k \cap B(s_k, Z_k); \lambda_k|X-X_k|\Big) \\ &\subset \Big(\Sigma \cap B(s_k, Z_k); \delta s_k + \lambda_k|X-X_k|\Big) \\ &\subset \Big(\lambda_k(\Sigma_\infty \cap B(1,X)) + Q_k; \varepsilon\lambda_k + \delta\lambda_k + \lambda_k|X-X_k|\Big). \end{aligned} \tag{2.9}$$

Therefore (2.7) and (2.9) yield

$$D[\Sigma_\infty \cap B(1,X), P_k \cap B(1,X_k)] \leq \varepsilon + 4\delta + 4|X-X_k|,$$

where $P_k = \frac{1}{\lambda_k}(L'_k - Q_k)$ is an $n$-plane containing $X_k$. Note that $\Lambda_k = P_k - X_k + X$ is an $n$-plane containing $X$, and satisfying

$$D[\Lambda_k \cap B(1,X), P_k \cap B(1,X_k)] \leq |X-X_k|,$$

and

$$D[\Sigma_\infty \cap B(1,X), \Lambda_k \cap B(1,X)] < \varepsilon + 4\delta + 5|X-X_k|.$$

Since the space of $n$-planes in $\mathbf{R}^{n+1}$ containing $X$ is compact, there exist an $n$-plane $\Lambda_X$ containing $X$ and a subsequence such that $\Lambda_{k'}$ converges to $\Lambda_X$ in the Hausdorff distance sense, uniformly on compact sets. Therefore for $X \in \Sigma_\infty$ there exists an $n$-plane $\Lambda_X$ containing $X$ such that

$$D[\Sigma_\infty \cap B(1,X), \Lambda_X \cap B(1,X)] \le 4\delta + \varepsilon.$$

By the Kowalski-Preiss theorem if $\Sigma_\infty$ is not an $n$-plane, $n \ge 3$ and there exist a translation and a rotation $\mathcal{R}$ such that

$$\mathcal{R}(\Sigma_\infty - Y_\infty) = \mathcal{C}, \quad \text{where } \mathcal{C} = \{(x_1, \dots, x_{n+1}) \in \mathbf{R}^{n+1} : x_4^2 = x_1^2 + x_2^2 + x_3^2\}.$$

Note that $Y_\infty \in \Sigma_\infty$, thus there exists an $n$-plane $\Lambda_\infty$ containing $Y_\infty$, so that if $\Lambda = \mathcal{R}(\Lambda_\infty - Y_\infty)$

$$D[\mathcal{C} \cap B(1,0), \Lambda \cap B(1,0)] \le 4\delta + \varepsilon.$$

Since $\delta \in (0, \frac{1}{4\sqrt{2}})$, $\varepsilon > 0$ can be chosen small enough so that

$$D[\mathcal{C} \cap B(1,0), \Lambda \cap B(1,0)] < \frac{1}{\sqrt{2}}.$$

On the other hand a simple computation shows that for any $n$-plane $L$ containing the origin

$$D[\mathcal{C} \cap B(1,0), L \cap B(1,0)] \ge \frac{1}{\sqrt{2}}.$$

Therefore for $n \ge 1$, $\Sigma_\infty$ is an $n$-plane. Now, (2.6) guarantees that given $\varepsilon > 0$ small enough there exists $k_0 \ge 1$ such that for $k \ge k_0$

$$D[\eta_{\lambda_k, Q_k}(\Sigma) \cap B(1,0), \Sigma_\infty \cap B(1,0)] \le \varepsilon.$$

Hence

$$\theta(\lambda_k, Q_k) \le \frac{1}{\lambda_k} D[\Sigma \cap B(\lambda_k, Q_k), \mathcal{L}_k \cap B(\lambda_k, Q_k)] \le \varepsilon,$$

where $\mathcal{L}_k = \lambda_k \Sigma_\infty + Q_k$ is an $n$-plane through $Q_k$. We conclude then that

$$\ell = \lim_{k \to \infty} \theta(\lambda_k, Q_k) = 0.$$

□

There is also a quantitative version of Theorem 2.1. In order to prove it we need to provide a quantitative version of the Kowalski and Preiss theorem. Namely:

Lemma 2.2. *Given $\varepsilon > 0$, there exists $\delta > 0$ such that, if $\nu$ is a nonzero Radon measure in $\mathbf{R}^{n+1}$ such that, for every $X \in \operatorname{spt} \nu$, and every $r \in (0, \infty)$*

$$\frac{1}{1+\delta} r^n \le \nu(B(r,X)) \le (1+\delta) r^n, \tag{2.10}$$

*then, if $\Sigma = \operatorname{spt} \nu$ and*

$$\theta_\Sigma(r) = \sup_{Q \in \Sigma} \inf_{L \ni Q} \left\{ \frac{1}{r} D[\Sigma \cap B(r,Q), L \cap B(r,Q)] \right\} :$$

1. *For* $n = 1, 2$, $\sup_{r>0} \theta_\Sigma(r) < \varepsilon$;

2. *For* $n \geq 3$, *if* $\sup_{r>0} \theta_\Sigma(r) < 1/4\sqrt{2}$ *then* $\sup_{r>0} \theta_\Sigma(r) < \varepsilon$.

*Proof.* Assume that the statement above is false. Then there exists $\varepsilon_0 > 0$ such that, for every $\delta > 0$, there is a Radon measure $\nu_\delta$ satisfying (2.10), but whose support $\Sigma_\delta = \operatorname{spt} \nu_\delta$ is such that $\sup_{r>0} \theta_\Sigma(r) \geq \varepsilon_0$. In the case $n \geq 3$, we also assume $\sup_{r>0} \theta_\Sigma(r) < 1/4\sqrt{2}$. Therefore, there exist sequences $\delta_i \downarrow 0$, $\lambda_i > 0$ and $Q_i \in \Sigma_i$ with $\Sigma_i = \Sigma_{\delta_i}$ such that

$$\theta_{\Sigma_i}(\lambda_i, Q_i) \geq \frac{\varepsilon_0}{2},$$

where, for $Q \in \Sigma_i$, $\theta_{\Sigma_i}(r, Q) = \inf_L \{\frac{1}{r} D[\Sigma_i \cap B(r, Q), L \cap B(r, Q)]\}$. Moreover, if

$$\nu_i(E) = \frac{\nu_{\delta_i}(\lambda_i E + Q_i)}{\lambda_i^n}, \quad \text{for each Borel set } E \subset \mathbf{R}^{n+1},$$

then, for $X \in \Gamma_i = \frac{1}{\lambda_i}(\Sigma_i - Q_i)$, and $r \in (0, \infty)$

$$\frac{1}{1 + \delta_i} r^n \leq \nu_i(B(r, X)) \leq (1 + \delta_i) r^n.$$

It is easy to see that for each $i \geq 1$

$$\theta_{\Gamma_i}(1, 0) \geq \frac{\varepsilon_0}{2}, \tag{2.11}$$

and that, if $n \geq 3$, $\sup_{r>0} \theta_{\Gamma_i}(r) < 1/4\sqrt{2}$. Similar arguments to the ones used to prove Lemma 2.1 and Theorem 2.2 insure that (modulo passing to a subsequence) $\nu_i \rightharpoonup \nu_\infty$ (weakly in the sense of Radon measures), $\Gamma_i \to \Gamma_\infty$ in the Hausdorff distance sense, uniformly on compact subsets, and $\Gamma_\infty = \operatorname{spt} \nu_\infty$. If $n \geq 3$, $\sup_{r>0} \theta_{\Gamma_\infty}(r) < 1/4\sqrt{2}$. Moreover, for $X \in \Gamma_\infty$ and $R > 0$, $\nu_\infty(B(r, X)) = r^n$. Therefore, $\Gamma_\infty$ is an $n$-plane by [KP]. The fact that $\Gamma_i \to \Gamma_\infty$ in the Hausdorff distance sense, uniformly on compact subsets, guarantees that, given $\varepsilon > 0$ there exists $i_0 \geq 1$ such that for $i \geq i_0$

$$\theta_{\Gamma_i}(1, 0) \leq D[\Gamma_i \cap B(1, 0), \Gamma_\infty \cap B(1, 0)] \leq \varepsilon.$$

This contradicts (2.11) whenever $\varepsilon < \varepsilon_0/2$. □

THEOREM 2.3. *Given* $\varepsilon > 0$, *there exists* $\delta > 0$ *such that, if* $\mu$ *is a doubling Radon measure in* $\mathbf{R}^{n+1}$, *and* $\Sigma = \operatorname{spt} \mu$ *then the following holds*:

1. *For* $n = 1, 2$, *if* $\mu$ *is* $\delta$*-approximately optimally doubling, then* $\Sigma$ *is* $\varepsilon$*-Reifenberg flat.*

2. *For* $n \geq 3$, *if* $\Sigma$ *is* $\eta$*-Reifenberg flat* (*for some* $\eta \in (0, 1/4\sqrt{2})$), *and* $\mu$ *is* $\delta$*-approximately optimally doubling then* $\Sigma$ *is* $\varepsilon$*-Reifenberg flat.*

*Proof.* Let $\delta > 0$, and $\mu$ be a doubling Radon measure in $\mathbf{R}^{n+1}$ which is $\delta$-approximately optimally doubling. If $n \geq 3$, the support of $\mu$, $\operatorname{spt}\mu = \Sigma$ is taken to be $\eta$-Reifenberg flat. Given a compact set $K \subset \mathbf{R}^{n+1}$, let $\ell = \limsup_{r\to 0}\theta_K(r)$. There exist sequences $\lambda_i \downarrow 0$ and $Q_i \in K \cap \Sigma$ with $Q_i \to Q \in K \cap \Sigma$ so that $\theta(\lambda_i, Q_i) \to \ell$. Let

$$\mu_i(E) = \frac{\mu(\lambda_i E + Q_i)}{\mu(B(\lambda_i, Q_i))}, \quad \text{for each Borel set } E \subset \mathbf{R}^{n+1},$$

and $\Gamma_i = \frac{1}{\lambda_i}(\Sigma - Q_i)$. Similar arguments to the ones used to prove Lemma 2.1 and Theorem 2.2 insure that (modulo passing to a subsequence), $\mu_i \rightharpoonup \nu$ (weakly in the sense of Radon measures), $\Gamma_i \to \Gamma \ni 0$ in the Hausdorff distance sense, uniformly on compact subsets, and $\Gamma = \operatorname{spt}\nu$. Note that $\nu$ is a pseudo-tangent measure of $\mu$ at $Q$. If $n \geq 3$, $\Gamma$ satisfies $\sup_{r>0}\theta_\Gamma(r) < 1/4\sqrt{2}$. Moreover, for $r > 0$ and $X \in \Gamma$

$$\frac{1}{1+\delta}r^n \leq \nu(B(r,X)) \leq (1+\delta)r^n.$$

Lemma 2.2 guarantees that, given $\varepsilon > 0$, there exists $\delta > 0$ so that $\sup_{r>0}\theta_\Gamma(r) < \varepsilon/2$. In particular, there exists an $n$-dimensional hyperplane containing 0 and satisfying

$$D[\Gamma \cap B(1,0), L \cap B(1,0)] \leq \frac{\varepsilon}{2}. \tag{2.12}$$

Since $\Gamma_i = \frac{1}{\lambda_i}(\Sigma - Q_i) \to \Gamma$ in the Hausdorff distance sense, uniformly on compact subsets, there exists $i_0 \geq 1$ such that, for $i \geq i_0$

$$D[\Gamma \cap B(1,0), \frac{1}{\lambda_i}(\Sigma - Q_i) \cap B(1,0)] \leq \frac{\varepsilon}{2}. \tag{2.13}$$

Thus combining (2.12) and (2.13), we conclude that, for $i \geq i_0$,

$$\theta(\lambda_i, Q_i) \leq D[L \cap B(1,0), \frac{1}{\lambda_i}(\Sigma - Q_i) \cap B(1,0)] \leq \varepsilon.$$

Thus $\ell = \limsup_{r\to 0}\theta_K(r) \leq \varepsilon$, which insures that $\operatorname{spt}\mu$ is $\varepsilon$-Reifenberg flat. □

Combining Reifenberg's Theorem with Theorem 2.3 we obtain the following regularity statement about the support of approximately optimally doubling Radon measures.

COROLLARY 2.1. *Let $\mu$ be a doubling Radon measure in $\mathbf{R}^{n+1}$, and $\Sigma = \operatorname{spt}\mu$. Then*

1. *If $n = 1, 2$, given $\beta \in (0,1)$, there exists $\delta(n,\beta) > 0$ such that, if $\mu$ is $\delta$-approximately optimally doubling, then $\Sigma$ is a $C^{0,\beta}$ $n$-dimensional submanifold.*

2. *If* $n \geq 3$ *and* $\Sigma$ *is* $\eta$*-Reifenberg flat, given* $\beta \in (0,1)$, *there exists* $\delta(n,\beta) > 0$ *such that, if* $\mu$ *is* $\delta$*-approximately optimally doubling, then* $\Sigma$ *is a* $C^{0,\beta}$ $n$*-dimensional submanifold.*

## 3. Doubling properties of the harmonic measure and geometry of the boundary

In the previous section we saw how the doubling properties of a measure determine the geometry of its support. In this section we focus on the relationship between the doubling properties of the harmonic measure of a domain $\Omega$ and the geometry of its boundary $\partial\Omega$. We start by recalling the definition of *non-tangentially accessible domains* (NTA) as well as the relevant theorems about the boundary behavior of the Green's function and the properties of the harmonic measure on these domains. An *M-non-tangential ball* $B(r,X)$, in a domain $\Omega$, is a ball in $\Omega$ whose distance to $\partial\Omega$ is comparable to its radius; i.e. $Mr > d(B(r,X),\partial\Omega) > M^{-1}r$. For $X_1,\ X_2 \in \Omega$ a *Harnack Chain* from $X_1$ to $X_2$ in $\Omega$ is a sequence of $M$-nontangential balls such that the first ball contains $X_1$, the last contains $X_2$, and such that consecutive balls intersect. The number of balls in this chain is called the length of the chain.

*Definition* 3.1 ([JK1]). A bounded (resp. unbounded) domain $\Omega$ in $\mathbf{R}^{n+1}$ is called *non-tangentially accessible* when there exist constants $M > 1$ and $R > 0$ (resp. $R = \infty$) such that:

1. *Corkscrew condition.* For any $Q \in \partial\Omega$, $r < R$ (resp. $r > 0$ ) there exists $A = A(r,Q) \in \Omega$ such that $M^{-1}r < |A - Q| < r$ and $d(A,\partial\Omega) > M^{-1}r$.

2. $\Omega^c$ satisfies the corkscrew condition.

3. *Harnack Chain Condition.* If $\varepsilon > 0$, and $X_1,\ X_2 \in \Omega \cap B(\frac{r}{4},Q)$ for some $Q \in \partial\Omega$, $r < R$ (resp $r > 0$), $d(X_j,\partial\Omega) > \varepsilon$ and $|X_1 - X_2| < 2^k\varepsilon$, then there exists a Harnack chain from $X_1$ to $X_2$ of length $Mk$ and such that the diameter of each ball is bounded below by $M^{-1}\min\{\mathrm{dist}(X_1,\partial\Omega), \mathrm{dist}(X_2,\partial\Omega)\}$.

When $\Omega$ is an unbounded domain, we also require as part of the definition of NTA domain that $\mathbf{R}^{n+1}\backslash\partial\Omega$ divide $\mathbf{R}^{n+1}$ into two distinct connected components $\Omega$ and $\operatorname{int}\Omega^c \neq \emptyset$.

Before stating further results we need to give a precise definition of harmonic measure. Let $\Omega \subset \mathbf{R}^{n+1}$ be a bounded domain, and $f$ be a function defined on $\partial\Omega$. Define the upper class of functions by

$$U_f = \{u : u \equiv +\infty \text{ on } \Omega \text{ or } \Delta u \leq 0 \text{ and } \liminf_{X\to Q} u(X) \geq f(Q)\ \forall Q \in \partial\Omega\},$$

and the lower class by $L_f = \{-u : u \in U_{-f}\}$. Let $\overline{H}f(X) = \inf\{u(X), u \in U_f\}$ (respectively $\underline{H}f(X) = \sup\{u(X), u \in L_f\}$) be the upper solution of the Dirichlet problem (resp. the lower solution). If $\overline{H}f(X) = \underline{H}f(X)$ for every $X \in \Omega$, and $\Delta(\overline{H}f) = 0$ in $\Omega$, $f$ is called a resolutive boundary function. In that case, we set $Hf(X) = \overline{H}f(X) = \underline{H}f(X)$. Weiner [W] showed that every continuous real-valued function on $\partial\Omega$ is resolutive. This fact, and the maximum principle make it possible to define harmonic measure.

*Definition* 3.2. The unique probability measure on $\partial\Omega$, denoted $\omega^X$, such that for all continuous functions $f$ on $\partial\Omega$, $Hf(X) = \int_{\partial\Omega} f\, d\omega^X$, is called the *harmonic measure* of $\Omega$ with pole at $X$.

If $\Omega \subset \mathbf{R}^{n+1}$ is an unbounded domain such that $\operatorname{int}\Omega^c \neq \emptyset$ the harmonic measure can be defined in a similar way (see [H, Ch. 9]). We note that as a consequence of Harnack's inequality, for any $X_1, X_2 \in \Omega$, and any domain $\Omega$, the measures $\omega^{X_1}$ and $\omega^{X_2}$ are mutually absolutely continuous. Thus, for a fixed point $X_0 \in \Omega$, the Radon-Nikodym derivative $K(X,Q) = (d\omega^X/d\omega^{X_0})(Q)$ exists, and is called the *kernel function*.

LEMMA 3.1 ([JK1], Lemma 4.9, 4.11). *Let $\Omega \subset \mathbf{R}^{n+1}$ be an* NTA *domain with constants $M > 1$ and $R > 0$, let $K \subset \mathbf{R}^{n+1}$ be a compact set, $Q \in \partial\Omega \cap K$, $0 < 2r < R$, and $X \in \Omega \backslash B(2Mr, Q)$. Then for $s \in [0, r]$,*

$$\omega^X(\Delta(2s,Q)) \leq C\omega^X(\Delta(s,Q)), \tag{3.1}$$

*where $C$ depends only on the* NTA *constants of $\Omega$ and on $K$. Here $\Delta(s,Q) = \partial\Omega \cap B(s,Q)$.*

In particular the harmonic measure on an NTA domain $\Omega$ is a doubling Radon measure on $\partial\Omega$. The pseudo-tangent measures of the harmonic measure deserve special attention as they are themselves harmonic measures, but with pole at infinity. In order to introduce the notion of harmonic measure with pole at infinity, for connected unbounded NTA domains, we need to recall some results concerning general NTA domains.

LEMMA 3.2 ([JK1], Lemma 4.1). *Let $\Omega$ be an* NTA *domain, let $K \subset \mathbf{R}^{n+1}$ be a compact set. There exists $\beta > 0$ such that for all $Q \in \partial\Omega \cap K$, $0 < 2r < R$, and every positive harmonic function $u$ in $\Omega \cap B(2r,Q)$, if $u$ vanishes continuously on $\Delta(2r,Q)$, then for $X \in \Omega \cap B(r,Q)$,*

$$u(X) \leq C\left(\frac{|X-Q|}{r}\right)^\beta \sup\{u(Y) : Y \in \partial B(2r,Q) \cap \Omega\},$$

*where $C$ depends only on $K$ and on the* NTA *constants.*

Corollary 3.1. *Let* $\Omega$ *be an* NTA *domain, let* $K \subset \mathbf{R}^{n+1}$ *be a compact set. Let* $Q \in \partial\Omega \cap K$, $0 < 2r < R$ *then*

$$\omega^{A(r,Q)}(\Delta(r,Q)) \geq C,$$

*where* $C$ *depends only on* $K$ *and on the* NTA *constants.*

Lemma 3.3 ([JK1], Lemma 4.4). *Let* $\Omega$ *be an* NTA *domain, let* $K \subset \mathbf{R}^{n+1}$ *be a compact set,* $Q \in \partial\Omega \cap K$, *and* $0 < 2r < R$. *If* $u \geq 0$ *is a harmonic function in* $\Omega$, *and* $u$ *vanishes continuously on* $\Delta(2r, Q)$, *then*

$$u(Y) \leq Cu(A(r,Q)),$$

*for all* $Y \in B(r,Q) \cap \Omega$. *Here* $C$ *depends only on* $K$ *and on the* NTA *constants.*

Lemma 3.4 ([JK1], Lemma 4.8). *Let* $\Omega$ *be an* NTA *domain, let* $K \subset \mathbf{R}^{n+1}$ *be a compact set,* $Q \in \partial\Omega \cap K$, $0 < 2r < R$, *and* $X \in \Omega \backslash B(2r, Q)$. *Then*

$$C^{-1} < \frac{\omega^X(\Delta(r,Q))}{r^{n-1} G(A(r,Q), X)} < C,$$

*where* $G(A(r,Q), -)$ *is the Green's function of* $\Omega$ *with pole at* $A(r,Q)$.

Lemma 3.5 ([JK1], Lemma 4.10, Comparison Principle). *Let* $\Omega$ *be an* NTA *domain,* $K \subset \mathbf{R}^{n+1}$ *be a compact set and* $0 < Mr < R$. *Suppose that* $u$ *and* $v$ *are positive harmonic functions in* $\Omega$ *vanishing continuously on* $\Delta(Mr, Q)$ *for some* $Q \in \partial\Omega \cap K$. *Then there exists a constant* $C > 1$ (*depending only on* $K$ *and on the* NTA *constants*) *such that for all* $X \in B(r,Q) \cap \Omega$,

$$C^{-1} \frac{u(A(r,Q))}{v(A(r,Q))} \leq \frac{u(X)}{v(X)} \leq C \frac{u(A(r,Q))}{v(A(r,Q))}.$$

Theorem 3.1 ([JK1], Theorem 7.1). *Let* $\Omega$ *be an* NTA *domain,* $K \subset \mathbf{R}^{n+1}$ *be a compact set. Let* $X_0, X \in \Omega$; *then for every* $Q \in \partial\Omega$

$$K(X,Q) = \frac{d\omega^X}{d\omega^{X_0}}(Q) = \lim_{r \to 0} \frac{\omega^X(\Delta(r,Q))}{\omega^{X_0}(\Delta(r,Q))} = \lim_{Z \to Q} \frac{G(X,Z)}{G(X_0,Z)}.$$

*There exist constants* $C > 1$, $N_0 > 1$ *and* $\alpha \in (0,1)$ *so that for* $s > 0$ *and* $Q_0 \in \partial\Omega$ *if* $X \in \Omega \backslash B(2Ns, Q_0)$, $N \geq N_0$, *then for every* $Q, Q' \in \Delta(s, Q_0)$

$$|K(X,Q') - K(X,Q)| \leq CK(X,Q) \left( \frac{|Q - Q'|}{s} \right)^{\alpha}.$$

Theorem 3.2 ([JK1], Theorem 7.9). *Let* $\Omega$ *be an* NTA *domain,* $K \subset \mathbf{R}^{n+1}$ *be a compact set. There exists a number* $\alpha > 0$, *such that for all* $Q \in \partial\Omega \cap K$, $0 < 2r < R$, *and all positive harmonic functions* $u$ *and* $v$ *in* $\Omega \cap B(2r, Q)$ *which vanish continuously on* $\Delta(2r, Q)$, *the function* $\frac{u(X)}{v(X)}$ *is*

*Hölder continuous of order* $\alpha$ *on* $\overline{\Omega} \cap \overline{B(r,Q)}$. *In particular, for every* $Q \in \partial\Omega$, $\lim_{X\to Q} \frac{u(X)}{v(X)}$ *exists, and for* $X, Y \in \Omega \cap B(r,Q)$,

$$\Big|\frac{u(X)}{v(X)} - \frac{u(Y)}{v(Y)}\Big| \leq C\frac{u(A(r,Q))}{v(A(r,Q))}\Big(\frac{|X-Y|}{r}\Big)^{\alpha}.$$

LEMMA 3.6 ([JK1], Lemma 4.11). *Let* $\Omega$ *be an* NTA *domain, and* $K \subset \mathbf{R}^{n+1}$ *be a compact set. Let* $Q_0, Q \in \partial\Omega \cap K$, $\Delta = \Delta(Q_0, r)$, *and* $r < R$. *Let* $\Delta' = \Delta(s,Q) \subset \Delta(\frac{r}{2}, Q_0)$. *If* $X \in \Omega\backslash B(Q_0, 2r)$, *then*

$$\omega^{A(r,Q_0)}(\Delta') \sim \frac{\omega^X(\Delta')}{\omega^X(\Delta)}.$$

($C_1 \sim C_2$ *means that the ratio between* $C_1$ *and* $C_2$ *is bounded above and below by a constant that depends only on* $M$, $R$ *and* $K$.)

The results quoted above from [JK1] were proved for bounded domains but they easily extend to unbounded domains.

LEMMA 3.7. *Let* $\Omega \subset \mathbf{R}^{n+1}$ *be an unbounded* NTA *domain, and* $Q \in \partial\Omega$. *There exists a unique function* $u$ *such that*

$$\begin{cases} \Delta u = 0 & \text{in } \Omega \\ u > 0 & \text{in } \Omega \\ u = 0 & \text{on } \partial\Omega, \end{cases} \tag{3.2}$$

*and*

$$u(A(1,Q)) = 1.$$

*Proof.* Note that without loss of generality we may assume that $Q = 0$. Let $A(1,0) = A$.

*Uniqueness.* Let $u$ and $v$ satisfy (3.2). By the comparison principle for $\rho > 1$ and $X \in B(\rho, 0) \cap \Omega$,

$$C^{-1}\frac{u(A(\rho,0))}{v(A(\rho,0))} \leq \frac{u(X)}{v(X)} \leq C\frac{u(A(\rho,0))}{v(A(\rho,0))},$$

where $C$ depends only on the NTA constants. Since $A \in B(\rho, 0)$, and $u(A) = v(A)$ then for $X \in B(\rho,0) \cap \Omega$,

$$C^{-1} \leq \frac{u(X)}{v(X)} \leq C. \tag{3.3}$$

Theorem 3.2 and (3.3) guarantee that for $X \in B(\rho,0) \cap \Omega$,

$$\Big|\frac{u(X)}{v(X)} - 1\Big| \leq C\frac{u(A(\rho,0))}{v(A(\rho,0))}\Big(\frac{|X-A|}{\rho}\Big)^{\alpha} \leq C\Big(\frac{|X-A|}{\rho}\Big)^{\alpha}.$$

Fixing $X$ and letting $\rho \to \infty$ we conclude that $u \equiv v$ in $\Omega$.

*Existence.* For $Y \in \Omega$ let $G(Y,-)$ denote the Green's function of $\Omega$ with pole at $Y$. Let

$$u_Y(X) = \frac{G(Y,X)}{G(Y,A)},$$

$u_Y$ is a nonnegative harmonic function on $B(|Y|,0) \cap \Omega$. Let $K \subset \mathbf{R}^{n+1}$ be a fixed compact set. Fix $\rho > 0$ such that $K \cap \Omega \subset B(\rho,0) \cap \Omega$, and let $|Y| \geq 2\rho$. Let $X \in K \cap \Omega$. Lemma 3.3 combined with the Harnack Principle yield

$$G(Y,X) \leq CG(Y,A(\rho,0)) \leq C_{K,n} G(Y,A). \tag{3.4}$$

Thus for $|Y| \geq 2\rho$, $\sup_{X \in K \cap \Omega} u_Y(X) \leq C_{K,n}$. Let $\{Y_j\}_j \subset \Omega$ be such that $|Y_j| \geq 2\rho$, and $|Y_j| \to \infty$ as $j \to \infty$. The corresponding $u_j = u_{Y_j}$ are nonnegative uniformly bounded harmonic functions on $B(\rho,0) \cap \Omega$. The Arzela-Ascoli theorem guarantees that there is a subsequence $u_{j'}$ which converges uniformly to a nonnegative harmonic function $u$ in $B(\rho,0) \cap \Omega$. Letting $\rho \to \infty$ and taking a diagonal subsequence we conclude that there is a subsequence $u_{j_k}$ which converges to the nonnegative harmonic function $u$, uniformly on compact sets of $\Omega$. Since $u(A) = 1$ and $u \equiv 0$ on $\partial\Omega$, $u > 0$ in $\Omega$. Therefore $u$ satisfies (3.2). By the uniqueness proved above we conclude that $u_{|Y|} \to u$ as $|Y| \to \infty$. □

COROLLARY 3.2. *Let $\Omega \subset \mathbf{R}^{n+1}$ be an unbounded* NTA *domain. For $Q \in \partial\Omega$, let $\Delta(1,Q) = \partial\Omega \cap B(1,Q)$. There exists a unique doubling Radon measure $\omega^\infty$, supported on $\partial\Omega$ satisfying*:

$$\int_{\partial\Omega} \varphi d\omega^\infty = \int_\Omega v\Delta\varphi \quad \textit{for all} \quad \varphi \in C_c^\infty(\mathbf{R}^{n+1})$$

*where*

$$\begin{cases} \Delta v = 0 & \text{in } \Omega \\ v > 0 & \text{in } \Omega \\ v = 0 & \text{on } \partial\Omega, \end{cases}$$

*and*

$$\omega^\infty(\Delta(1,Q)) = 1.$$

*$\omega^\infty$ is called the harmonic measure of $\Omega$ with pole at infinity normalized at $Q \in \partial\Omega$.*

*Proof.* Again without loss of generality we may assume that $Q = 0$ and $A(1,0) = A$.

*Uniqueness.* Let $\omega_1$, $\omega_2$ be two such measures. Let $v_1$, $v_2$ be the corresponding nonnegative harmonic functions in $\Omega$ which vanish on $\partial\Omega$. By Lemma 3.7 $v_i(X) = v_i(A)u(X)$ where $u$ is the unique function satisfying (3.2). Thus for $\varphi \in C_c^\infty(\mathbf{R}^{n+1})$,

$$\int_{\partial\Omega} \varphi d\omega_i = v_i(A) \int_\Omega u\Delta\varphi,$$

which implies that

$$\frac{1}{v_1(A)}\int_{\partial\Omega}\varphi d\omega_1 = \frac{1}{v_2(A)}\int_{\partial\Omega}\varphi d\omega_2. \tag{3.5}$$

(3.5) asserts that

$$\frac{d\omega_1}{v_1(A)} = \frac{d\omega_2}{v_2(A)} \quad \Rightarrow \frac{\omega_1(\Delta(1,0))}{v_1(A)} = \frac{\omega_2(\Delta(1,0))}{v_2(A)} \Rightarrow v_1(A) = v_2(A).$$

Thus $v_1 \equiv v_2$ and $\omega_1 = \omega_2$.

*Existence.* Fix $\rho > 0$, and let $\varphi \in C_c^\infty(B(\rho, 0))$, let $Y \in \Omega$ be such that $|Y| \geq 2\rho$. Let $\omega^Y$ denote the harmonic measure of $\Omega$ with pole at $Y$. Then

$$\int_{\partial\Omega}\varphi(Q)\frac{d\omega^Y(Q)}{G(Y,A)} = \int_\Omega \Delta\varphi(X)\frac{G(Y,X)}{G(Y,A)}dX = \int_\Omega \Delta\varphi(X)u_Y(X)dX.$$

Lemma 3.4 guarantees that

$$\frac{\omega^Y(\Delta(\rho,0))}{G(Y,A)} \sim \rho^{n-1}\frac{G(Y,A(\rho,0))}{G(Y,A)} = \rho^{n-1}u_Y(A(\rho,0)).$$

Therefore the Radon measures $\frac{\omega^Y}{G(Y,A)}$ are uniformly bounded on $B(\rho, 0)$ (see (3.4)). Given a sequence $\{Y_j\}_j \subset \Omega$ such that $|Y_j| \geq 2\rho$, and $|Y_j| \to \infty$ as $j \to \infty$, there exists a subsequence $\{Y_{j'}\}$ and a Radon measure $\mu$ such that for $\phi \in C_c^\infty(B(\rho, 0))$

$$\int \phi\frac{d\omega^{Y_{j'}}}{G(Y_{j'},A)} \to \int \phi d\mu.$$

Letting $\rho \to \infty$, and taking a diagonal subsequence $\{j_k\}$ we conclude that

$$\frac{\omega^{Y_{j_k}}}{G(Y_{j_k},A)} \rightharpoonup \mu.$$

Since $u_{Y_{j_k}}$ converges uniformly to $u$ (as in the proof of Lemma 3.7), in compact subsets of $\Omega$ we also have that

$$\int_{\partial\Omega}\varphi d\mu = \int_\Omega u\Delta\varphi dX \quad \text{ for all } \varphi \in C_c^\infty(\mathbf{R}^{n+1}).$$

If $\omega^\infty = \frac{\mu}{\mu(\Delta(1,0))}$ then $\omega^\infty(\Delta(1,0)) = 1$, and

$$\int_{\partial\Omega}\varphi d\omega^\infty = \int_\Omega v\Delta\varphi \quad \text{ for all } \varphi \in C_c^\infty(\mathbf{R}^{n+1}).$$

Here $v = \frac{u}{\mu(\Delta(1,0))}$ is a nonnegative harmonic function in $\Omega$ which vanishes in $\partial\Omega$. In order to prove that $\omega^\infty$ is a doubling measure let $K \subset \mathbf{R}^{n+1}$ be a compact set and let $Q \in K \cap \partial\Omega$. Given $r > 0$, there exists $j_K \geq 1$ such that if $j_k \geq j_K$, $Y_{j_k} \in \Omega\backslash B(2Mr, Q)$. Then for $s \in [0, r]$ (3.1) guarantees that

$$\omega^{Y_{j_k}}(\Delta(2s,Q)) \leq C\omega^{Y_{j_k}}(\Delta(s,Q)),$$

where $C$ depends only on $K$ and the NTA constants of $\Omega$ (see Lemma 3.1). Hence

$$\begin{aligned}\omega^{\infty}(\Delta(2s,Q)) &\leq \liminf_{j_k\to\infty} \frac{\omega^{Y_{j_k}}(\Delta(2s,Q))}{\mu(\Delta(1,0))G(Y_{j_k},A)} \\ &\leq C\liminf_{j_k\to\infty} \frac{\omega^{Y_{j_k}}(\overline{\Delta}(\frac{s}{2},Q))}{\mu(\Delta(1,0))G(Y_{j_k},A)} \\ &\leq C\omega^{\infty}(\overline{\Delta}(\frac{s}{2},Q)) \\ &\leq C\omega^{\infty}(\Delta(s,Q)).\end{aligned}$$ □

The next lemma states that the class of NTA domains in $\mathbf{R}^{n+1}$ is compact under blow-ups, and that in some sense so is the class of harmonic measures. Namely:

LEMMA 3.8. *Let $\Omega \subset \mathbf{R}^{n+1}$ be an* NTA *domain. Let $\{\lambda_i\}$ be a sequence of positive numbers such that $\lambda_i \downarrow 0$, let $Q_i \in \partial\Omega$ be such that $Q_i \to Q \in \partial\Omega$, and let $Y \in \Omega$. For $E$ a Borel set in $\mathbf{R}^{n+1}$ let $\omega^Y(E) = \omega^Y(E \cap \partial\Omega)$, and $\mu_j(E) = \frac{\omega^Y(\lambda_j E + Q_j)}{\omega^Y(B(\lambda_j,Q_j))}$. There exist subsequences $\lambda_j \to 0$ and $Q_j \to Q$ such that*

$$\eta_{\lambda_j,Q_j}(\Omega) \to \Omega_\infty \qquad \text{in the Hausdorff distance sense, uniformly on compact sets,}$$

*and*

$$\eta_{\lambda_j,Q_j}(\partial\Omega) \to \partial\Omega_\infty \qquad \text{in the Hausdorff distance sense, uniformly on compact sets,}$$

*where $\Omega_\infty$ is an unbounded* NTA *domain. Moreover*

$$\mu_j \rightharpoonup \nu,$$

*where $\nu$, a pseudo-tangent measure of $\omega^Y$ at $Q \in \partial\Omega$, is a constant multiple of the harmonic measure of $\Omega_\infty$ with pole at infinity normalized at $0 \in \partial\Omega_\infty$.*

The proof of this lemma is straightforward although slightly technical. Since we do not need the result in the sequel, we only present a rough outline of the proof. Note that since $\Omega$ is an NTA domain its harmonic measure with pole at $Y$ is doubling. This combined with Remark 2.1 guarantees that there exists a subsequence (which we relabel) such that $\eta_{\lambda_j,Q_j}(\Omega) \to \Omega_\infty$, $\eta_{\lambda_j,Q_j}(\partial\Omega) \to \Sigma_\infty$, and $\mu_j \rightharpoonup \nu$. An argument similar to the one presented in Lemma 2.1 guarantees that $\Sigma_\infty = \operatorname{spt}\nu$. In order to show that $\Sigma_\infty = \partial\Omega_\infty$, we combine a connectivity argument with the fact that both $\Omega$ and $\Omega^c$ satisfy the corkscrew condition. If $\Omega$ is an NTA domain with constants $M$ and $R > 0$, by inspection, one can check that $\eta_{\lambda_j,Q_j}(\Omega)$ is an NTA domain with constants $M$ and $R/\lambda_j$. It is straightforward then that $\Omega_\infty$ is an unbounded NTA domain.

The fact that $\nu$ is a multiple of $\omega^\infty$ is a consequence of Corollary 3.2 and Lemma 3.4 applied to the harmonic measure of $\eta_{\lambda_j,Q_j}(\Omega)$ with pole at $(Y-Q_j)/\lambda_j$.

An immediate consequence of Theorem 2.1, Theorem 2.3 and Theorem 4.3 in [KT1] is the characterization of Reifenberg flat domains in terms of the doubling properties of their harmonic measure.

THEOREM 3.3.

1. *Let $n = 1, 2$, and let $\Omega \subset \mathbf{R}^{n+1}$ be an* NTA *domain. Given $\varepsilon > 0$ there exists $\delta > 0$ such that if $\omega$ is $\delta$-approximately optimally doubling, then $\Omega$ is an $\varepsilon$-Reifenberg flat domain.*

2. *Let $n \geq 3$ and let $\Omega \subset \mathbf{R}^{n+1}$ be an $\eta$-Reifenberg flat domain. Given $\varepsilon > 0$ there exists $\delta > 0$ such that if $\omega$ is $\delta$-approximately optimally doubling, then $\Omega$ is an $\varepsilon$-Reifenberg flat domain.*

*Here $\omega$ denotes either the harmonic measure with pole $X \in \Omega$ or the harmonic measure with pole at infinity in the case where $\Omega$ is an unbounded domain.*

*Definition* 3.3. Let $\Omega \subset \mathbf{R}^{n+1}$. We say that $\Omega$ is a *Reifenberg flat domain with vanishing constant* if it is a $\delta$-Reifenberg flat domain for some $\delta \in (0, \delta_n)$, and $\partial\Omega$ is a Reifenberg flat set with vanishing constant.

THEOREM 3.4.

1. *Let $n = 1, 2$, and let $\Omega \subset \mathbf{R}^{n+1}$ be an* NTA *domain. $\Omega$ is a Reifenberg flat domain with vanishing constant if and only if $\omega$ is asymptotically optimally doubling.*

2. *Let $n \geq 3$ and let $\Omega \subset \mathbf{R}^{n+1}$ be an $\eta$-Reifenberg flat domain. $\Omega$ is a Reifenberg flat domain with vanishing constant if and only if $\omega$ is asymptotically optimally doubling.*

In the last part of this section we show that the Kowalski-Preiss cone gives rise to an example of an NTA domain $\Omega$ which is not $\delta$-Reifenberg flat but for which the harmonic measure with pole at infinity (normalized appropriately) and the surface measure of the boundary coincide. Moreover the surface measure of the boundary, and therefore its harmonic measure, are uniform.

PROPOSITION 3.1. *Let $n \geq 3$. Let*

$$\Omega = \left\{ (x_1, \dots, x_{n+1}) \in \mathbf{R}^{n+1} : |x_4| < \sqrt{x_1^2 + x_2^2 + x_3^2} \right\}.$$

*$\Omega$ is an unbounded* NTA *domain whose harmonic measure $\omega^\infty$ with pole at infinity and normalized at the origin satisfies*

$$\omega^\infty = \frac{1}{\omega_n} \mathcal{H}^n \llcorner \partial\Omega.$$

*Moreover for $Q \in \partial\Omega$ and $r > 0$*

$$\omega^\infty(B(r,Q)) = \frac{1}{\omega_n}\mathcal{H}^n(\partial\Omega \cap B(r,Q)) = r^n.$$

*Remark* 3.1. By uniformly smoothing out the corners of the domain $B(R,0) \cap \Omega$ at $\partial B(R,0) \cap \Omega$ we produce a bounded NTA domain whose harmonic measure at each point is asymptotically optimally doubling, but whose boundary is not Reifenberg flat with vanishing constant.

*Proof of Proposition* 3.1. The fact that $\Omega$ is an NTA domain is straightforward, thus the proof is omitted. Since

$$\partial\Omega = \left\{(x_1, \dots, x_{n+1}) \subset \mathbf{R}^{n+1} : x_4^2 = x_1^2 + x_2^2 + x_3^2\right\},$$

the Kowalski-Preiss theorem guarantees that for $Q \in \partial\Omega$ and $r > 0$,

$$\mathcal{H}^n \llcorner \partial\Omega(B(r,Q)) = \omega_n r^n.$$

We introduce a new set of coordinates. Let $r = \sqrt{x_1^2 + x_2^2 + x_3^2 + x_4^2}$. For $\theta \in [\frac{\pi}{4}, \frac{3\pi}{4}]$, let $x_4 = r\cos\theta$, then $X = (x_1, x_2, x_3, x_4, \dots, x_{n+1}) \in \overline{\Omega}$. Let

$$u(X) = -\frac{r}{2\sqrt{2}\omega_n}\frac{\cos 2\theta}{\sin\theta}.$$

Note that $u = 0$ if $\theta = \frac{\pi}{4}$ or $\theta = \frac{3\pi}{4}$, i.e. $u = 0$ on $\partial\Omega$. Moreover if $X \in \Omega$ then $u(X) > 0$. Computing the Laplacian of $u$ in terms of the coordinates $r$ and $\theta$ one shows that $u$ is a harmonic function in $\Omega$. Let $\overrightarrow{n}$ denote the inward unit normal vector to $\partial\Omega$. We have that either

$$\frac{\partial u}{\partial n} = \langle \overrightarrow{n}, \nabla u\rangle = \frac{1}{r}\frac{\partial u}{\partial\theta}\Big|_{\theta=\frac{\pi}{4}} = \frac{1}{\omega_n},$$

or

$$\frac{\partial u}{\partial n} = -\langle \overrightarrow{n}, \nabla u\rangle = \frac{1}{r}\frac{\partial u}{\partial\theta}\Big|_{\theta=\frac{3\pi}{4}} = \frac{1}{\omega_n}.$$

For $\varphi \in C_c^\infty(\mathbf{R}^{n+1})$,

$$\int_\Omega u\Delta\varphi = \int_{\partial\Omega}\frac{\partial u}{\partial n}\varphi d\mathcal{H}^n = \frac{1}{\omega_n}\int_{\partial\Omega}\varphi d\mathcal{H}^n.$$

Corollary 3.2 guarantees that $\omega^\infty = \frac{1}{\omega_n}\mathcal{H}^n \llcorner \partial\Omega$. □

## 4. Some characterizations of chord-arc domains with small constant

The notion of chord-arc surface with small constant (CASSC) was first introduced by S. Semmes (see [Se1, Se2]). A CASSC is a smooth hypersurface which divides $\mathbf{R}^{n+1}$ into two distinct connected components and whose unit normal vector has small mean oscillation (i.e. the normal vector has small BMO norm). They were characterized in terms of the behavior of the corresponding Clifford-Cauchy integrals. Motivated by Semmes' work we introduced the notion of chord arc domains (see [KT1]). In this section we provide several new measure theoretic characterizations of chord arc domains with small (and vanishing) constant (see Definition 4.1 and Theorem 4.6). These results play an important role in the proof of the Main Theorem (Section 5), which combined with the results in [KT1] (Section 5 also) provide a characterization of chord arc domains with vanishing constant via potential theory. One of the technical ingredients needed in Sections 5 and 6 is Lemma 4.1. This lemma asserts that the surface balls on the boundary of a chord arc domain which is Reifenberg flat with vanishing constant can be approximated by a the graph of a Lipschitz function with small constant and with very small $C^0$ norm.

Recall the definition of chord arc domain given in the introduction.

*Definition* 1.11. Let $\delta \in (0, \delta_n)$. A set of locally finite perimeter $\Omega$ (see [EG]) is said to be a *$\delta$-chord arc domain* or a *chord arc domain with small constant* if $\Omega$ is a $\delta$-Reifenberg flat domain, $\partial\Omega$ is Ahlfors regular and for each compact set $K \subset \mathbf{R}^{n+1}$ there exists $R > 0$ so that

$$\sup_{Q \in \partial\Omega \cap K} \|\overrightarrow{n}\|_*(B(R, Q)) < \delta.$$

Here $\overrightarrow{n}$ denotes the unit normal vector to the boundary.

*Definition* 4.1. Let $\Omega \subset \mathbf{R}^{n+1}$ be a domain. $\Omega$ is said to be a *chord arc domain with vanishing constant* if $\Omega$ is a chord arc domain with small constant and $\overrightarrow{n} \in \mathrm{VMO}(\partial\Omega)$.

*Remark* 4.1. Assume that for $Q \in \partial\Omega$ and $r > 0$, there are an $n$-dimensional plane $L(r, Q)$ containing $Q$ and a choice of unit normal vector to $L(r, Q)$, $\overrightarrow{n}_{r,Q}$ such that

$$\frac{1}{r} D[\partial\Omega \cap B(r, Q), L(r, Q) \cap B(r, Q)] < \delta,$$

for some $\delta \in (0, 1/4\sqrt{2})$, and

$$\mathcal{T}^+(r, Q) = \{X = (x, t) = x + t\overrightarrow{n}_{r,Q} \in B(r, Q) : x \in L(r, Q),\ t > 2\delta r\} \subset \Omega,$$

and

$$\mathcal{T}^-(r, Q) = \{X = (x, t) = x + t\overrightarrow{n}_{r,Q} \in B(r, Q) : x \in L(r, Q),\ t < -2\delta r\} \subset \Omega^c.$$

For $s \in [\frac{r}{2}, r]$, let

$$\begin{aligned}\Gamma(s, Q) = \partial\Omega \cap \{X &= (x, t) \\ &= x + t\overrightarrow{n}_{r,Q} \in \mathbf{R}^{n+1} : x \in L(r, Q), |x - Q| \le s,\ |t| \le s\}.\end{aligned}$$

Note that if $\Delta(r,Q) = B(r,Q) \cap \partial\Omega$,

$$\Gamma(r\sqrt{1-4\delta^2}, Q) \subset \Delta(r,Q) \subset \Gamma(r,Q).$$

In fact if $X = (x,t) \in \Gamma(r\sqrt{1-4\delta^2}, Q)$, $|x| \le r\sqrt{1-4\delta^2}$ and $|t| < 2\delta r$. Let $\Pi : \mathbf{R}^{n+1} \to L(r,Q)$ denote the orthogonal projection onto $L(r,Q)$. The separation property, and a connectivity argument guarantee that for $\tau \in [1/2, \sqrt{1-4\delta^2}]$,

$$\Pi(\Gamma(\tau r, Q)) = \{x \in L(r,Q),\ |x-Q| \le \tau r\}.$$

Thus for $\delta > 0$ small if $\sigma = \mathcal{H}^n \llcorner \partial\Omega$,

$$\begin{aligned} \sigma(\Delta(r,Q)) &\ge \sigma(\Gamma(r\sqrt{1-4\delta^2}, Q)) \\ &\ge |\Pi(\Gamma(r\sqrt{1-4\delta^2}, Q))| = \omega_n r^n (1-4\delta^2)^{\frac{n}{2}} \\ &\ge (1+\delta)^{-1}\omega_n r^n. \end{aligned}$$

In particular, if $\Omega$ is a Reifenberg flat domain, given a compact set $K \subset \mathbf{R}^{n+1}$ there exists an $R > 0$, such that for every $Q \in \partial\Omega \cap K$ and $r \in (0,R]$, there exists $n$-dimensional plane $L(r,Q)$ containing $Q$ and a choice of a unit normal vector to $L(r,Q)$ satisfying the hypothesis above. Thus the conclusion also holds. Namely if $\Omega$ is a $\delta$-Reifenberg flat domain for $\delta$ small enough, then given any compact set $K \subset \mathbf{R}^{n+1}$ there is $R > 0$ so that for every $Q \in \partial\Omega \cap K$, and $r \in (0,R]$, $\mathcal{H}^n(B(r,Q) \cap \partial\Omega) \ge (1+\delta)^{-1}\omega_n r^n$.

*Remark* 4.2. Note that if $\Omega$ is a set of locally finite perimeter which is a Reifenberg flat domain then the topological boundary of $\Omega$ and its measure theoretic boundary agree. In fact if $\Omega$ is Reifenberg flat for each compact set $K \subset \mathbf{R}^{n+1}$ there exists $R > 0$ such that for all $Q \in \partial\Omega \cap K$ and $r \in (0,R]$ there exists an $n$-dimensional plane $\mathcal{L}(r,Q)$ containing $Q$ and a choice of unit normal vector to $\mathcal{L}(r,Q)$, $\overrightarrow{n}_{r,Q}$ satisfying

$$\mathcal{T}^+(r,Q) = \{X = (x,t) = x + t\,\overrightarrow{n}_{r,Q} \in B(r,Q) : x \in \mathcal{L}(r,Q),\ t > \frac{1}{4}r\} \subset \Omega,$$

and

$$\mathcal{T}^-(r,Q) = \{X = (x,t) = x + t\,\overrightarrow{n}_{r,Q} \in B(r,Q) : x \in \mathcal{L}(r,Q),\ t < -\frac{1}{4}r\} \subset \Omega^c.$$

Therefore for all $Q \in \partial\Omega$,

$$\limsup_{r\to 0} \frac{\mathcal{H}^{n+1}(B(r,Q)\cap\Omega)}{\omega_{n+1}r^{n+1}} \ge \left(\frac{1}{4}\right)^{n+1},$$

and

$$\limsup_{r\to 0} \frac{\mathcal{H}^{n+1}(B(r,Q)\cap\Omega^c)}{\omega_{n+1}r^{n+1}} \ge \left(\frac{1}{4}\right)^{n+1},$$

where $\omega_{n+1}$ denotes the volume of the unit ball in $\mathbf{R}^{n+1}$. Therefore for $\mathcal{H}^n$ a.e. $Q \in \partial\Omega$ there exists a unique measure theoretic inner unit normal $\overrightarrow{n}(Q)$ such that the Gauss-Green theorem holds,

$$\int_\Omega \operatorname{div}\varphi\, dx = -\int_{\partial\Omega} \varphi \cdot \overrightarrow{n}(Q)\, d\sigma(Q),$$

for all $\varphi \in C_c^1(\mathbf{R}^{n+1}, \mathbf{R}^{n+1})$. For further discussion of sets of locally finite perimeter, and the measure theoretic unit normal see [EG, Ch. 5]. Note that if $\Omega$ is a set of locally finite perimeter and if $\Omega$ and $\Omega^c$ satisfy the corkscrew condition a similar argument proves that the topological boundary of $\Omega$ and its measure theoretic boundary agree.

In his study of chord arc surfaces with small constant Semmes proved a very useful decomposition lemma. The result also holds in the more general setting.

*Definition* 4.2. Let $\Omega \subset \mathbf{R}^{n+1}$ be a set of locally finite perimeter which is a Reifenberg flat domain. Let $\delta > 0$. We say that $\Omega$ is *$\delta$-Semmes decomposable*, if for each compact set $K \subset \mathbf{R}^{n+1}$ there exists $R > 0$ such that for $P \in \partial\Omega \cap K$ and $0 < r \leq R$ there exists a Lipschitz function $h : L(r, P) = \langle \overrightarrow{n}(r, Q)\rangle^{\perp} \to \mathbf{R}$, $\|\nabla h\|_\infty \leq \delta$, whose graph $\mathcal{G} = \{X = (x, t) \in \mathbf{R}^{n+1} : \ t = h(x)\}$ approximates $\partial\Omega$ in the cylinder $C(r, P) = \{(x, t) : x \in L(r, P) \cap B(r, P), |t| \leq r\}$ in the sense that

$$\sigma(C(r, P) \cap (\{\partial\Omega\backslash\mathcal{G}\} \cup \{\mathcal{G}\backslash\partial\Omega\})) \leq C_1 \exp(-\frac{C_2}{\delta})\omega_n r^n,$$

for some $C_1, C_2 > 0$. Moreover, $C(r, P) \cap \partial\Omega = G \cup B$, where $\sigma(B) \leq C_1 \exp(-\frac{C_2}{\delta})\omega_n r^n$, $G \subset \mathcal{G}$, and where $Q \in B$ implies

$$|Q - (\Pi(Q), h(\Pi(Q)))| \leq \delta \operatorname{dist}(\Pi(Q), \Pi(G)).$$

Also $\Pi(C(r, P) \cap \partial\Omega) = \{|x - P| \leq r\}$.

Note that if $\Omega \subset \mathbf{R}^{n+1}$ is $\delta$-Semmes decomposable, then $\partial\Omega$ is a $2\delta$-Reifenberg flat set. Thus when talking about $\delta$-Semmes decomposable domains we always assume that $\delta \in (0, \delta_n/2)$, where $\delta_n$ is the constant that appears in Definition 1.10.

Theorem 4.1 ([Se1], Semmes Decomposition). *There exists $\delta(n) > 0$ so that for $\delta \in (0, \delta(n))$, and $\Omega$ a set of locally finite perimeter, if $\Omega$ is a $\delta$-Reifenberg flat domain and if for each compact set $K \subset \mathbf{R}^{n+1}$, there exist $R > 0$ and $C > 1$ such that*

$$\sigma(\Delta(r, Q)) \leq C r^n \quad \text{for } Q \in K \cap \partial\Omega, \quad \text{and} \quad r \in (0, R], \tag{4.1}$$

*and*

$$\sup_{Q \in \partial\Omega \cap K} \|\overrightarrow{n}\|_*(B(R, Q)) \leq \delta^2,$$

*then $\Omega$ is $\delta$-Semmes decomposable.*

Note that, in particular, by Remark 4.1 this implies that for $r \in (0, R]$, and $P \in K \cap \partial\Omega$

$$\Delta(r, P) = G(r, P) \cup E(r, P),$$

where $G(r,P) \subset \mathcal{G}$ and $\sigma(E(r,P)) \le C\exp(-C_2/\delta)\sigma(\Delta(r,P))$. Here $\mathcal{G}$ is the graph of the Lipschitz function $h$ that appears in the Semmes decomposition.

Condition (4.1) and Remark 4.1 guarantee that under the hypothesis of Theorem 4.1, $\partial\Omega$ is a space of homogeneous type in the sense of Coifman and Weiss [CW]. Therefore the $L^p$ estimates for the Hardy-Littlewood maximal function hold on $\partial\Omega$. Using the same Calderón-Zygmund decomposition for the unit normal $\overrightarrow{n}$ (which exists almost everywhere in $\partial\Omega$), and the same degree theory argument as in the smooth case, we conclude that the Semmes decomposition is valid in this context. It is clear then that:

Theorem 4.2. *There exists $\delta(n) > 0$ so that for $\delta \in (0,\delta(n))$, and $\Omega$ a set of locally finite perimeter, if $\Omega$ is a $\delta$-Reifenberg flat domain, and if for each compact set $K \subset \mathbf{R}^{n+1}$ there exist $R > 0$ and $C > 1$ such that*

$$\sigma(\Delta(r,Q)) \le Cr^n \quad \text{for } Q \in K \cap \partial\Omega, \quad \text{and} \quad r \in (0,R],$$

*and*

$$\sup_{Q\in\partial\Omega\cap K} \|\overrightarrow{n}\|_*(B(R,Q)) \le \delta^2,$$

*then for every $Q \in \partial\Omega \cap K$ and every $r \in (0,R]$,*

$$(1+\delta)^{-1}\omega_n r^n \le \sigma(\Delta(r,Q)) \le (1+\delta)\omega_n r^n.$$

Theorem 4.3 ([KT1], Theorem 2.1). *There exists $\delta(n) > 0$ so that for $\delta \in (0,\delta(n))$ and $\Omega$ a set of locally finite perimeter, if $\Omega$ is a $\delta$-Reifenberg flat domain, and for each compact set $K \subset \mathbf{R}^{n+1}$ there exists $R > 0$ such that*

$$\sigma(\Delta(r,Q)) \le (1+\delta)\omega_n r^n \quad \forall Q \in \partial\Omega \cap K \quad \text{and} \quad \forall r \in (0,R],$$

*then*

$$\sup_{Q\in K\cap\partial\Omega} \|\overrightarrow{n}\|_*(B(\frac{R}{2},Q)) \le C\sqrt{\delta}, \tag{4.2}$$

*i.e. $\Omega$ is a chord arc domain with small constant. Here $C$ is a constant that depends only on the dimension.*

Combining the results obtained in Section 2, and in this section, we provide a new characterization for chord arc domains in terms of surface measure and Reifenberg flatness. Namely:

Theorem 4.4. *Let $n = 1,2$. Given $\delta > 0$, there exists $\eta > 0$, such that if $\Omega$ is a set of locally finite perimeter whose topological boundary agrees with its measure theoretic boundary, and if for each compact set $K \subset \mathbf{R}^{n+1}$ there exists $R > 0$ such that*

$$(1+\eta)^{-1}\omega_n r^n \le \sigma(\Delta(r,Q)) \le (1+\eta)\omega_n r^n \quad \forall Q \in \partial\Omega \cap K \quad \text{and} \quad \forall r \in (0,R], \tag{4.3}$$

*then $\Omega$ is a $\delta$-chord arc domain.*

*Remark* 4.3. In particular the above theorem asserts that an Ahlfors regular domain (see [D]) with constant close to 1 in $\mathbf{R}^2$ or in $\mathbf{R}^3$ is a chord arc domain. On the other hand a bounded domain with a cusp, smooth away from the cusp, is an Ahlfors regular domain, which is not a chord arc.

THEOREM 4.5. *Let $n \geq 3$. Given $\delta > 0$, there exists $\eta > 0$, such that if $\Omega$ is a set of locally finite perimeter which is a Reifenberg flat domain, and if for each compact set $K \subset \mathbf{R}^{n+1}$ there exists $R > 0$ such that*

(4.4)
$$(1+\eta)^{-1}\omega_n r^n \leq \sigma(\Delta(r,Q)) \leq (1+\eta)\omega_n r^n \quad \forall Q \in \partial\Omega \cap K \quad \text{and} \quad \forall r \in (0,R],$$

*then $\Omega$ is a $\delta$-chord arc domain.*

*Proof.* Note that in both cases the topological boundary of $\Omega$ agrees with its measure theoretic boundary, i.e. $\sigma = \mathcal{H}^n \llcorner \partial\Omega$. Moreover conditions (4.3) and (4.4) guarantee that as long as $\eta$ is small enough $\sigma$ is a $3\eta$-approximately optimally doubling measure measure. Theorem 2.3 asserts that there exists $\eta(\delta, n) \in (0, \delta)$ so that if $\eta \leq \eta(\delta, n)$ then $\partial\Omega$ is $\delta$-Reifenberg flat. Moreover, when $C\sqrt{\eta} \leq \delta$ where $C$ is as in (4.2), Theorem 4.3 insures that for each compact set $K \subset \mathbf{R}^{n+1}$ there exists $R_0 > 0$ so that $\sup_{Q\in\partial\Omega\cap K} \|\overrightarrow{n}\|(B(R_0, Q)) \leq \delta$. Thus the only remaining thing to check is that $\Omega$ satisfies the separation property (see Definition 1.9), in the case where $n = 1, 2$. We proceed by contradiction. If the separation property did not hold, a connectivity argument, plus the fact that $\partial\Omega$ is $\delta$-Reifenberg flat would allow us to show that there exist $s > 0$ and $Q \in \partial\Omega$ such that for every $P \in \Delta(s, Q)$ either

$$\limsup_{r\to 0} \frac{\mathcal{H}^{n+1}(B(r,P)\cap\Omega)}{\omega_{n+1} r^{n+1}} = 0$$

or

$$\limsup_{r\to 0} \frac{\mathcal{H}^{n+1}(B(r,P)\cap\Omega^c)}{\omega_{n+1} r^{n+1}} = 0.$$

This contradicts the assumption that the topological boundary of $\Omega$ agrees with its measure theoretic boundary. □

THEOREM 4.6. *Let $\Omega \subset \mathbf{R}^{n+1}$ be a set of locally finite perimeter which is a $\delta$-Reifenberg flat domain. The following statements are equivalent*:

1. *$\Omega$ is a chord arc domain with vanishing constant.*

2. *$\Omega$ is a Reifenberg flat domain with vanishing constant, and*

$$\lim_{r\to 0} \sup_{Q\in\partial\Omega\cap K} \frac{\sigma(\Delta(r,Q))}{\omega_n r^n} = 1.$$

3. *For each compact set* $K \subset \mathbf{R}^{n+1}$

$$\lim_{r\to 0} \inf_{Q\in\partial\Omega\cap K} \frac{\sigma(\Delta(r,Q))}{\omega_n r^n} = \lim_{r\to 0} \sup_{Q\in\partial\Omega\cap K} \frac{\sigma(\Delta(r,Q))}{\omega_n r^n} = 1.$$

We finish this section by giving a version of the Semmes decomposition lemma for chord arc domains which are Reifenberg flat with vanishing constant. This geometric lemma plays a key role in Sections 5 and 6.

LEMMA 4.1. *There exists* $\delta_0(n) > 0$ *such that, if* $\Omega$ *is a* $\delta$*-Semmes decomposable domain with* $\delta \le \delta_0$, *and for each compact set* $K \subset \mathbf{R}^{n+1}$ *there exists* $r_0 > 0$ *such that, for* $Q_0 \in \partial\Omega \cap K$ *and* $r \in (0, r_0)$ *there is an n-plane* $L(r, Q_0)$ *containing* $Q_0$, *and satisfying*

$$\frac{1}{r} D[\partial\Omega \cap B(r,Q_0), L(r,Q_0) \cap B(r,Q_0)] < \varepsilon,$$

*where* $\varepsilon > 0$ *is small (in particular* $\varepsilon < \delta$*), then for* $r \in (0, r_0)$ *there exists a Lipschitz function* $\phi : L(r,Q_0) \to \mathbf{R}$ *such that* $\|\nabla\phi\|_\infty \le C(n)(\delta + \varepsilon)$, *whose graph* $\mathcal{G} = \{(x,t) \in \mathbf{R}^{n+1} : t = \phi(x)\}$ *approximates* $\partial\Omega$ *in the ball* $B(r, Q_0)$ *in the sense that*

$$\Delta(r,Q_0) = \mathcal{G}(r,Q_0) \cup \mathcal{E}(r,Q_0),$$

*where*

$$\mathcal{G}(r,Q_0) \subset \mathcal{G} \quad \text{and} \quad \sigma(\mathcal{E}(r,Q_0)) \le C_1 \exp(-\frac{C_2}{\delta})\sigma(\Delta(r,Q_0)),$$

*for some* $C_1, C_2 > 0$. *Moreover if* $\Pi : \mathbf{R}^{n+1} \to L(r,Q_0)$ *denotes the orthogonal projection then*

$$\sup_{\Pi(\mathcal{G}(r,Q_0))} |\phi| \le \varepsilon r, \quad \text{and} \quad \sup_{B(r,Q_0)\cap L(r,Q_0)} |\phi| \le C(n)(\varepsilon + \delta) r.$$

In particular, Lemma 4.1 asserts that if $\Omega$ is a chord arc domain which is Reifenberg flat with vanishing constant, Semmes decomposition can be improved in the sense that we get a better estimate for the supremum of the function whose graph approximates $\partial\Omega$. This apparently minor detail plays an important role in the proof of the Main Lemma in Section 5.

*Proof.* Let $\varepsilon' = \varepsilon'(\varepsilon) \in (0, \frac{1}{16})$. Let $\delta \in (0, \delta_n)$ and assume that $\Omega$ is a $\delta$-Semmes decomposable domain. In particular there exists $R > 0$ such that for $r \in (0, R)$ and $Q_0 \in K \cap \partial\Omega$ there exist an $n$-dimensional plane $\Lambda(r, Q_0)$ and a Lipschitz function $h : \Lambda(r, Q_0) \to \mathbf{R}$ such that $\|\nabla h\|_\infty \le \delta$ and whose graph $\mathcal{G} = \{(x,t) \in \mathbf{R}^{n+1} : t = h(x)\}$ approximates $\partial\Omega$ in the ball $B(r, Q_0)$ in the sense that

$$\Delta(r,Q_0) = \mathcal{G}(r,Q_0) \cup \mathcal{E}(r,Q_0),$$

where

$$\mathcal{G}(r,Q_0)\subset\mathcal{G}\ \text{ and }\ \sigma(\mathcal{E}(r,Q_0))\le C_1\exp(-\frac{C_2}{\delta})\sigma(\Delta(r,Q_0)),$$

for some $C_1, C_2 > 0$. Recall that being $\delta$-Semmes decomposable also guarantees that if $Q \in C(r,Q_0)\cap\partial\Omega$ where $C(r,Q_0)=\{(x,t): x\in\Lambda(r,Q_0)\cap B(r,Q_0), |t|\le r\}$ then

$$|Q-(\Pi_0(Q),h(\Pi_0(Q)))|\le\delta\,\mathrm{dist}(\Pi_0(Q),\Pi_0(\mathcal{G}(r,Q_0))).$$

Here $\Pi_0$ denotes the orthogonal projection of $\mathbf{R}^{n+1}$ onto $\Lambda(r,Q_0)$. In particular $|h(\Pi_0(Q_0))|=|h(Q_0)|\le\delta r$ which implies that for $Q\in C(r,Q_0)\cap\partial\Omega$, $|h(\Pi_0(Q))|\le\delta|\Pi_0(Q)-Q_0|+\delta r\le 2\delta r$. Hence

$$\partial\Omega\cap B(r,Q_0)\subset\partial\Omega\cap C(r,Q_0)\subset(\Lambda(r,Q_0)\cap B(r,Q_0),2\delta r)\,.$$

Moreover $\Pi_0(C(r,Q_0)\cap\partial\Omega)=\Lambda(r,Q_0)\cap B(r,Q_0)$. For $X\in\Lambda(r,Q_0)\cap B(r,Q_0)$ there exists $Y\in\Lambda(r,Q_0)\cap B((1-4\delta)r,Q_0)$ such that $|X-Y|\le 4\delta r$. There exists $Q\in C(r,Q_0)\cap\partial\Omega$ such that $\Pi_0(Q)=Y$. Thus $|Y-Q|=|\Pi_0(Q)-Q|\le|(\Pi_0(Q),h(\Pi_0(Q)))-Q|+|h(\Pi_0(Q))|\le\delta r+2\delta r\le 3\delta r$, and $|Q-Q_0|\le|\Pi_0(Q)-Q_0|+3\delta r<r$. Hence $|X-Q|\le|Y-Q|+|X-Y|\le 7\delta r$. We conclude that

$$\Lambda(r,Q_0)\cap B(r,Q_0)\subset(\partial\Omega\cap B(r,Q_0),7\delta r)\,,$$

and

$$\frac{1}{r}D[\partial\Omega\cap B(r,Q_0);\Lambda(r,Q_0)\cap B(r,Q_0)]\le 7\delta. \tag{4.5}$$

Our hypotheses guarantee that for every $r\in(0,r_0)$ there exists an $n$-dimensional plane $L(r,Q_0)$ such that

$$\frac{1}{r}D[\partial\Omega\cap B(r,Q_0);L(r,Q_0)\cap B(r,Q_0)]\le\varepsilon. \tag{4.6}$$

Combining (4.5) and (4.6) we obtain that the angle between $\Lambda(r,Q_0)$ and $L(r,Q_0)$ is small. More specifically

$$d(\Lambda(r,Q_0),L(r,Q_0))=D[\Lambda(r,Q_0)\cap B(1,Q_0);L(r,Q_0)\cap B(1,Q_0)]\le 7\delta+\varepsilon.$$

In fact if $X\in\Lambda(r,Q_0)\cap B(1,Q_0)$, (4.5) guarantees that there exists $Q\in\partial\Omega\cap B(r,Q_0)$ such that $|r(X-Q_0)+Q_0-Q|\le 7\delta r$ and (4.6) insures that there exists $Y\in L(r,Q_0)\cap B(1,Q_0)$ such that $|r(Y-Q_0)+Q_0-Q|\le\varepsilon r$. Therefore $|X-Y|\le 7\delta+\varepsilon$. This proves that

$$\Lambda(r,Q_0)\cap B(1,Q_0)\subset(L(r,Q_0)\cap B(1,Q_0),7\delta+\varepsilon)\,.$$

The other inclusion is proved in exactly the same manner. Our goal is to show that $\Delta(r,Q_0)$ can be well approximated by the graph of a Lipschitz function $\phi: L(r,Q_0)\to\mathbf{R}$. Let $\overrightarrow{\nu}_0$ be the unit normal vector of $\Lambda(r,Q_0)$ and $\overrightarrow{n}_0$ be the unit normal vector of $L(r, Q_0)$. Provided that we choose the appropriate orientation, $\langle \overrightarrow{\nu}_0, \overrightarrow{n}_0 \rangle = \sqrt{1-(7\delta+\varepsilon)^2} \geq 1/2$ for $\delta \in (0, 1/32)$ and $\varepsilon \in (0, 1/16)$. Let $\Pi : \mathbf{R}^{n+1} \to L(r, Q_0)$ denote the orthogonal projection. We claim that $\Pi$ is one-to-one on $\mathcal{G}(r, Q_0)$. In order to simplify the computations we may assume, without loss of generality, that $Q_0 = 0$. If $P, Q \in \mathcal{G}(r, 0)$ then $P = p + h(p)\overrightarrow{\nu}_0$, and $Q = q + h(q)\overrightarrow{\nu}_0$, where $p, q \in \Lambda(r, 0) \cap B(r, 0)$. If $\Pi(Q) = \Pi(P)$ then

$$\begin{aligned}
\Pi(P) &= P - \langle P, \overrightarrow{n}_0\rangle(\overrightarrow{n}_0 - \langle \overrightarrow{n}_0, \overrightarrow{\nu}_0\rangle \overrightarrow{\nu}_0) - \langle P, \overrightarrow{n}_0\rangle\langle \overrightarrow{n}_0, \overrightarrow{\nu}_0\rangle \overrightarrow{\nu}_0,\\
\Pi(Q) &= Q - \langle Q, \overrightarrow{n}_0\rangle(\overrightarrow{n}_0 - \langle \overrightarrow{n}_0, \overrightarrow{\nu}_0\rangle \overrightarrow{\nu}_0) - \langle Q, \overrightarrow{n}_0\rangle\langle \overrightarrow{n}_0, \overrightarrow{\nu}_0\rangle \overrightarrow{\nu}_0.
\end{aligned}$$

Hence,

$$\begin{aligned}
(4.7)\qquad & p - q - \langle P-Q, \overrightarrow{n}_0\rangle(\overrightarrow{n}_0 - \langle \overrightarrow{n}_0, \overrightarrow{\nu}_0\rangle \overrightarrow{\nu}_0)\\
&\qquad = \langle P-Q, \overrightarrow{n}_0\rangle\langle \overrightarrow{n}_0, \overrightarrow{\nu}_0\rangle \overrightarrow{\nu}_0 + (h(q)-h(p))\overrightarrow{\nu}_0.
\end{aligned}$$

Note that

$$(4.8)\qquad p - q - \langle P-Q, \overrightarrow{n}_0\rangle(\overrightarrow{n}_0 - \langle \overrightarrow{n}_0, \overrightarrow{\nu}_0\rangle \overrightarrow{\nu}_0) \in \Lambda(r, Q_0) = \langle \overrightarrow{\nu}_0\rangle^{\perp}.$$

Combining (4.7) and (4.8) we have that

$$\langle P-Q, \overrightarrow{n}_0\rangle\langle \overrightarrow{n}_0, \overrightarrow{\nu}_0\rangle = h(p) - h(q),$$

which implies

$$(4.9)\qquad \langle p-q, \overrightarrow{n}_0\rangle = \frac{(h(p)-h(q))}{\langle \overrightarrow{n}_0, \overrightarrow{\nu}_0\rangle}(1 - \langle \overrightarrow{n}_0, \overrightarrow{\nu}_0\rangle^2).$$

Combining (4.7), (4.8) and (4.9) we have

$$\begin{aligned}
p - q &= \langle p-q, \overrightarrow{n}_0\rangle(\overrightarrow{n}_0 - \langle \overrightarrow{n}_0, \overrightarrow{\nu}_0\rangle \overrightarrow{\nu}_0)\\
&\quad + (h(p)-h(q))\langle \overrightarrow{\nu}_0, \overrightarrow{n}_0\rangle(\overrightarrow{n}_0 - \langle \overrightarrow{n}_0, \overrightarrow{\nu}_0\rangle \overrightarrow{\nu}_0),
\end{aligned}$$

and

$$\begin{aligned}
|p-q| &\leq |\langle p-q, \overrightarrow{n}_0\rangle| + |h(p)-h(q)|,\\
|p-q| &\leq \left(\frac{1}{|\langle \overrightarrow{n}_0, \overrightarrow{\nu}_0\rangle|} + 1\right)|h(p)-h(q)|,\\
|p-q| &\leq 3|h(p)-h(q)| \leq 3\delta|p-q| \leq \frac{3}{32}|p-q|,
\end{aligned}$$

which implies that $p = q$, and therefore $P = Q$. Hence for $\widetilde{p} \in \Pi(\mathcal{G}(r, Q_0)) \subset L(r, Q_0)$ there exists a unique $P \in \mathcal{G}(r, Q_0)$ such that $\Pi(P) = \widetilde{p}$. We can define $\varphi : \Pi(\mathcal{G}(r, Q_0)) \to \mathbf{R}$ by $\varphi(\widetilde{p}) = \langle P - \widetilde{p}, \overrightarrow{n}_0\rangle$. Let $\widetilde{p}, \widetilde{q} \in \Pi(\mathcal{G}(r, Q_0))$; by the definition of $\varphi$ we know that there exist $P, Q \in \mathcal{G}(r, Q_0)$ such that

$$P = \widetilde{p} + \varphi(\widetilde{p})\overrightarrow{n}_0 = p + h(p)\overrightarrow{\nu}_0 \quad \text{and} \quad Q = \widetilde{q} + \varphi(\widetilde{q})\overrightarrow{n}_0 = q + h(q)\overrightarrow{\nu}_0,$$

where $p, q \in \Lambda(r, Q_0) \cap B(r, Q_0)$. In order to prove that $\varphi$ is Lipschitz in $\Pi(\mathcal{G}(r, Q_0))$ we compute

$$\begin{aligned}
|\varphi(\widetilde{p}) - \varphi(\widetilde{q})| &\leq |\langle p - q, \overrightarrow{n}_0\rangle| + |h(p) - h(q)||\langle \overrightarrow{n}_0, \overrightarrow{\nu}_0\rangle \\
&\leq |\langle p - q, \overrightarrow{n}_0 - \overrightarrow{\nu}_0\rangle| + |h(p) - h(q)| \\
&\leq |p - q|(\delta + |\overrightarrow{n}_0 - \overrightarrow{\nu}_0|) \\
&\leq |p - q| \left(\delta + \sqrt{2(1 - \langle \overrightarrow{\nu}_0, \overrightarrow{n}_0\rangle)}\right) \\
&\leq (15\delta + 2\varepsilon)|p - q| \\
&\leq (15\delta + 2\varepsilon)|P - Q| \\
&\leq (15\delta + 2\varepsilon)\sqrt{|\widetilde{p} - \widetilde{q}|^2 + |\varphi(\widetilde{p}) - \varphi(\widetilde{q})|^2} \\
&\leq (15\delta + 2\varepsilon)\left(|\widetilde{p} - \widetilde{q}| + |\varphi(\widetilde{p}) - \varphi(\widetilde{q})|\right).
\end{aligned}$$

Thus

$$|\varphi(\widetilde{p}) - \varphi(\widetilde{q})| \leq 4(15\delta + 2\varepsilon)|\widetilde{p} - \widetilde{q}|.$$

The extension theorem for Lipschitz functions guarantees that there exists a Lipschitz function $\phi : L(r, Q_0) \to \mathbf{R}$ such that $\varphi = \phi_{|\Pi(\mathcal{G}(r,Q_0))}$ and so that $\|\nabla\varphi\|_\infty = \|\nabla\phi\|_\infty \leq 4(15\delta + 2\varepsilon)$ (see [EG, Ch. 3, Th. 1]). If $\widetilde{p} \in \Pi(\mathcal{G}(r, Q_0))$, $P = \widetilde{p} + \varphi(\widetilde{p})\overrightarrow{n}_0 \in \mathcal{G}(r, Q_0) \subset (L(r, Q_0), \varepsilon r)$ which implies that $|\varphi(\widetilde{p})| \leq \varepsilon r$. □

## 5. Poisson kernel estimates on chord arc domains

In Sections 5 and 6 we provide a characterization of chord arc domains with small constant in terms of the doubling properties of their harmonic measure and the oscillation of the logarithm of their Poisson kernel. We prove that if the harmonic measure of a chord arc domain with small constant is asymptotically optimally doubling and the logarithm of its Poisson kernel has vanishing mean oscillation then the domain is a chord arc domain with vanishing constant. In [KT1] we established the converse of this result, i.e. that on a chord arc domain with vanishing constant the harmonic measure is asymptotically optimally doubling and the logarithm of the Poisson kernel has vanishing mean oscillation. In Section 5 we establish several estimates that are satisfied by the Poisson kernel of a chord arc domain with small constant $\Omega$, whenever its logarithm has vanishing mean oscillation, and the corresponding harmonic measure is asymptotically optimally doubling. Some of these estimates allow us to compare the Poisson kernel of $\Omega$ to the Poisson kernel of an appropriate half space (see Corollaries 5.12 and 5.13). In Section 6 we exploit this relationship by means of Rellich's identity (see [JK3]) applied to both $\Omega$ and the half space mentioned above. The Main Theorem should be understood as a result about the regularity of a free boundary under weak conditions. From information about the measure whose density is the normal derivative of the Green's function (i.e the Poisson kernel), and this normal derivative on the boundary of the domain, we obtain results about the regularity of the boundary.

Initially we recall some results proved in [KT1]. If $\Omega$ is a $\delta$-chord arc domain for $\delta$ small enough, we know that the surface measure and the harmonic measure are mutually absolutely continuous (see [DJ], [KT1], [Se3]). For $X \in \Omega$, $k_X = \frac{d\omega^X}{d\sigma}$ denotes the Poisson kernel. The statements appear here in their scale invariant form. Since the proofs are identical to the ones presented in [KT1], we omit them.

*Definition* 5.1. Let $\Omega \subset \mathbf{R}^{n+1}$ be an unbounded domain. Let $\delta \in (0, \delta_n)$. We say that $\Omega$ is a $(\delta, \infty)$-*chord arc domain* if $\Omega$ is a set of locally finite perimeter such that

$$\sup_{r>0} \sup_{Q \in \partial\Omega} \theta(r, Q) < \delta, \tag{5.1}$$

and

$$\sigma(\Delta(r, Q)) \le (1 + \delta)\omega_n r^n \quad \text{for } Q \in \partial\Omega \text{ and } r > 0. \tag{5.2}$$

Here $\sigma = \mathcal{H}^n \llcorner \partial\Omega$, where $\mathcal{H}^n$ denotes the $n$-dimensional Hausdorff measure, and $\omega_n$ denotes the volume of the $n$-dimensional unit ball in $\mathbf{R}^n$.

Note that this definition could also have been stated in terms of the properties of the unit normal vector to the boundary (see Theorems 4.4 and 4.5).

THEOREM 5.1 ([KT1], Corollary 5.2). *Given $\varepsilon > 0$, there exist $\delta_1(\varepsilon) > 0$ and $N(\varepsilon) > 1$ such that, if $\Omega$ is a $(\delta_1, \infty)$-chord arc domain, for $N \ge N(\varepsilon)$, $Q \in \partial\Omega$, $s > 0$ and $X \in \Omega \backslash B(Ns, Q)$,*

$$\frac{1}{\sigma(\Delta(s, Q))} \int_{\Delta(s,Q)} |\log k_X - (\log k_X)_{s,Q}| d\sigma \le \varepsilon.$$

*Here* $(\log k_X)_{s,Q} = \frac{1}{\sigma(\Delta(s,Q))} \int_{\Delta(s,Q)} \log k_X d\sigma$.

THEOREM 5.2 ([KT1], Corollary 5.2]). *Given $\varepsilon > 0$, there exists $\delta_1(\varepsilon) > 0$ so that, if $\Omega$ is a bounded $\delta_1$-chord arc domain, there exist $N(\varepsilon) > 1$ and $s(\varepsilon, \operatorname{diam}\Omega) > 0$ such that for $N \ge N(\varepsilon)$, $s \in (0, s(\varepsilon, \operatorname{diam}\Omega))$, $Q \in \partial\Omega$, and $X \in \Omega \backslash B(Ns, Q)$*

$$\frac{1}{\sigma(\Delta(s, Q))} \int_{\Delta(s,Q)} |\log k_X - (\log k_X)_{s,Q}| d\sigma \le \varepsilon.$$

The next statements follow directly from Theorem 5.1 and Theorem 5.2 above. For a proof of this general fact see [GCRdF].

Corollary 5.1. *There exist $\delta_1 > 0$ and $N_0 > 1$ such that, if $\Omega$ is a $(\delta_1,\infty)$-chord arc domain, for $N \geq N_0$, $s > 0$, $Q \in \partial\Omega$ and $X \in \Omega\backslash B(Ns,Q)$,*

$$\left(\frac{1}{\sigma(\Delta(s,Q))}\int_{\Delta(s,Q)} k_X^2 d\sigma\right)^{\frac{1}{2}} \leq 2\frac{1}{\sigma(\Delta(s,Q))}\int_{\Delta(s,Q)} k_X d\sigma. \tag{5.3}$$

Corollary 5.2. *There exists $\delta_1 > 0$ so that if $\Omega$ is a bounded $\delta_1$-chord arc domain, there exist $N_0 > 1$ and $s_0 > 0$ so that for $N \geq N_0$, $s \in (0,s_0)$, $Q \in \partial\Omega$ and $X \in \Omega\backslash B(Ns,Q)$,*

$$\left(\frac{1}{\sigma(\Delta(s,Q))}\int_{\Delta(s,Q)} k_X^2 d\sigma\right)^{\frac{1}{2}} \leq 2\frac{1}{\sigma(\Delta(s,Q))}\int_{\Delta(s,Q)} k_X d\sigma. \tag{5.4}$$

From now on every time we talk about a $(\delta,\infty)$-chord arc domain, it should be understood that $\delta \leq \delta_1$. In particular, such an $\Omega$ is an unbounded NTA domain (see [KT1] for a proof) and (5.3) holds. Corollary 3.2 guarantees that for $Q_* \in \partial\Omega$ there exists a unique doubling Radon measure $\omega^\infty = \omega$ supported on $\partial\Omega$ satisfying:

$$\int_{\partial\Omega} \varphi d\omega = \int_\Omega u\Delta\varphi \quad \forall\varphi \in C_c^\infty(\mathbf{R}^{n+1}) \tag{5.5}$$

where

$$\begin{cases} \Delta u = 0 & \text{in } \Omega \\ u > 0 & \text{in } \Omega \\ u = 0 & \text{on } \partial\Omega, \end{cases} \tag{5.6}$$

and

$$\omega(\Delta(1,Q_*)) = 1.$$

Furthermore if $Q \in \partial\Omega$ and $s > 0$ for $Y \in \Omega\backslash B(Ns,Q)$, $N \geq N_0$, Corollary 5.1 implies that

$$\begin{aligned}\left(\int_{\Delta(s,Q)} \left(\frac{k_Y}{G(Y,A_*)}\right)^2 d\sigma\right)^{\frac{1}{2}} &\leq 2\left(\sigma(\Delta(s,Q))\right)^{-\frac{1}{2}}\int_{\Delta(s,Q)}\frac{k_Y}{G(Y,A_*)}d\sigma, \\ &\leq 2\left(\sigma(\Delta(s,Q))\right)^{-\frac{1}{2}}\frac{\omega^Y(\Delta(s,Q))}{G(Y,A_*)},\end{aligned}$$

where $A_* = A(1,Q_*)$. Corollary 3.2 and its proof guarantee that

$$\frac{\omega^Y}{G(Y,A_*)} \rightharpoonup \mu = \kappa_0\omega, \quad \text{as } |Y| \to \infty \tag{5.7}$$

where $\kappa_0 = \mu(\Delta(1,Q_*))$. Hence for a compact set $K \subset B(s,Q)$,

$$\sup_{Y\in\Omega\backslash B(Ns,Q)} \|\frac{k_Y}{G(Y,A_*)}\|_{L^2(K\cap\partial\Omega)} \leq C_K.$$

Given a sequence $\{Y_j\} \subset \Omega \backslash B(Ns, Q)$ such that $|Y_j| \to \infty$ as $j \to \infty$ there exist $h_0 \in L^2_{\mathrm{loc}}(d\sigma)$ and a subsequence $\{Y_{j'}\}$ such that

$$\frac{k_{Y_{j'}}}{G(Y_{j'}, A_*)} \rightharpoonup h_0 \quad \text{in} \quad L^2(\Delta(s, Q)).$$

Letting $s \to \infty$, and taking a diagonal subsequence $\{Y_{j_k}\}$ of $\{Y_{j'}\}$, we conclude that

$$\frac{k_{Y_{j_k}}}{G(Y_{j_k}, A_*)} \rightharpoonup h_0 \quad \text{in} \quad L^2_{\mathrm{loc}}(\partial\Omega). \tag{5.8}$$

Combining (5.7) and (5.8) we conclude that if $h = \frac{h_0}{\kappa_0}$

$$\frac{k_Y}{G(Y, A_*)} \rightharpoonup \kappa_0 h \ \text{ in } \ L^2_{\mathrm{loc}}(\partial\Omega) \ \text{ as } \ |Y| \to \infty \ \text{ and } \ d\omega = h\, d\sigma. \tag{5.9}$$

In particular (5.5) becomes

$$\int_{\partial\Omega} \varphi h d\sigma = \int_\Omega u \Delta\varphi \quad \forall \varphi \in C_c^\infty(\mathbf{R}^{n+1}). \tag{5.10}$$

Here $h$ denotes the Poisson kernel of $\Omega$ with pole at $\infty$. In particular $\omega \in A_\infty(\sigma)$, and for $s > 0$ and $Q \in \partial\Omega$

$$\left(\frac{1}{\sigma(\Delta(s, Q))} \int_{\Delta(s,Q)} h^2 d\sigma\right)^{\frac{1}{2}} \le 2 \frac{1}{\sigma(\Delta(s, Q))} \int_{\Delta(s,Q)} h d\sigma.$$

This follows from the corresponding inequality for $\frac{k_Y}{G(Y, A_*)}$ above, and the lower semi-continuity of the $L^2$ norm under weak convergence.

Main Theorem. *Let $\Omega \subset \mathbf{R}^{n+1}$ be a $(\delta_1, \infty)$-chord arc domain. Assume that $\omega$ is asymptotically optimally doubling and that $\log h \in \mathrm{VMO}(d\sigma)$. Then there exists $\delta(n) > 0$ such that, if $\Omega$ is a $\delta$-chord arc domain, $\Omega$ is a chord arc domain with vanishing constant* (*i.e.* $\overrightarrow{n} \in \mathrm{VMO}(\partial\Omega)$).

We note that the condition that $\Omega$ be a $(\delta_1, \infty)$-chord arc domain is a global condition on $\partial\Omega$, which should be understood as a technical condition. It is only used to guarantee that $\Omega$ is an unbounded NTA domain and that Corollary 5.1 holds in its scale invariant form. The essence of the theorem is local, in the sense that both the hypotheses and the conclusions are local statements that are satisfied uniformly on compact sets (see for instance Theorem 5.4). The main theorem is an easy consequence of the following lemma which is a decay type estimate. (See Definition 1.8 (vector-valued version) for the definition of $\sup_{Q \in K \cap \partial\Omega} \|\overrightarrow{n}\|_*(B(R, Q)))$.

MAIN LEMMA. *Let $\Omega \subset \mathbf{R}^{n+1}$ be a $(\delta_1, \infty)$-chord arc domain. Assume that $\omega$ is asymptotically optimally doubling and that $\log h \in \mathrm{VMO}(d\sigma)$. Then there exists $\delta(n) > 0$ so that for $\delta \in (0, \delta(n)]$, $K \subset \mathbf{R}^{n+1}$ a compact set and $R > 0$ if*

$$\sup_{Q \in K \cap \partial\Omega} \|\overrightarrow{n}\|_*(B(R, Q)) \leq \delta;$$

*then there exists $r > 0$ such that*

$$\sup_{Q \in K \cap \partial\Omega} \|\overrightarrow{n}\|_*(B(r, Q)) \leq \frac{\delta}{2}.$$

We now indicate how the main theorem follows from the main lemma. Assume that $\Omega$ satisfies the hypotheses of the main theorem, and that $\Omega$ is a $\delta$-chord arc domain, for $\delta \in (0, \delta(n)]$ (where $\delta(n)$ is as in the Main Lemma). For each compact set $K \subset \mathbf{R}^{n+1}$, there exists $r_0 > 0$, such that

$$\sup_{Q \in K \cap \partial\Omega} \|\overrightarrow{n}\|_*(B(r_0, Q)) \leq \delta.$$

Then there exists $r_1 \in (0, r_0)$ such that

$$\sup_{Q \in K \cap \partial\Omega} \|\overrightarrow{n}\|_*(B(r_1, Q)) \leq \frac{\delta}{2}.$$

Applying the main lemma inductively, we construct a sequence $\{r_k\}_k$, with $r_k \downarrow 0$ satisfying

$$\sup_{Q \in K \cap \partial\Omega} \|\overrightarrow{n}\|_*(B(r_k, Q)) \leq \frac{\delta}{2^k}.$$

Since

$$\sup_{Q \in K \cap \partial\Omega} \|\overrightarrow{n}\|_*(B(s, Q)) \leq \sup_{Q \in K \cap \partial\Omega} \|\overrightarrow{n}\|_*(B(r, Q)),$$

wherever $0 < s \leq r$, we conclude that

$$\lim_{r \to 0} \sup_{Q \in K \cap \partial\Omega} \|\overrightarrow{n}\|_*(B(r, Q)) = 0. \qquad \square$$

The rest of Section 5 and all of Section 6 are devoted to the proof of the main lemma. There are several ideas behind this proof. First, since $\omega$ is asymptotically optimally doubling, and $\Omega$ is Reifenberg flat, then $\Omega$ is a Reifenberg flat domain with vanishing constant (see Theorem 3.4). Therefore, locally, its boundary can be well approximated in the Hausdorff distance sense by $n$-dimensional planes. This provides a way to compare the harmonic measure of $\Omega$ to the harmonic measure of the appropriate half space (see Theorem 5.5). It also gives us an improved version of the Semmes decomposition lemma (see Lemma 4.1). Second, the fact that $\log h \in \mathrm{VMO}(\partial\Omega)$ is telling us that, locally, the oscillation of $h$ is small (except for a very small set) (see Lemma 5.6). This allows us to compare the Poisson kernel of $\Omega$ to the Poisson kernel of the appropriate half space, in a large set, which is also "regular" from the geometric point of view (see Corollary 5.13). Third, we follow Jerison's idea (see [J]) to use Rellich's identity (see [JK3]) for $\partial\Omega$, and for the boundary of the half space mentioned above (see Corollaries 6.2 and 6.3). An integration by parts, combined with our improved version of the Semmes decomposition lemma (see Theorem 6.1 and the proof of inequality (6.26)) and the fact that the density function (see [Si] for the definition) is bounded below by a constant very close to 1 on $\partial\Omega$, lead to the proof of the main lemma (see inequality 6.55 and Corollary 6.6).

Note that there is nothing special about the harmonic measure with pole at infinity and its corresponding Poisson kernel. In particular, similar results hold for bounded and unbounded chord arc domains with small constant for which the harmonic measure with finite pole is asymptotically optimally doubling and the logarithm of the corresponding Poisson kernel has vanishing mean oscillation. Roughly speaking, since all the results are local, once we focus our attention on small balls centered at the boundary, the distance from the pole to the boundary becomes very large with respect to the radius of the balls, and hence similar results hold.

THEOREM 5.3. *Let $\Omega \subset \mathbf{R}^{n+1}$ be a $(\delta_1, \infty)$-chord arc domain. There exists $\delta > 0$ such that if $\Omega$ is a $\delta$-chord arc domain, $\omega^X$ is asymptotically optimally doubling and $\log k_X \in \mathrm{VMO}(d\sigma)$, for some $X \in \Omega$, then $\Omega$ is a chord arc domain with vanishing constant.*

THEOREM 5.4. *There exists $\delta > 0$, such that if $\Omega$ is a bounded $\delta$-chord arc domain, $\omega^X$ is asymptotically optimally doubling, and $\log k_X \in \mathrm{VMO}(d\sigma)$, for some $X \in \Omega$, then $\Omega$ is a chord arc domain with vanishing constant.*

The proofs of Theorems 5.3 and 5.4 are very similar to the proof of the main theorem. Thus, as we go along, we indicate how to modify the arguments where the proofs differ.

Let $Q_0 \in \mathbf{R}^{n+1}$, $M > 1$, $\mathcal{K} > 4$, $s > 0$. Let $L(M\mathcal{K}s, Q_0)$ be an $n$-dimensional plane containing $Q_0$. Let $\overrightarrow{n}(M\mathcal{K}s, Q_0)$ be a unit normal vector to $L(M\mathcal{K}s, Q_0)$. We denote by

$$\begin{aligned}
\mathcal{C}(M\mathcal{K}s, Q_0) &= \{(x,t) = x + t\,\overrightarrow{n}(M\mathcal{K}s, Q_0) : x \in L(M\mathcal{K}s, Q_0), \\
&\qquad\qquad |x - Q_0| \le M\mathcal{K}s,\ |t| \le M\mathcal{K}s\}, \\
\mathcal{C}^+(M\mathcal{K}s, Q_0) &= \{(x,t) = x + t\,\overrightarrow{n}(M\mathcal{K}s, Q_0) : x \in L(M\mathcal{K}s, Q_0), \\
&\qquad\qquad |x - Q_0| \le M\mathcal{K}s,\ 0 \le t \le M\mathcal{K}s\}, \\
A(M\mathcal{K}s, Q_0) &= A = Q_0 + s\,\overrightarrow{n}(M\mathcal{K}s, Q_0) \in \mathcal{C}^+(M\mathcal{K}s, Q_0),
\end{aligned}$$

and by $\omega_+$ the harmonic measure of $\mathcal{C}^+(M\mathcal{K}s, Q_0)$. We denote by $\omega_0$ the harmonic measure of the half space

$$\{(x,t) = x + t\overrightarrow{n}(M\mathcal{K}s, Q_0) \in \mathbf{R}^{n+1} : x \in L(M\mathcal{K}s, Q_0), t \geq 0\}.$$

Recall that

$$a \underset{\varepsilon}{\sim} b \quad \text{if and only if} \quad \frac{1}{1+\varepsilon} \leq \frac{a}{b} \leq 1+\varepsilon.$$

The results in Lemmas 5.1, 5.2, 5.3 and Theorem 5.5 are essentially contained in [KT1], although the dependence on the parameters is slightly different. The following lemma describes the doubling properties of the harmonic measures of a cylinder and of a half space.

LEMMA 5.1. *Given $\varepsilon > 0$ and $M > 1$, there exists $\theta_1(\varepsilon, M) \in (0,1)$ such that for every $s > 0$, $\mathcal{K} > 4$, $Q \in L(M\mathcal{K}s, Q_0) \cap B(Ms, Q_0)$, and $r_1, r_2 \in (0, \theta_1 s)$ the following hold*:

$$\frac{\omega_+^A(\Delta_+(r_1, Q))}{\omega_+^A(\Delta_+(r_2, Q))} \underset{\varepsilon}{\sim} \left(\frac{r_1}{r_2}\right)^n,$$

*and*

$$\frac{\omega_0^A(\Delta_+(r_1, Q))}{\omega_0^A(\Delta_+(r_2, Q))} \underset{\varepsilon}{\sim} \left(\frac{r_1}{r_2}\right)^n. \tag{5.11}$$

*Here $\Delta_+(r, Q) = B(r, Q) \cap L(M\mathcal{K}s, Q_0)$.*

*Proof of Lemma* 5.1. Without loss of generality we may assume that $Ms = 1$. Let $r_i \in (0, \frac{\theta}{M})$ for some $\theta \in (0,1)$ to be determined. Let $G_+(A, -)$ denote the Green's function of $\mathcal{C}^+(\mathcal{K}, Q_0)$ with pole at $A = Q_0 + \frac{1}{M}\overrightarrow{n}(M\mathcal{K}s, Q_0)$ and $\frac{\partial G_+(A,-)}{\partial n}$ denote its normal derivative. Let $A_{1/2} = A = Q_0 + \frac{1}{2M}\overrightarrow{n}(M\mathcal{K}s, Q_0)$. Note that $\sup_{\{Y=(x,t):|x|\leq 2, 0\leq t\leq \frac{1}{2M}\}} G_+(A, Y) \leq C_n M^{-(n-2)}$, and $G_+(A, A_{1/2}) \sim M^{-(n-2)}$. The boundary regularity theory guarantees that since $\mathcal{K} > 4$, $\log \frac{\partial G_+(A,-)}{\partial n}$ is a $C^\infty$ function on $\Delta_+(2, Q_0)$. The Hopf boundary point lemma asserts that for $Q \in \Delta_+(1, Q_0)$, $\frac{\partial G_+(A,Q)}{\partial n} \geq C(n, M)$. For $Q \in \Delta_+(1, Q_0)$ and $P \in \Delta_+(r_i, Q) \subset \Delta_+(2, Q_0)$,

$$\frac{\left|\frac{\partial G_+(A,P)}{\partial n} - \frac{\partial G_+(A,Q)}{\partial n}\right|}{\frac{\partial G_+(A,Q)}{\partial n}} \leq C(n,M)|P - Q| \leq C(n,M)\frac{\theta}{M} = C(n,M)\theta.$$

Since

$$\frac{\omega_+^A(\Delta_+(r_1, Q))}{\omega_+^A(\Delta_+(r_2, Q))} = \frac{\int_{\Delta_+(r_1,Q)} \frac{\partial G_+(A,P)}{\partial n} d\sigma}{\int_{\Delta_+(r_2,Q)} \frac{\partial G_+(A,P)}{\partial n} d\sigma},$$

then

$$\left(\frac{1 - C(n,M)\theta}{1 + C(n,M)\theta}\right)\left(\frac{r_1}{r_2}\right)^n \leq \frac{\omega_+^A(\Delta_+(r_1, Q))}{\omega_+^A(\Delta_+(r_2, Q))} \leq \left(\frac{1 + C(n,M)\theta}{1 - C(n,M)\theta}\right)\left(\frac{r_1}{r_2}\right)^n.$$

Choose $\theta_1(\varepsilon, M) \in (0,1)$ so that for $\theta \le \theta_1$, $\frac{1+C(n,M)\theta}{1-C(n,M)\theta} \le 1+\varepsilon$. The proof of (5.11) is identical to the previous one. □

*Remark* 5.1. Note that if $\omega$ is asymptotically optimally doubling and $\Omega$ is Reifenberg flat then $\Omega$ is Reifenberg flat with vanishing constant. Let $K \subset \mathbf{R}^{n+1}$ be a fixed compact set. Let $K' = \overline{(K,1)}$, note that $K'$ is a compact set (namely the closure of the 1-neighborhood of $K$). Then for $\zeta > 0$, $M > 1$ and $\mathcal{K} > 4$ there exists $s = s(\zeta, M, \mathcal{K}) > 0$ such that $\sqrt{n+1}M\mathcal{K}s \le 1$,

$$\sup_{0<r\le\sqrt{n+1}M\mathcal{K}s} \theta_{K'}(r) = \sup_{0<r\le\sqrt{n+1}M\mathcal{K}s}\ \sup_{Q\in\partial\Omega\cap K'} \theta(r,Q) < \frac{1}{16M\mathcal{K}}, \tag{5.12}$$

and

$$\sup_{0<r\le\sqrt{n+1}M\mathcal{K}s} \theta_{K'}(r) < \frac{\zeta}{\sqrt{n+1}}.$$

If $Q_0 \in K \cap \partial\Omega$ there exists an $n$-dimensional plane $L(M\mathcal{K}s, Q_0)$ containing $Q_0$ and such that

$$\frac{1}{M\mathcal{K}s} D[\partial\Omega \cap B(\sqrt{n+1}M\mathcal{K}s, Q_0), L(M\mathcal{K}s, Q_0) \cap B(\sqrt{n+1}M\mathcal{K}s, Q_0)] < \zeta.$$

Note that $\mathcal{C}(M\mathcal{K}s, Q_0) \subset B(\sqrt{n+1}M\mathcal{K}s, Q_0)$. We assume that $\overrightarrow{n}(M\mathcal{K}s, Q_0)$, the unit normal vector to $L(M\mathcal{K}s, Q_0)$, has been chosen so that if $(x,t) = x + t\,\overrightarrow{n}(M\mathcal{K}s, Q_0)$, then

$$B(\sqrt{n+1}M\mathcal{K}s, Q_0) \cap \{(x,t) \in \mathbf{R}^{n+1} : x \in L(M\mathcal{K}s, Q_0), t \ge 2\zeta M\mathcal{K}s\} \subset \Omega.$$

We denote by $\widetilde{\Omega}(M\mathcal{K}s, Q_0) = \Omega \cap \mathcal{C}(M\mathcal{K}s, Q_0)$ and by $\widetilde{\omega}$ its harmonic measure. Even though $\widetilde{\Omega}(M\mathcal{K}s, Q_0)$ might not be an NTA domain, all points in $\Delta(\frac{M\mathcal{K}s}{2}, Q_0)$ are nontangentially accessible. Moreover the Harnack chain condition holds on $B(\frac{M\mathcal{K}s}{2}, Q_0) \cap \Omega$. Thus the results quoted at the beginning of Section 3 hold for $\widetilde{\omega}$ on $B(\frac{M\mathcal{K}s}{2}, Q_0) \cap \Omega$. Let $\Pi : \mathbf{R}^{n+1} \to L(M\mathcal{K}s, Q_0)$ be the orthogonal projection. Note that since (5.12) holds $A = Q_0 + s\,\overrightarrow{n}(M\mathcal{K}s, Q_0) \in \widetilde{\Omega}(M\mathcal{K}s, Q_0)$. Note also that if $\theta \in (0, \frac{1}{64})$, $Q \in \Delta(Ms, Q_0)$ and $L(\theta s, Q)$ is an $n$-plane containing $Q$, satisfying the conditions described in Definition 1.9, and such that

$$\frac{1}{\theta s} D[\partial\Omega \cap B(\theta s, Q), L(\theta s, Q) \cap B(\theta s, Q)] < \zeta,$$

then

$$\begin{aligned} &d(L(M\mathcal{K}s, Q_0), L(\theta s, Q)) \\ &\qquad = D[(L(M\mathcal{K}s, Q_0) - Q_0) \cap B(1,0), (L(\theta s, Q) - Q) \cap B(1,0)] < \frac{3M\mathcal{K}\zeta}{4\theta}, \end{aligned}$$

(see proof of Lemma 3.1 in [T]). Therefore if $\Lambda(M\mathcal{K}s, Q) = L(M\mathcal{K}s, Q_0) - (Q_0 - Q)$,

$$\frac{1}{\theta s} D[\partial\Omega \cap B(\theta s, Q), \Lambda(M\mathcal{K}s, Q) \cap B(\theta s, Q)] < \frac{M\mathcal{K}\zeta}{\theta}.$$

Moreover, the same argument proves that for $1/2 \leq \tau \leq \sqrt{1-\left(\frac{4M\mathcal{K}\zeta}{\theta}\right)^2}$,

$$\frac{1}{\theta\tau s} D[\partial\Omega \cap B(\theta\tau s, Q), \Lambda(M\mathcal{K}s, Q) \cap B(\theta\tau s, Q)] \leq \frac{2M\mathcal{K}\zeta}{\theta}.$$

It is not difficult to see that, for $\zeta > 0$ small, and the appropriate choice of unit normal,

$$\{X = (x,t) = x + t\,\overrightarrow{n}(M\mathcal{K}s, Q_0) : \\ x \in \Lambda(M\mathcal{K}s, Q),\ t > 4\frac{M\mathcal{K}\zeta}{\theta}s\} \cap\ B(\theta\tau s, Q) \subset \Omega,$$

and

$$\{X = (x,t) = x + t\,\overrightarrow{n}(M\mathcal{K}s, Q_0) : \\ x \in \Lambda(M\mathcal{K}s, Q),\ t < -4\frac{M\mathcal{K}\zeta}{\theta}s\} \cap\ B(\theta\tau s, Q) \subset \Omega^c,$$

Remark 4.1 guarantees that

$$\Gamma(\theta\tau s, Q) \subset \Delta(\theta s, Q) \subset \Gamma(\theta s, Q),$$

where

$$\Gamma(\theta\tau s, Q) = \partial\Omega \cap \{(x,t) :\ x \in L(M\mathcal{K}s, Q_0), |x - \Pi(Q)| \leq \theta\tau s, |t| \leq \theta\tau s\}.$$

Since $L(M\mathcal{K}s, Q_0)$ and $\Lambda(M\mathcal{K}s, Q)$ are parallel planes, Remark 4.1 also guarantees that if $\Pi : \mathbf{R}^{n+1} \to L(M\mathcal{K}s, Q_0)$ denotes the orthogonal projection onto $L(M\mathcal{K}s, Q_0)$, then

$$\Pi\left(\Gamma(\theta\tau s, Q)\right) = \{x \in L(M\mathcal{K}s, Q_0),\ |x - \Pi(Q)| \leq \theta\tau s\}.$$

The notation introduced above is adopted for the rest of the paper. In particular, we also fix $\zeta > 0$ as above. We would like to warn the reader that in several of the results stated in Sections 5 and 6, a parameter $\mathcal{K} > 4$ appears in the hypothesis but seems to vanish in the conclusions. The reader should keep in mind that the point $A$ defined above always depends on $\mathcal{K}$. Thus in most such statements the dependence on $\mathcal{K}$ is hidden in $A$. The key fact to remember about the point $A$ is that it can be neither too close nor too far from the boundary as controlled by the various parameters which appear above, namely $\theta$, $M$, and $\mathcal{K}$. Note for example that this balance in the relative distance of the point $A$ to the boundary is what allows us to insure that $A = Q_0 + s\,\overrightarrow{n}(M\mathcal{K}s, Q_0) \in \widetilde{\Omega}(M\mathcal{K}s, Q_0)$. One should also remember that all the constants depend implicitly on the compact set $K \subset \mathbf{R}^{n+1}$ which is fixed, and on the modulus of flatness in the Reifenberg vanishing condition.

LEMMA 5.2. *Let $\Omega$ be a Reifenberg flat domain with vanishing constant. Given $\varepsilon > 0$, $M > 1$ and $\theta \in (0,1)$ there exists $\mathcal{K}(\varepsilon, M, \theta) > 4$ such that if $\mathcal{K} \geq \mathcal{K}(\varepsilon, M, \theta)$ there is $s(M, \theta, \mathcal{K}) > 0$ so that for $s \in (0, s(M, \theta, \mathcal{K}))$ and $Q \in \Delta(Ms, Q_0)$,*

$$\frac{\widetilde{\omega}^A(\Delta(\theta s, Q))}{\omega^A(\Delta(\theta s, Q))} \underset{\varepsilon}{\sim} 1. \tag{5.13}$$

Note that in particular this lemma applies to the half space $\{(x,t) \in \mathbf{R}^{n+1} : x \in L(M\mathcal{K}s, Q_0), t \geq 0\}$ and the corresponding cylinder $\mathcal{C}^+(M\mathcal{K}s, Q_0)$. That is, using the notation introduced above we have:

COROLLARY 5.3. *Given $\varepsilon > 0$, $M > 1$ and $\theta \in (0,1)$ there exists $\mathcal{K}(\varepsilon, M, \theta) > 4$ such that if $\mathcal{K} \geq \mathcal{K}(\varepsilon, M, \theta)$ then for every $Q \in B(Ms, Q_0) \cap L(M\mathcal{K}s, Q_0)$ and every $s > 0$,*

$$\frac{\omega_+^A(\Delta_+(\theta s, Q))}{\omega_0^A(\Delta_+(\theta s, Q))} \underset{\varepsilon}{\sim} 1. \tag{5.14}$$

Note that Corollary 5.3 is true for every $s > 0$ because $\theta_{K'}(s) = 0$, and therefore the statements in Remark 5.1 hold at any scale.

*Proof of Lemma* 5.2. Since $\widetilde{\Omega}(M\mathcal{K}s, Q_0) \subset \Omega$ and for $Q \in \Delta(Ms, Q_0)$, $\Delta(\theta s, Q) \subset \partial\widetilde{\Omega}(M\mathcal{K}s, Q_0)$, the maximum principle guarantees that for every $X \in \widetilde{\Omega}(M\mathcal{K}s, Q_0)$,

$$\widetilde{\omega}^X(\Delta(\theta s, Q)) \leq \omega^X(\Delta(\theta s, Q)).$$

Therefore $v(X) = \omega^X(\Delta(\theta s, Q)) - \widetilde{\omega}^X(\Delta(\theta s, Q))$ is a nonnegative harmonic function on $\widetilde{\Omega}(M\mathcal{K}s, Q_0)$ bounded above by 1. Moreover $v$ vanishes on $\Delta(\frac{M\mathcal{K}s}{2}, Q_0)$. Lemma 3.2 asserts that

$$v(A) \leq C \left( \sup_{X \in B(\frac{M\mathcal{K}s}{2}, Q_0) \cap \Omega} v(X) \right) \left( \frac{|A - Q_0|}{M\mathcal{K}s} \right)^\alpha \leq C^{-1} \left( \frac{1}{M\mathcal{K}} \right)^\alpha. \tag{5.15}$$

Since $|Q + s\overrightarrow{n}(M\mathcal{K}s, Q_0) - A| = |Q - Q_0| < Ms \leq 2^k s$ the Harnack chain condition insures that $A$ can be joined to $Q + s\overrightarrow{n}(M\mathcal{K}s, Q_0) \in \Omega$ by a chain of $Ck$ nontangential balls. Therefore Harnack's principle and Corollary 3.1 guarantee that for every $Q \in \Delta(Ms, Q_0)$,

$$\begin{aligned} \omega^A(\Delta(s, Q)) &\geq 3^{-(n+1)Ck} \omega^{Q + s\overrightarrow{n}(M\mathcal{K}s, Q_0)}(\Delta(s, Q_0)), \\ \omega^A(\Delta(s, Q)) &\geq \frac{C}{M^{C(n+1)}} \geq \frac{C}{M^C}. \end{aligned} \tag{5.16}$$

Combining Lemma 3.1, an iteration argument, and (5.16) we have that there exists $q \geq 1$ depending on the NTA constants, and on the compact set

$K \subset \mathbf{R}^{n+1}$, so that

$$\omega^A(\Delta(\theta s, Q)) \geq C^{-1}\theta^q \omega^A(\Delta(s,Q)) \geq C^{-1}\frac{\theta^q}{M^C}. \tag{5.17}$$

(5.15) and (5.17) yield

$$\begin{aligned} \omega^A(\Delta(\theta s, Q)) &\geq C^{-1}\frac{\theta^q}{M^C}(M\mathcal{K})^\alpha v(A), \\ \omega^A(\Delta(\theta s, Q)) - \widetilde{\omega}^A(\Delta(\theta s, Q)) &\leq C\frac{M^C}{\theta^q(M\mathcal{K})^\alpha}\omega^A(\Delta(\theta s, Q)), \\ \left(1 - C\frac{M^C}{\theta^q(M\mathcal{K})^\alpha}\right)\omega^A(\Delta(\theta s, Q)) &\leq \widetilde{\omega}^A(\Delta(\theta s, Q)). \end{aligned} \tag{5.18}$$

Choose $\mathcal{K}(\varepsilon, M, \theta) > 4$ large enough so that for $\mathcal{K} \geq \mathcal{K}(\varepsilon, M, \theta)$, $(1+\varepsilon)^{-1} \leq 1 - C\frac{M^C}{\theta^q(M\mathcal{K})^\alpha}$. Combining (5.18) and the maximum principle we conclude that $Q \in \Delta(Ms, Q_0)$ and $\mathcal{K} \geq \mathcal{K}(\varepsilon, M, \theta)$,

$$(1+\varepsilon)^{-1}\omega^A(\Delta(\theta s, Q)) \leq \widetilde{\omega}^A(\Delta(\theta s, Q)) \leq \omega^A(\Delta(\theta s, Q)). \qquad \square$$

LEMMA 5.3. *Let $\Omega$ be a Reifenberg flat domain with vanishing constant. Given $\varepsilon > 0$ and $M > 1$, there exists $\theta(\varepsilon, M) \in (0,1)$ such that, for $\theta \in (0, \theta(\varepsilon, M))$ and $\mathcal{K} > 4$, there exists $s(\varepsilon, M, \theta, \mathcal{K}) > 0$ so that, if $s \in (0, s(\varepsilon, M, \theta, \mathcal{K}))$ and $Q \in \Delta(Ms, Q_0)$,*

$$\frac{\widetilde{\omega}^A(\Delta(\theta s, Q))}{\omega_+^A(\Delta_+(\theta s, \Pi(Q)))} \underset{\varepsilon}{\sim} 1. \tag{5.19}$$

*Proof.* We introduce first the notation that will be used in the proof. Let $\delta = \sqrt{n+1}\sup_{0<r\leq M\mathcal{K}\sqrt{n+1}s}\theta_{K'}(r)$, where $K' = \overline{(K,1)}$. Let

$$\begin{aligned} &\underline{\mathcal{C}}(M\mathcal{K}s, Q_0) \\ &\quad = \{(x,t) : x \in L(M\mathcal{K}s, Q_0),\ |x - Q_0| \leq M\mathcal{K}s,\ 2\delta M\mathcal{K}s \leq t \leq M\mathcal{K}s\}, \\ &\overline{\mathcal{C}}(M\mathcal{K}s, Q_0) \\ &\quad = \{(x,t) : x \in L(M\mathcal{K}s, Q_0),\ |x - Q_0| \leq M\mathcal{K}s,\ -2\delta M\mathcal{K}s \leq t \leq M\mathcal{K}s\}. \end{aligned}$$

We denote by $\underline{\omega}$ (resp. $\overline{\omega}$) the harmonic measure of $\underline{\mathcal{C}}(M\mathcal{K}s, Q_0)$ (resp. $\overline{\mathcal{C}}(M\mathcal{K}s, Q_0)$). Since $M\mathcal{K}\delta < \frac{1}{16}$ (see (5.12))

$$A \in \underline{\mathcal{C}}(M\mathcal{K}s, Q_0) \subset \widetilde{\Omega}(M\mathcal{K}s, Q_0) \subset \overline{\mathcal{C}}(M\mathcal{K}s, Q_0).$$

We denote by

$$\begin{aligned} \underline{Q} &= \Pi(Q) + 2\delta M\mathcal{K}s\,\overrightarrow{n}(M\mathcal{K}s, Q_0), \\ \overline{Q} &= \Pi(Q) - 2\delta M\mathcal{K}s\,\overrightarrow{n}(M\mathcal{K}s, Q_0), \\ \underline{\Delta}(\rho, \underline{Q}) &= B(\rho, \underline{Q}) \cap \{(x,t) : t = 2\delta M\mathcal{K}s\}, \end{aligned}$$

and

$$\overline{\Delta}(\rho, \overline{Q}) = B(\rho, \overline{Q}) \cap \{(x,t) : t = -2\delta M\mathcal{K}s\}.$$

Let $\varepsilon' = \varepsilon'(\varepsilon) > 0$ to be chosen later, let $\theta(\varepsilon', M) > 0$ be as in Lemma 5.1. Let $\theta \in (0, \frac{1}{2}\theta(\varepsilon', M))$ and $\eta = \eta(\varepsilon') \in (0, \frac{1}{2})$ to be chosen later. Since $\Omega$ is a Reifenberg flat domain with vanishing constant we can choose $s = s(\varepsilon', M, \mathcal{K}, \theta)$ so that $\delta > 0$ satisfies $\sqrt{1 - (\frac{4\delta M\mathcal{K}}{\theta})^2} \geq (1-\eta)^2$, $\frac{16 M\mathcal{K}\delta}{\theta} \leq \eta < \frac{1}{2}$, and Remark 5.1 holds for $\zeta = \delta$. For $Y \in \underline{\Delta}(\theta s(1-\eta)^2, \underline{Q})$ there exists $Z \in \Gamma(\theta s(1-\eta)^2, Q)$, such that $\Pi(Y) = \Pi(Z)$ and $|Y - Z| \leq 4\delta M\mathcal{K}s \leq \frac{\eta\theta s}{4}$ (see Remark 5.1). Here

$$\begin{aligned}\Gamma(\theta s(1-\eta)^2, Q) \quad = \quad & \partial\Omega \cap \{(x,t) : \ x \in L(M\mathcal{K}s, Q_0), \\ & |x - \Pi(Q)| \leq \ \theta s(1-\eta)^2, |t| \leq \ \theta s(1-\eta)^2\}.\end{aligned}$$

Since $\Gamma(\frac{\eta\theta s}{2}, Z) \subset \Gamma(\theta s(1-\eta), Q)$, for $Z' \in \Gamma(\frac{\eta\theta s}{2}, Z)$, $1 - \widetilde{\omega}^{Z'}(\Gamma(\theta s(1-\eta), Q)) = 0$. Lemma 3.2 guarantees then that for $Y' \in \{(x,t) : |x - \Pi(Z)| \leq \frac{\eta\theta s}{4}, |t| \leq \frac{\eta\theta s}{4}\} \cap \widetilde{\Omega}(M\mathcal{K}s, Q_0)$,

$$1 - \widetilde{\omega}^{Y'}(\Gamma(\theta s(1-\eta), Q)) \leq C\left(\frac{|Y' - Z|}{\eta\theta s}\right)^{\beta}. \tag{5.20}$$

By our choice of $s$, inequality (5.20) holds for any $Y \in \underline{\Delta}(\theta s(1-\eta)^2, \underline{Q})$; hence

$$1 - \widetilde{\omega}^{Y}(\Gamma(\theta s(1-\eta), Q)) \leq C\left(\frac{4\delta M\mathcal{K}}{\eta\theta}\right)^{\beta},$$

and

$$\left(1 - C\left(\frac{4\delta M\mathcal{K}}{\eta\theta}\right)^{\beta}\right)\underline{\omega}^{Y}(\underline{\Delta}(\theta s(1-\eta)^2, \underline{Q})) \leq \widetilde{\omega}^{Y}(\Gamma(\theta s(1-\eta), Q)). \tag{5.21}$$

Thus by the maximum principle, (5.21) holds for all $X \in \underline{\mathcal{C}}(M\mathcal{K}s, Q_0)$. Note that Remark 5.1 guarantees that if $\Lambda(M\mathcal{K}s, Q) = L(M\mathcal{K}s, Q_0) - (Q_0 - Q) \ni Q$ then

$$\frac{1}{\theta s} D[\partial\Omega \cap B(\theta s, Q), \Lambda(M\mathcal{K}s, Q) \cap B(\theta s, Q)] \leq \frac{M\mathcal{K}\delta}{\theta}.$$

Thus

$$\Gamma(\theta s(1-\eta), Q)) \subset \Gamma(\theta s\sqrt{1 - (\frac{4\delta M\mathcal{K}}{\theta})^2}, Q) \subset \Delta(\theta s, Q) \subset \Gamma(\theta s, Q).$$

Therefore

$$\left(1 - C\left(\frac{4\delta M\mathcal{K}}{\eta\theta}\right)^{\beta}\right)\underline{\omega}^{A}(\underline{\Delta}(\theta s(1-\eta)^2, \underline{Q})) \leq \widetilde{\omega}^{A}(\Delta(\theta s, Q)).$$

Lemma 5.1 applied on $\underline{\mathcal{C}}(M\mathcal{K}s, Q_0)$ for $r_1 = \theta s(1-\eta)^2$ and $r_2 = \theta s$ yields

$$\left(1 - C\left(\frac{4\delta M\mathcal{K}}{\eta\theta}\right)^{\beta}\right)(1+\varepsilon')^{-1}(1-\eta)^{2n}\underline{\omega}^{A}(\underline{\Delta}(\theta s, \underline{Q})) \leq \widetilde{\omega}^{A}(\Delta(\theta s, Q)). \tag{5.22}$$

Now for $Y \in \Delta(\theta s, Q)$ there is $Z \in \overline{\Delta}(\theta s, \overline{Q})$ such that $\Pi(Y) = \Pi(Z)$, and $|Y - Z| \leq 4\delta M\mathcal{K}s$. For $Z' \in \overline{\Delta}(\frac{\eta\theta s}{2}, Z) \subset \overline{\Delta}(\frac{\theta s}{1-\eta}, \overline{Q})$, $1 - \overline{\omega}^{Z'}(\overline{\Delta}(\frac{\theta s}{1-\eta}, \overline{Q})) = 0$. Lemma 3.2 guarantees that for $Y' \in \{(x,t) \in \mathbf{R}^{n+1} : |x - \Pi(Z)| \leq \frac{\eta\theta s}{4}, |t| \leq \frac{\eta\theta s}{4}\} \cap \overline{\mathcal{C}}(M\mathcal{K}s, Q_0)$,

$$1 - \overline{\omega}^{Y'}(\overline{\Delta}(\frac{\theta s}{1-\eta}, \overline{Q})) \leq C \left(\frac{|Y' - Z|}{\eta\theta s}\right)^{\beta}. \tag{5.23}$$

By our choice of $s$ inequality (5.23) holds for $Y \in \Delta(\theta s, Q)$. Hence

$$\left(1 - C\left(\frac{\delta M\mathcal{K}}{\eta\theta}\right)^{\beta}\right) \widetilde{\omega}^{Y}(\Delta(\theta s, Q)) \leq \overline{\omega}^{Y}(\overline{\Delta}(\frac{\theta s}{1-\eta}, \overline{Q})). \tag{5.24}$$

By the maximum principle inequality (5.24) is valid for all $X \in \widetilde{\Omega}(M\mathcal{K}s, Q_0)$. In particular for $X = A$, applying Lemma 5.1 on $\overline{\mathcal{C}}(M\mathcal{K}s, Q_0)$ for $r_1 = \frac{\theta s}{1-\eta}$ and $r_2 = \theta s$ we have that

$$\left(1 - C\left(\frac{\delta M\mathcal{K}}{\eta\theta}\right)^{\beta}\right) \widetilde{\omega}^{A}(\Delta(\theta s, Q)) \leq (1-\eta)^{-n}(1+\varepsilon')\overline{\omega}^{A}(\overline{\Delta}(\theta s, \overline{Q})). \tag{5.25}$$

Now choose $\eta > 0$ so that $(1-\eta)^n \geq \frac{1}{1+\varepsilon'}$, and $s(M, \mathcal{K}, \varepsilon', \theta) > 0$ such that if $s \in (0, s(M, \mathcal{K}, \varepsilon', \theta))$, and if $\delta = \sqrt{n+1}\sup_{0<r\leq\sqrt{n+1}M\mathcal{K}s} \theta_{K'}(r)$ then $1 - C(\frac{\delta M\mathcal{K}}{\eta\theta})^{\beta} \geq \frac{1}{1+\varepsilon'}$. Under these conditions, combining (5.22) and (5.25) we obtain

$$(1+\varepsilon')^{-4}\underline{\omega}^{A}(\underline{\Delta}(\theta s, \underline{Q})) \leq \widetilde{\omega}^{A}(\Delta(\theta s, Q)) \leq (1+\varepsilon')^{3}\overline{\omega}^{A}(\overline{\Delta}(\theta s, \overline{Q})). \tag{5.26}$$

Our last step is to compare $\underline{\omega}^{A}(\underline{\Delta}(\theta s, \underline{Q}))$ and $\overline{\omega}^{A}(\overline{\Delta}(\theta s, \overline{Q}))$. Define

$$u_1(x,t) = \overline{\omega}^{(x,t)}(\overline{\Delta}(\theta s, \overline{Q})) \quad \text{for} \quad (x,t) \in \overline{\mathcal{C}}(M\mathcal{K}s, Q_0),$$

and

$$u_2(x,t) = \underline{\omega}^{(x,t)}(\underline{\Delta}(\theta s, \underline{Q})) \quad \text{for} \quad (x,t) \in \underline{\mathcal{C}}(M\mathcal{K}s, Q_0).$$

We want to compare $u_1(x, t - 4M\mathcal{K}\delta s)$ and $u_2(x,t)$ for $(x,t) \in \partial\underline{\mathcal{C}}(M\mathcal{K}s, Q_0)$. Note that if $t = 2M\mathcal{K}\delta s$ or $|x - Q_0| = M\mathcal{K}\delta s$ then $u_1(x, t - 4M\mathcal{K}\delta s) = u_2(x,t)$. Since $\overline{\mathcal{C}}(M\mathcal{K}s, Q_0)$ is an NTA domain, and $u_1$ is a nonnegative harmonic function on $\overline{\mathcal{C}}(M\mathcal{K}s, Q_0)$ which vanishes on $\partial\overline{\mathcal{C}}(M\mathcal{K}s, Q_0) \cap \{(x,t) \in \mathbf{R}^{n+1} : t \geq \frac{3}{4}M\mathcal{K}s\}$, and is bounded by 1, Lemma 3.2 guarantees that $u_1(x, M\mathcal{K}s - 4\delta M\mathcal{K}s) \leq C_1\delta^{\beta}$. Therefore for $(x,t) \in \partial\underline{\mathcal{C}}(M\mathcal{K}s, Q_0)$

$$u_1(x, t - 4\delta M\mathcal{K}s) \leq u_2(x,t) + C_1\delta^{\beta}\frac{t - 2\delta M\mathcal{K}s}{M\mathcal{K}s(1-2\delta)}. \tag{5.27}$$

By the maximum principle (5.27) holds for every $(x,t) \in \underline{\mathcal{C}}(M\mathcal{K}s, Q_0)$. Thus for $\delta \leq \frac{1}{64}$, (5.27), combined with an interior estimate (used to compare $u_1(A)$ with $u_1(A - 4M\mathcal{K}\delta s\overrightarrow{n}(M\mathcal{K}s, Q_0))$), gives

$$\overline{\omega}^{A}(\overline{\Delta}(\theta s, \overline{Q}) \leq \underline{\omega}^{A}(\underline{\Delta}(\theta s, \underline{Q}) + C_2\delta^{\beta}. \tag{5.28}$$

The Harnack principle and the doubling properties of the harmonic measure of NTA domains lead, as in the proof of Lemma 5.2 (see (5.17)), to the following inequality

$$\underline{\omega}^A(\underline{\Delta}(\theta s,\underline{Q}) \geq C^{-1}\frac{\theta^q}{M^C}, \tag{5.29}$$

for some $q \geq 1$. Combining (5.28) and (5.29) we have that

$$\overline{\omega}^A(\overline{\Delta}(\theta s,\overline{Q}) \leq (1 + C\delta^\beta \frac{M^C}{\theta^q})\underline{\omega}^A(\underline{\Delta}(\theta s,\underline{Q}). \tag{5.30}$$

For $\theta$ as above, choose $s(M,\mathcal{K},\varepsilon',\theta)$ so that for $s \in (0, s(M,\mathcal{K},\varepsilon',\theta))$, $C\delta^\beta\frac{M^C}{\theta^q} \leq \varepsilon'$ where $\delta = \sqrt{n+1}\sup_{0<r<\sqrt{n+1}M\mathcal{K}s}\theta_{K'}(r)$. Under these hypotheses, (5.26) and (5.30) combined yield

$$\frac{1}{(1+\varepsilon')^5}\overline{\omega}^A(\overline{\Delta}(\theta s,\overline{Q})) \leq \widetilde{\omega}^A(\Delta(\theta s,Q)) \leq (1+\varepsilon')^4\underline{\omega}^A(\underline{\Delta}(\theta s,\underline{Q})).$$

The maximum principle implies that

$$\begin{aligned} \omega_+^A(\Delta_+(\theta s,\Pi(Q))) &\leq \overline{\omega}^A(\overline{\Delta}(\theta s,\overline{Q})), \\ \underline{\omega}^A(\underline{\Delta}(\theta s,\underline{Q})) &\leq \omega_+^A(\Delta_+(\theta s,\Pi(Q))). \end{aligned}$$

Choosing $\varepsilon' > 0$ so that $(1+\varepsilon')^5 < 1+\varepsilon$, then $\theta$ and $s$ as indicated above, we conclude that

$$\frac{1}{1+\varepsilon}\omega_+^A(\Delta_+(\theta s,\Pi(Q))) \leq \widetilde{\omega}^A(\Delta(\theta s,Q)) \leq (1+\varepsilon)\omega_+^A(\Delta_+(\theta s,\Pi(Q))). \qquad \square$$

THEOREM 5.5. *Let $\Omega$ be a Reifenberg flat domain with vanishing constant. Given $\varepsilon > 0$ and $M > 1$ there exists $\theta(\varepsilon, M) \in (0,1)$ so that, for $\theta \in (0,\theta(\varepsilon,M))$, there exists $\mathcal{K}(\varepsilon,M,\theta) > 4$ such that if $\mathcal{K} \geq \mathcal{K}(\varepsilon,M,\theta)$, there exists $s(\varepsilon,M,\theta,\mathcal{K}) > 0$ so that for $s \in (0, s(\varepsilon,M,\theta,\mathcal{K}))$ and every $Q \in \Delta(Ms,Q_0)$,*

$$\frac{\omega^A(\Delta(\theta s,Q))}{\omega_0^A(\Delta_+(\theta s,\Pi(Q)))} \underset{\varepsilon}{\sim} 1.$$

*Proof.* Let $\varepsilon' = \varepsilon'(\varepsilon) > 0$ and $M > 1$. Let $\theta \in (0,\theta(\varepsilon',M))$, where $\theta(\varepsilon',M) > 0$ has been chosen as in the statement of Lemma 5.3. For such $\theta$, let $\mathcal{K}(\varepsilon',M,\theta) > 4$ be such that the statements of both Lemma 5.2 and Corollary 5.3 are satisfied. Let $\mathcal{K} \geq \mathcal{K}(\varepsilon',M,\theta)$. Let $s(\varepsilon',M,\theta,\mathcal{K}) > 0$ be so that, for $s \in (0, s(\varepsilon',M,\theta,\mathcal{K}))$ the statements of Lemmas 5.2 and 5.3 and Corollary 5.3 are verified. Then (5.13), (5.19) and (5.14) become

$$\frac{\omega^A(\Delta(\theta s,Q))}{\widetilde{\omega}^A(\Delta(\theta s,Q))} \underset{\varepsilon'}{\sim} 1,$$

$$\frac{\widetilde{\omega}^A(\Delta(\theta s,Q))}{\omega_+^A(\Delta_+(\theta s,\Pi(Q)))} \underset{\varepsilon'}{\sim} 1,$$

and

$$\frac{\omega_+^A(\Delta_+(\theta s, \Pi(Q)))}{\omega_0^A(\Delta_+(\theta s, \Pi(Q)))} \underset{\varepsilon'}{\sim} 1.$$

We conclude the proof by choosing $\varepsilon' > 0$ so that $(1+\varepsilon')^3 \leq 1+\varepsilon$, and $\theta(\varepsilon, M)$, $\mathcal{K}(\varepsilon, M, \theta)$ and $s(\varepsilon, M, \theta, \mathcal{K})$ accordingly. □

Theorem 5.5 allows us to compare the harmonic measure (with finite pole) of a Reifenberg flat domain with vanishing constant $\Omega$, with the harmonic measure of an appropriate half space. So far, the fact that the harmonic measure is asymptotically optimally doubling has only been used as a way to guarantee that $\Omega$ is a Reifenberg flat domain with vanishing constant. Thus, the same arguments can be used in the proof of the corresponding results for the harmonic measure with finite pole. In order to be able to use the hypothesis that $\log h \in \mathrm{VMO}(\partial\Omega)$ (resp. $\log k_X \in \mathrm{VMO}(\partial\Omega)$), we need to compare the Poisson kernel of $\Omega$ with pole at $A$ and $h$ (resp. $k_X$). To achieve this, we look at the kernel function with pole at infinity. The next lemma summarizes some of the properties of this kernel function on $(\delta_1, \infty)$-chord arc domains. Here $\delta_1$ has to be chosen so that Corollary 5.1 holds. Note that Lemma 5.4 is the analog of Theorem 3.1 in the case that the kernel function has pole at infinity.

LEMMA 5.4. *Let $\Omega \subset \mathbf{R}^{n+1}$ be a $(\delta_1, \infty)$-chord arc domain. Let $X \in \Omega$; then for almost every $Q \in \partial\Omega$,*

$$\frac{d\omega^X}{d\omega}(Q) = \frac{k_X(Q)}{h(Q)} = \lim_{r\to 0} \frac{\omega^X(\Delta(r,Q))}{\omega(\Delta(r,Q))} = \lim_{Z\to Q} \frac{G(X,Z)}{u(Z)}. \tag{5.31}$$

*Here $\omega^X$ denotes the harmonic measure, $G(X,-)$ denotes the Green's function, and $k_X$ the Poisson kernel for $\Omega$ with pole at $X$. Now, $\omega$, $u$ and $h$ satisfy* (5.6), (5.9) *and* (5.10). *Let $K(X,Q) = \frac{k_X(Q)}{h(Q)}$. There exist constants $C > 1$, $N_0 > 1$ and $\alpha \in (0,1)$ so that, for $s > 0$ and $Q_0 \in \partial\Omega$, if $X \in \Omega \backslash B(2Ns, Q_0)$, $N \geq N_0$, then for every $Q, Q' \in \Delta(s, Q_0)$,*

$$|K(X,Q') - K(X,Q)| \leq CK(X,Q)\left(\frac{|Q-Q'|}{s}\right)^\alpha.$$

*Proof.* Recall that $\sigma$, $\omega$ and $\omega^X$ are doubling Radon measures on $\partial\Omega$. Moreover, $d\omega^X = k_X d\sigma$, $d\omega = h\, d\sigma$ and $k_X$, $h \in L^2_{\mathrm{loc}}(d\sigma)$. Therefore, $\omega$ and $\omega^X$ are mutually absolutely continuous, and $\frac{d\omega^X}{d\omega}(Q) = \frac{k_X(Q)}{h(Q)}$. Because the Lebesgue differentiation theorem holds, the Radon-Nikodym derivative, for $\omega$-almost every $Q \in \partial\Omega$ satisfies (see [Si, §4])

$$\frac{d\omega^X}{d\omega}(Q) = \lim_{r\to 0} \frac{\omega^X(\Delta(r,Q))}{\omega(\Delta(r,Q))}.$$

Let $\varphi \in C_c^\infty(\mathbf{R}^{n+1})$, $0 \le \varphi \le 1$, $\varphi \equiv 1$ in $B(1,0)$, $\varphi \equiv 0$ in $\mathbf{R}^{n+1}\backslash B(2,0)$, $|\nabla\varphi| \le C$, and $|\Delta\varphi| \le C$. Let $\psi_r(Z) = \varphi(\frac{1}{r}(Z-Q))$. Let $u_r$ satisfy $\Delta u_r = 0$ in $\Omega$ and $u_r = \psi_r$ on $\partial\Omega$. Then for $r > 0$ small, $X \notin B(2r,Q)$, and

$$u_r(X) = \int_{\partial\Omega} \psi_r(Q)d\,\omega^X(Q) = \int_\Omega G(X,Z)\Delta\psi_r(Z)d\,Z.$$

(5.5) guarantees that

$$\begin{aligned}(5.32)\qquad \int_{\partial\Omega}\psi_r(Q)d\,\omega(Q) &= \int_\Omega u(Z)\Delta\psi_r(Z)d\,Z \\ &= \int_\Omega \frac{u(Z)}{G(X,Z)}G(X,Z)\Delta\psi_r(Z)d\,Z.\end{aligned}$$

A similar argument to the one used above shows that, for $\omega$-almost every $Q \in \partial\Omega$ (see [St, Ch. 1])

$$\frac{d\omega^X}{d\omega}(Q) = \lim_{r\to 0}\frac{u_r(X)}{\int_{\partial\Omega}\psi_r(Q)d\,\omega(Q)}.$$

Theorem 3.2 guarantees that $\frac{u(-)}{G(X,-)}$ is a $C^\alpha$ function on $\overline{B}(Ns,Q)\cap\overline{\Omega}$, where $Q \in \partial\Omega$, $2Ns \le \operatorname{dist}(X,\partial\Omega)$, and $N \ge N_0$ where $N_0$ is chosen so that for $X \in \Omega\backslash B(2Ns,Q)$, and $r \in [0,s]$, $\omega^X(\Delta(2r,Q)) \le C\omega^X(\Delta(r,Q))$ (see Lemma 3.1). Thus if

$$\lim_{Z\to Q}\frac{u(Z)}{G(X,Z)} = \ell(Q),$$

then for $Z \in B(2r,Q)\cap\Omega$, and $r < s$

$$(5.33)\qquad \Big|\frac{u(Z)}{G(X,Z)} - \ell(Q)\Big| \le C\frac{u(A(Ns,Q))}{G(X,A(Ns,Q))}\left(\frac{|Q-Z|}{Ns}\right)^\alpha \le C\left(\frac{r}{Ns}\right)^\alpha.$$

Thus combining (5.32) and (5.33) we obtain

$$\begin{aligned}(5.34)\quad \int_\Omega u(Z)\Delta\psi_r(Z)d\,Z &\ge \ell(Q)\int_\Omega G(X,Z)\Delta\psi_r(Z)d\,Z \\ &\quad - C\left(\frac{r}{Ns}\right)^\alpha\int_\Omega G(X,Z)|\Delta\psi_r(Z)|d\,Z,\end{aligned}$$

and

$$\begin{aligned}(5.35)\quad \int_\Omega u(Z)\Delta\psi_r(Z)d\,Z &\le \ell(Q)\int_\Omega G(X,Z)\Delta\psi_r(Z)d\,Z \\ &\quad + C\left(\frac{r}{Ns}\right)^\alpha\int_\Omega G(X,Z)|\Delta\psi_r(Z)|d\,Z.\end{aligned}$$

Combining Lemmas 3.3 and 3.4 we have that for $Z \in B(2r,Q)\cap\Omega$,

$$G(X,Z) \le CG(X,A(2r,Q)) \le C\frac{\omega^X(\Delta(2r,Q))}{r^{n-1}}.$$

Since $|\Delta\psi_r| \leq Cr^{-2}$, (5.34) and (5.35) become

$$\ell(Q)\int_\Omega G(X,Z)\Delta\psi_r(Z)d\,Z - C\left(\frac{r}{Ns}\right)^\alpha \omega^X(\Delta(2r,Q)) \leq \int_\Omega u(Z)\Delta\psi_r(Z)d\,Z,$$

and

$$\int_\Omega u(Z)\Delta\psi_r(Z)d\,Z \leq \ell(Q)\int_\Omega G(X,Z)\Delta\psi_r(Z)d\,Z + C\left(\frac{r}{Ns}\right)^\alpha \omega^X(\Delta(2r,Q)).$$

By the maximum principle we have that $u_r(X) \geq \omega^X(\Delta(r,Q))$; this together with Lemma 3.1 guarantees that $u_r(X) \geq C^{-1}\omega^X(\Delta(2r,Q))$. Hence,

$$\ell(Q) - C\left(\frac{r}{Ns}\right)^\alpha \leq \frac{\int_{\partial\Omega}\psi_r(Q)d\,\omega(Q)}{u_r(X)} \leq \ell(Q) + C\left(\frac{r}{Ns}\right)^\alpha.$$

Taking limits as $r \to 0$, we conclude that

$$\lim_{r\to 0}\frac{\int_{\partial\Omega}\psi_r(Q)d\,\omega(Q)}{u_r(X)} = \ell(Q) = \lim_{Z\to Q}\frac{u(Z)}{G(X,Z)},$$

which is (5.31). When $X \in \Omega\backslash B(2Ns,Q_0)$ for $Q, Q' \in \Delta(2s,Q_0)$, Theorem 3.2 combined with (5.31) asserts that

$$|K(X,Q) - K(X,Q')| \leq C\frac{G(X,A(s,Q_0))}{u(A(s,Q_0))}\left(\frac{|Q-Q'|}{s}\right)^\alpha. \tag{5.36}$$

Lemma 3.5 guarantees that for $Z \in B(2s,Q_0)\cap\Omega$,

$$C^{-1}\frac{G(X,A(s,Q_0))}{u(A(s,Q_0))} \leq \frac{G(X,Z)}{u(Z)} \leq C\frac{G(X,A(s,Q_0))}{u(A(s,Q_0))},$$

where $C > 1$ depends only on the NTA constants. Therefore since the function $\frac{G(X,-)}{u(-)}$ is Hölder continuous in $\overline{\Omega}\cap\overline{B}(2s,Q_0)$, letting $Z \to Q$ we conclude that

$$\frac{G(X,A(s,Q_0))}{u(A(s,Q_0))} \sim K(X,Q). \tag{5.37}$$

Combining (5.36) and (5.37) we obtain

$$|K(X,Q) - K(X,Q')| \leq CK(X,Q)\left(\frac{|Q-Q'|}{s}\right)^\alpha. \qquad \square$$

From now on the constant depends on the compact set $K$ that was chosen in Remark 5.1. In fact, recall that $Q_0 \in K\cap\partial\Omega$. The reader should also keep in mind that the point $A$ always depends on $\mathcal{K}$. Thus in statements the dependence on $\mathcal{K}$ is hidden in $A$.

Corollary 5.4. *Let $\Omega \subset \mathbf{R}^{n+1}$ be a $(\delta_1,\infty)$-chord arc domain and a Reifenberg flat domain with vanishing constant. Given $\varepsilon > 0$, $M > 1$, and $\mathcal{K} \geq 4$ there exist $\theta(\varepsilon) > 0$ and $s(M,\mathcal{K}) > 0$ such that, for $\theta \in (0,\theta(\varepsilon))$ and $s \in (0, s(M,\mathcal{K}))$, if $Q \in \Delta(Ms,Q_0)$ and $Q' \in \Delta(\theta s,Q)$ then*

$$\frac{K(A,Q')}{K(A,Q)} \underset{\varepsilon}{\sim} 1.$$

*Proof.* Choose $\theta' > 0$ so that $4N_0\theta' = 1$, where $N_0$ is as in Lemma 5.4. If $Q \in \Delta(Ms, Q_0)$ and $s(M, \mathcal{K}) > 0$ is chosen so that (5.12) holds then $|A - Q| > \frac{1}{2}s = 2N_0\theta' s$. Lemma 5.4 guarantees that for $\theta \le \theta'$, and $Q' \in \Delta(\theta s, Q)$,

$$|\frac{K(A,Q')}{K(A,Q)} - 1| \le C\left(\frac{|Q-Q'|}{\theta' s}\right)^\alpha \le C\left(\frac{\theta}{\theta'}\right)^\alpha.$$

Choose $\theta(\varepsilon) > 0$ so that $1 + C\left(\frac{\theta}{\theta'}\right)^\alpha \le 1 + \varepsilon$ and $1 - C\left(\frac{\theta}{\theta'}\right)^\alpha \ge \frac{1}{1+\varepsilon}$. We conclude then that $\frac{1}{1+\varepsilon} \le \frac{K(A,Q')}{K(A,Q)} \le 1 + \varepsilon$. □

Corollary 5.5. *Let $\Omega \subset \mathbf{R}^{n+1}$ be a $(\delta_1, \infty)$-chord arc domain and a Reifenberg flat domain with vanishing constant. Given $\varepsilon > 0$, $M > 1$, and $\mathcal{K} \ge 4$ there exist $\theta(\varepsilon) > 0$ and $s(M, \mathcal{K}) > 0$ such that for $\theta \in (0, \theta(\varepsilon))$ and $s \in (0, s(M, \mathcal{K}))$ if $Q \in \Delta(Ms, Q_0)$ then*

$$\frac{\omega^A(\Delta(\theta s, Q))}{\omega(\Delta(\theta s, Q))} \underset{\varepsilon}{\sim} \frac{k_A(Q)}{h(Q)}.$$

*Proof.* Let $\varepsilon > 0$, and let $\theta(\varepsilon) > 0$ be chosen as in the statement of Corollary 5.4. Let $\theta \in (0, \theta(\varepsilon))$. If $Q \in \Delta(Ms, Q_0)$ then for every $Q' \in \Delta(\theta s, Q)$

$$\frac{1}{1+\varepsilon} K(A, Q') \le K(A, Q) \le (1 + \varepsilon) K(A, Q).$$

Multiplying by $h(Q')$ and integrating with respect to $d\sigma$ over $\Delta(\theta s, Q)$ we have that

$$\begin{aligned} \frac{1}{1+\varepsilon} \int_{\Delta(\theta s, Q)} k_A(Q') d\sigma(Q') &\le \frac{k_A(Q)}{h(Q)} \omega(\Delta(\theta s, Q)) \\ &\le (1+\varepsilon) \int_{\Delta(\theta s, Q)} k_A(Q') d\sigma(Q'). \end{aligned}$$ □

Corollary 5.6. *Let $\Omega \subset \mathbf{R}^{n+1}$ be a $(\delta_1, \infty)$-chord arc domain and a Reifenberg flat domain with vanishing constant. Given $\varepsilon > 0$, and $M > 1$, there exists $\theta(\varepsilon, M) > 0$ such that for $\theta \in (0, \theta(\varepsilon, M))$, there exists $\mathcal{K}(\varepsilon, M, \theta) \ge 4$, so that for $\mathcal{K} \ge \mathcal{K}(\varepsilon, M, \theta)$, there exists $s(\varepsilon, M, \theta, \mathcal{K}) > 0$, so that for $s \in (0, s(\varepsilon, M, \theta, \mathcal{K}))$, if $Q \in \Delta(Ms, Q_0)$ then,*

$$\frac{\omega_0^A(\Delta_+(\theta s, \Pi(Q)))}{\omega(\Delta(\theta s, Q))} \underset{\varepsilon}{\sim} \frac{k_A(Q)}{h(Q)}.$$

*Here $\Pi$ denotes the orthogonal projection from $\mathbf{R}^{n+1}$ onto $L(M\mathcal{K}s, Q_0)$, and $\omega_0$ the harmonic measure of the half space containing $A$ and with boundary $L(M\mathcal{K}s, Q_0)$.*

*Proof.* Let $\varepsilon' = \varepsilon'(\varepsilon) > 0$ to be determined. Choose $\theta(\varepsilon, M) = \min\{\theta(\varepsilon', M), \theta(\varepsilon')\}$, with $\theta(\varepsilon', M)$ as in Theorem 5.5, and $\theta(\varepsilon')$ as in Corollary 5.5. Then for $\theta \in (0, \theta(\varepsilon, M))$ there exists $\mathcal{K}(\varepsilon', M, \theta) > 4$ such that if $\mathcal{K} \geq \mathcal{K}(\varepsilon', M, \theta)$ there exists $s(\varepsilon', M, \theta, \mathcal{K}) > 0$ so that for every $s \in (0, s(\varepsilon', M, \theta, \mathcal{K}))$ and every $Q \in \Delta(Ms, Q_0)$,

$$\frac{\omega^A(\Delta(\theta s, Q))}{\omega_0^A(\Delta_+(\theta s, \Pi(Q)))} \underset{\varepsilon'}{\sim} 1, \tag{5.38}$$

and

$$\frac{\omega^A(\Delta(\theta s, Q))}{\omega(\Delta(\theta s, Q))} \underset{\varepsilon'}{\sim} \frac{k_A(Q)}{h(Q)}. \tag{5.39}$$

Combining (5.38) and (5.39) and choosing $\varepsilon' > 0$ so that $(1+\varepsilon')^2 \leq 1 + \varepsilon$ we have that

$$(1+\varepsilon)^{-1}\frac{k_A(Q)}{h(Q)} \leq \frac{\omega_0^A(\Delta_+(\theta s, \Pi(Q)))}{\omega(\Delta(\theta s, Q))} \leq (1+\varepsilon)\frac{k_A(Q)}{h(Q)}. \qquad \square$$

*Remark* 5.2. Recall that

$$\omega_0^A(\Delta_+(\theta s, \Pi(Q))) = \int_{\Delta_+(\theta s, \Pi(Q))} P_A(x)dx,$$

where $P_A$ denotes the Poisson kernel of $\{(x,t) \in \mathbf{R}^{n+1} : x \in L(M\mathcal{K}s, Q_0), t \geq 0\}$ with pole at $A$. Lemma 5.1 guarantees that, given $\varepsilon > 0$ and $M > 1$, there exists $\theta_1(\varepsilon, M) \in (0,1)$ such that for every $s > 0$, $\mathcal{K} > 4$, $Q \in \Delta(Ms, Q_0)$, and $r, r_1 \in (0, \theta_1 s)$,

$$\frac{\omega_0^A(\Delta_+(r_1, \Pi(Q)))}{\omega_n r_1^n} \underset{\varepsilon}{\sim} \frac{\omega_0^A(\Delta_+(r, \Pi(Q)))}{\omega_n r^n}.$$

Letting $r_1 \to 0$ we conclude that, given $\varepsilon > 0$ and $M > 1$, there exists $\theta_1(\varepsilon, M) \in (0,1)$ such that, for any $s > 0$ and any $\mathcal{K} > 4$, and for every $Q \in \Delta(Ms, Q_0)$ and every $r \in (0, \theta_1 s)$,

$$P_A(\Pi(Q)) \underset{\varepsilon}{\sim} \frac{\omega_0^A(\Delta_+(r, \Pi(Q)))}{\omega_n r^n}.$$

Corollary 5.7. *Let $\Omega \subset \mathbf{R}^{n+1}$ be a $(\delta_1, \infty)$-chord arc domain and a Reifenberg flat domain with vanishing constant. Given $\varepsilon > 0$ and $M > 1$, there exists $\theta(\varepsilon, M) > 0$ such that, for $\theta \in (0, \theta(\varepsilon, M))$, there exists $\mathcal{K}(\varepsilon, M, \theta) \geq 4$ so that, for $\mathcal{K} \geq \mathcal{K}(\varepsilon, M, \theta)$ there exists $s(\varepsilon, M, \theta, \mathcal{K}) > 0$, such that for $s \in (0, s(\varepsilon, M, \theta, \mathcal{K}))$ and $Q \in \Delta(Ms, Q_0)$*

$$\frac{k_A(Q)}{h(Q)} \underset{\varepsilon}{\sim} \frac{\omega_n(\theta s)^n}{\omega(\Delta(\theta s, Q))} P_A(\Pi(Q)).$$

*Proof.* Let $\varepsilon' = \varepsilon'(\varepsilon) > 0$ to be determined. Choose $\theta(\varepsilon, M) = \min\{\theta(\varepsilon', M), \theta(\varepsilon')\}$, with $\theta(\varepsilon', M)$ as in Corollary 5.6, and $\theta(\varepsilon')$ as in Remark 5.2 above. Then for $\theta \in (0, \theta(\varepsilon, M))$ there exists $\mathcal{K}(\varepsilon', M, \theta) > 4$ such that, if $\mathcal{K} \geq \mathcal{K}(\varepsilon', M, \theta)$, there exists $s(\varepsilon', M, \theta, \mathcal{K}) > 0$ so that, for every $s \in (0, s(\varepsilon', M, \theta, \mathcal{K}))$, and $Q \in \Delta(Ms, Q_0)$

$$\frac{\omega_0^A(\Delta_+(\theta s, \Pi(Q)))}{\omega(\Delta(\theta s, Q))} \underset{\varepsilon'}{\sim} \frac{k_A(Q)}{h(Q)}, \tag{5.40}$$

and

$$P_A(\Pi(Q)) \underset{\varepsilon'}{\sim} \frac{\omega_0^A(\Delta_+(\theta s, \Pi(Q)))}{\omega_n(\theta s)^n}. \tag{5.41}$$

Combining (5.40) and (5.41), and choosing $\varepsilon' > 0$ so that $(1+\varepsilon')^2 \leq 1 + \varepsilon$ we have

$$(1+\varepsilon)^{-1} \frac{k_A(Q)}{h(Q)} \leq \frac{\omega_n(\theta s)^n}{\omega(\Delta(\theta s, Q))} P_A(\Pi(Q)) \leq (1+\varepsilon) \frac{k_A(Q)}{h(Q)}. \qquad \square$$

Theorem 3.1, and similar arguments to the ones presented above guarantee that:

COROLLARY 5.8. *Let $\Omega \subset \mathbf{R}^{n+1}$ be a Reifenberg flat domain with vanishing constant. There exists $\delta > 0$ such that, if $\Omega$ is a $\delta$-chord arc domain and $X \in \Omega$, then the following statement holds: Given $\varepsilon > 0$ and $M > 1$, there exists $\theta(\varepsilon, M) > 0$ such that, for $\theta \in (0, \theta(\varepsilon, M))$ there exists $\mathcal{K}(\varepsilon, M, \theta) \geq 4$ so that, for $\mathcal{K} \geq \mathcal{K}(\varepsilon, M, \theta)$ there exists $s_1 = s(\varepsilon, M, \theta, \mathcal{K}, \operatorname{dist}(X, \partial\Omega)) > 0$ such that, for $s \in (0, s_1)$ and $Q \in \Delta(Ms, Q_0)$,*

$$\frac{k_A(Q)}{k_X(Q)} \underset{\varepsilon}{\sim} \frac{\omega_n(\theta s)^n}{\omega^X(\Delta(\theta s, Q))} P_A(\Pi(Q)).$$

*Remark* 5.3. Recall that under the assumption that $\Omega$ is Reifenberg flat, the fact that $\omega$ is asymptotically optimally doubling is equivalent to the fact that $\Omega$ is a Reifenberg flat domain with vanishing constant . It is an easy consequence of the fact that $\omega$ is an asymptotically optimally doubling measure that for $\varepsilon > 0$, $M > 1$ and $\theta \in (0,1)$, there exists $s = s(\varepsilon, M, \theta) > 0$ such that for every $Q \in \Delta(Ms, Q_0)$

$$\frac{1}{1+\varepsilon} \left(\frac{\theta}{M}\right)^n \leq \frac{\omega(\Delta(\theta s, Q))}{\omega(\Delta(Ms, Q))} \leq (1+\varepsilon) \left(\frac{\theta}{M}\right)^n.$$

LEMMA 5.5. *Let $\Omega \subset \mathbf{R}^{n+1}$ be a Reifenberg flat domain. Assume that $\omega$ is asymptotically optimally doubling. Given $\varepsilon > 0$ and $M > 1$, there exists $s(\varepsilon, M) > 0$ so that, for $s \in (0, s(\varepsilon, M))$, and $Q \in \Delta(Ms, Q_0)$,*

$$\frac{\omega(\Delta(Ms, Q_0))}{\omega(\Delta(Ms, Q))} \underset{\varepsilon}{\sim} 1. \tag{5.42}$$

Corollary 5.9. *Let $\Omega \subset \mathbf{R}^{n+1}$ be a Reifenberg flat domain. Assume that $\omega$ is asymptotically optimally doubling. Given $\varepsilon > 0$, $M > 1$, and $\theta \in (0,1)$, there exists $s(\varepsilon, M, \theta) > 0$ so that, for $s \in (0, s(\varepsilon, M, \theta))$, and $Q \in \Delta(Ms, Q_0)$*

$$\frac{\omega(\Delta(\theta s, Q))}{\omega(\Delta(Ms, Q_0))} \underset{\varepsilon}{\sim} \left(\frac{\theta}{M}\right)^n .$$

*Proof.* Let $\varepsilon' = \varepsilon'(\varepsilon) > 0$ to be determined. For $M > 1$ and $\theta \in (0,1)$, there exists $s(\varepsilon', M, \theta) > 0$ so that for $s \in (0, s(\varepsilon', M, \theta))$, and $Q \in \Delta(Ms, Q_0)$, Remark 5.3 and Lemma 5.5 hold. Namely for $Q \in \Delta(Ms, Q_0)$,

$$\left(\frac{\theta}{M}\right)^n \underset{\varepsilon'}{\sim} \frac{\omega(\Delta(\theta s, Q))}{\omega(\Delta(Ms, Q))}, \tag{5.43}$$

and

$$\frac{\omega(\Delta(Ms, Q))}{\omega(\Delta(Ms, Q_0))} \underset{\varepsilon'}{\sim} 1. \tag{5.44}$$

Combining (5.43) and (5.44) and choosing $\varepsilon' > 0$ so that $(1+\varepsilon')^2 \le 1 + \varepsilon$ we have that

$$(1+\varepsilon)^{-1} \frac{\omega(\Delta(\theta s, Q))}{\omega(\Delta(Ms, Q_0))} \le \left(\frac{\theta}{M}\right)^n \le (1+\varepsilon) \frac{\omega(\Delta(\theta s, Q))}{\omega(\Delta(Ms, Q_0))}. \qquad \square$$

*Proof of Lemma* 5.5. Let $\varepsilon' = \varepsilon'(\varepsilon) > 0$ and $N = N(\varepsilon') > 2$, to be determined. Assume that $Ms < 1$. Let $\tau_1 = (N+1)^{-1}$, $\tau_2 = (N-1)^{-1}$, and $\tau_3 = N^{-1}$. Since $\omega$ is asymptotically optimally doubling, there exists $R(\varepsilon', \tau_1, \tau_2, \tau_3) > 0$ so that, for $0 < r \le R(\varepsilon', \tau_1, \tau_2, \tau_3)$, and $Q \in \Delta(Ms, Q_0)$,

$$\frac{\omega(\Delta(\tau_i r, Q))}{\omega(\Delta(r, Q))} \underset{\varepsilon'}{\sim} \tau_i^n, \tag{5.45}$$

for $i = 1, 2, 3$. When $s \le \frac{R(\varepsilon', \tau_1, \tau_2, \tau_3)}{2MN}$, then

$$\begin{aligned} \frac{\omega(\Delta(Ms, Q))}{\omega(\Delta(Ms, Q_0))} &= \frac{\omega(\Delta(Ms, Q))}{\omega(\Delta(MNs, Q))} \cdot \frac{\omega(\Delta(MNs, Q))}{\omega(\Delta(Ms, Q_0))}, \\ \frac{\omega(\Delta(Ms, Q))}{\omega(\Delta(Ms, Q_0))} &\sim_{\varepsilon'} \frac{1}{N^n} \cdot \frac{\omega(\Delta(MNs, Q))}{\omega(\Delta(Ms, Q_0))}. \end{aligned} \tag{5.46}$$

Moreover, for $Q \in \Delta(Ms, Q_0)$,

$$\frac{\omega(\Delta(M(N-1)s, Q_0))}{\omega(\Delta(Ms, Q_0))} \le \frac{\omega(\Delta(MNs, Q))}{\omega(\Delta(Ms, Q_0))} \le \frac{\omega(\Delta(M(N+1)s, Q_0))}{\omega(\Delta(Ms, Q_0))}. \tag{5.47}$$

Applying (5.45) to (5.47), we obtain

$$(1+\varepsilon')^{-1}(N-1)^n \le \frac{\omega(\Delta(MNs, Q))}{\omega(\Delta(Ms, Q_0))} \le (1+\varepsilon')(N+1)^n. \tag{5.48}$$

Combinination of (5.46) and (5.48) yields

$$(1+\varepsilon')^{-2}\left(1-\frac{1}{N}\right)^n \le \frac{\omega(\Delta(Ms,Q))}{\omega(\Delta(Ms,Q_0))} \le (1+\varepsilon')^2\left(1+\frac{1}{N}\right)^n.$$

Choosing $N > 2$ so that $(1-\frac{1}{N})^n \ge (1+\varepsilon')^{-1}$ and $(1+\frac{1}{N})^n \le 1+\varepsilon'$, and $\varepsilon' > 0$ so that $(1+\varepsilon')^3 \le 1+\varepsilon$ we conclude that for $s \le s(\varepsilon, M) = \frac{R(\varepsilon',\tau_1,\tau_2,\tau_3)}{2MN}$, and $Q \in \Delta(Ms, Q_0)$, (5.42) is satisfied. □

Note that Lemma 5.5 and Corollary 5.9 are results about asymptotically optimally doubling measures, therefore similar results hold for the harmonic measure with finite pole, under the hypotheses of Theorems 5.3 and 5.4.

Corollary 5.10. *Let $\Omega \subset \mathbf{R}^{n+1}$ be a $(\delta_1, \infty)$-chord arc domain. Assume that $\omega$ is asymptotically optimally doubling. Given $\varepsilon > 0$ and $M > 1$, there exists $\mathcal{K}(\varepsilon, M) > 4$ such that, for $\mathcal{K} \ge \mathcal{K}(\varepsilon, M)$ there is $s(\varepsilon, M, \mathcal{K}) > 0$ so that, for $s \in (0, s(\varepsilon, M, \mathcal{K}))$ if $Q \in \Delta(Ms, Q_0)$, then*

$$\frac{k_A(Q)}{h(Q)} \underset{\varepsilon}{\sim} \frac{\omega_n (Ms)^n}{\omega(\Delta(Ms,Q_0))} P_A(\Pi(Q)).$$

Corollary 5.11. *Let $\Omega \subset \mathbf{R}^{n+1}$ be a Reifenberg flat domain. There exists $\delta > 0$ such that if $\Omega$ is a $\delta$-chord arc domain, $X \in \Omega$ and $\omega^X$ is asymptotically optimally doubling then given $\varepsilon > 0$ and $M > 1$, there exists $\mathcal{K}(\varepsilon, M) > 4$ such that, for $\mathcal{K} \ge \mathcal{K}(\varepsilon, M)$ there is $s_1 = s(\varepsilon, M, \mathcal{K}, \operatorname{dist}(X, \partial\Omega)) > 0$ so that, for $s \in (0, s_1)$ if $Q \in \Delta(Ms, Q_0)$, then*

$$\frac{k_A(Q)}{k_X(Q)} \underset{\varepsilon}{\sim} \frac{\omega_n (Ms)^n}{\omega^X(\Delta(Ms,Q_0))} P_A(\Pi(Q)).$$

*Proof of Corollary* 5.10. Under the hypotheses above, $\Omega$ is a Reifenberg flat domain with vanishing constant . Let $\varepsilon' = \varepsilon'(\varepsilon) > 0$ to be determined later. Corollary 5.7 guarantees that there exist $\theta = \theta(\varepsilon', M) > 0$, $\mathcal{K} = \mathcal{K}(\varepsilon', M, \theta(\varepsilon', M)) = \mathcal{K}(\varepsilon', M) \ge 4$, and $s_1(\varepsilon', M) > 0$ so that, for $s \in (0, s_1(\varepsilon', M))$, if $Q \in \Delta(Ms, Q_0)$, then

$$\frac{k_A(Q)}{h(Q)} \underset{\varepsilon'}{\sim} \frac{\omega_n (\theta s)^n}{\omega(\Delta(\theta s,Q))} P_A(\Pi(Q)). \tag{5.49}$$

Corollary 5.9 guarantees that there exists $s(\varepsilon', M, \theta) = s_2(\varepsilon', M) > 0$ so that for $s \in (0, s_2(\varepsilon', M))$ and $Q \in \Delta(Ms, Q_0)$,

$$\frac{\omega(\Delta(\theta s,Q))}{\omega(\Delta(Ms,Q_0))} \underset{\varepsilon'}{\sim} \left(\frac{\theta}{M}\right)^n. \tag{5.50}$$

For $s \le \min\{s_1(\varepsilon', M), s_2(\varepsilon', M)\}$, combining (5.49) and (5.50) we obtain

$$\begin{aligned}(1+\varepsilon')^{-2}\frac{\omega_n(Ms)^n}{\omega(\Delta(Ms,Q_0))}P_A(\Pi(Q)) &\le \frac{k_A(Q)}{h(Q)}\\ \frac{k_A(Q)}{h(Q)} &\le (1+\varepsilon')^2\frac{\omega_n(Ms)^n}{\omega(\Delta(Ms,Q_0))}P_A(\Pi(Q)).\end{aligned}$$

Choosing $\varepsilon' > 0$ so that $(1+\varepsilon')^2 \le 1+\varepsilon$ we finish the proof of Corollary 5.10. □

So far we have not used the hypothesis that $\log h \in \mathrm{VMO}(d\sigma)$ ($\log k_X \in \mathrm{VMO}(d\sigma)$). The fact that $\log h \in \mathrm{VMO}(d\sigma)$ guarantees that $h(Q)$ can be well approximated by its average over a small ball around $Q$, at least for a large set of $Q$'s. More precisely:

LEMMA 5.6. *Let $\Omega$ be a $(\delta_1,\infty)$-chord arc domain. Assume that $\log h \in \mathrm{VMO}(d\sigma)$. Given $\varepsilon > 0$, there exists $r(\varepsilon) > 0$ such that, for every $r \in (0, r(\varepsilon))$, there exists $G(r,Q_0) \subset \Delta(r,Q_0)$ such that $\sigma(\Delta(r,Q_0)) \le (1+\varepsilon)\sigma(G(r,Q_0))$, and for all $P \in G(r,Q_0)$*

$$h(P) \underset{\varepsilon}{\sim} \frac{\omega(\Delta(r,Q_0))}{\sigma(\Delta(r,Q_0))}.$$

*Proof.* Let $\varepsilon' = \varepsilon'(\varepsilon) \in (0,1)$ to be determined later. Since $\log h \in \mathrm{VMO}(\partial\Omega)$, there exists $r(\varepsilon') = r_1 > 0$ so that

$$\sup_{Q_0\in\partial\Omega\cap K}\ \sup_{0<r\le r_1} \frac{1}{\sigma(\Delta(r,Q_0))}\int_{\Delta(r,Q_0)} |\log h - (\log h)_{r,Q_0}|\,d\sigma < \varepsilon'$$

where

$$(\log h)_{r,Q_0} = \frac{1}{\sigma(\Delta(r,Q_0))}\int_{\Delta(r,Q_0)} \log h\,d\sigma.$$

For $r \in (0, r_1]$ let

$$G(r,Q_0) = \{P \in \Delta(r,Q_0) : |\log h(P) - (\log h)_{r,Q_0}| < \sqrt{\varepsilon'}\}.$$

Then

$$\sigma(G(r,Q_0)) \ge (1-\sqrt{\varepsilon'})\sigma(\Delta(r,Q_0)), \tag{5.51}$$

and if $P \in G(r,Q_0)$

$$|\log h(P) - (\log h)_{r,Q_0}| \le \sqrt{\varepsilon'}.$$

Thus

$$e^{-\sqrt{\varepsilon'}}e^{(\log h)_{r,Q_0}} \le h(P) \le e^{\sqrt{\varepsilon'}}e^{(\log h)_{r,Q_0}}. \tag{5.52}$$

Integrating (5.52) over $G(r,Q_0)$ and using (5.51) we obtain

$$(1-\sqrt{\varepsilon'})e^{-\sqrt{\varepsilon'}}e^{(\log h)_{r,Q_0}} \le \frac{1}{\sigma(\Delta(r,Q_0))}\int_{G(r,Q_0)} h\,d\sigma \le e^{\sqrt{\varepsilon'}}e^{(\log h)_{r,Q_0}}. \tag{5.53}$$

We now need to compare $\frac{1}{\sigma(\Delta(r,Q_0))}\int_{G(r,Q_0)} h\,d\sigma$ with $\frac{1}{\sigma(\Delta(r,Q_0))}\int_{\Delta(r,Q_0)} h\,d\sigma$. Clearly,

$$\frac{1}{\sigma(\Delta(r,Q_0))}\int_{G(r,Q_0)} h\,d\sigma \le \frac{1}{\sigma(\Delta(r,Q_0))}\int_{\Delta(r,Q_0)} h\,d\sigma,$$

and therefore

$$(1-\sqrt{\varepsilon'})e^{-\sqrt{\varepsilon'}}e^{(\log h)_{r,Q_0}} \le \frac{1}{\sigma(\Delta(r,Q_0))}\int_{\Delta(r,Q_0)} h\,d\sigma. \tag{5.54}$$

Note that

$$\frac{1}{\sigma(\Delta(r,Q_0))}\int_{G(r,Q_0)} h\,d\sigma = \frac{\omega(\Delta(r,Q_0))}{\sigma(\Delta(r,Q_0))} - \frac{1}{\sigma(\Delta(r,Q_0))}\int_{\Delta(r,Q_0)\backslash G(r,Q_0)} h\,d\sigma. \tag{5.55}$$

Since $\omega \in A_\infty(d\sigma)$ (see remark after (5.10) and [GCRdF]), there exist $C > 1$ and $\gamma \in (0,1)$, so that

$$\begin{aligned} \frac{\omega(\Delta(r,Q_0)\backslash G(r,Q_0))}{\omega(\Delta(r,Q_0))} &\le C\left(\frac{\sigma(\Delta(r,Q_0)\backslash G(r,Q_0))}{\sigma(\Delta(r,Q_0))}\right)^\gamma, \\ \frac{\omega(\Delta(r,Q_0)\backslash G(r,Q_0))}{\omega(\Delta(r,Q_0))} &\le C(\sqrt{\varepsilon'})^\gamma. \end{aligned} \tag{5.56}$$

Combining (5.55) and (5.56) we obtain

$$(1-C(\sqrt{\varepsilon'})^\gamma)\frac{1}{\sigma(\Delta(r,Q_0))}\int_{\Delta(r,Q_0)} h\,d\sigma \le \frac{1}{\sigma(\Delta(r,Q_0))}\int_{G(r,Q_0)} h\,d\sigma. \tag{5.57}$$

From (5.53) and (5.57) we deduce that

$$\frac{1}{\sigma(\Delta(r,Q_0))}\int_{\Delta(r,Q_0)} h\,d\sigma \le (1-C(\sqrt{\varepsilon'})^\gamma)^{-1}e^{\sqrt{\varepsilon'}}e^{(\log h)_{r,Q_0}}. \tag{5.58}$$

From (5.52), (5.54), and (5.58) we conclude that for $P \in G(r,Q_0)$

$$e^{-2\sqrt{\varepsilon'}}(1-C(\sqrt{\varepsilon'})^\gamma)\frac{1}{\sigma(\Delta(r,Q_0))}\int_{\Delta(r,Q_0)} h\,d\sigma \le h(P),$$

and

$$h(P) \le e^{2\sqrt{\varepsilon'}}(1-\sqrt{\varepsilon'})^{-1}\frac{1}{\sigma(\Delta(r,Q_0))}\int_{\Delta(r,Q_0)} h\,d\sigma.$$

Choose $\varepsilon'(\varepsilon) > 0$ small enough so that $e^{-2\sqrt{\varepsilon'}}(1-C(\sqrt{\varepsilon'})^\gamma) \ge (1+\varepsilon)^{-1}$ and $e^{2\sqrt{\varepsilon'}}(1-\sqrt{\varepsilon'})^{-1} \le 1+\varepsilon$. Under these assumptions we have proved that there exists $r(\varepsilon) > 0$ so that for $r \in (0, r(\varepsilon))$ and $Q_0 \in \partial\Omega \cap K$ the conclusion of Lemma 5.6 holds. □

Note that Lemma 5.6 is a result about functions whose logarithm is in VMO. Therefore an analogous statement holds for the Poisson kernel with finite pole, $k_X$, provided that $\log k_X \in \mathrm{VMO}(d\sigma)$. In particular the following corollaries also hold under the assumptions of Theorems 5.3 and 5.4.

The following statement is a straightforward consequence of Corollary 5.10 and Lemma 5.6. In Corollaries 5.12 and 5.13, $s$ is replaced by $2s$ for notational convenience.

COROLLARY 5.12. *Let $\Omega$ be a $(\delta_1, \infty)$-chord arc domain. Assume that $\omega$ is asymptotically optimally doubling, and that $\log h \in \mathrm{VMO}(d\sigma)$. Given $\varepsilon > 0$ and $M > 1$, there exists $\mathcal{K}(\varepsilon, M) > 4$ so that for $\mathcal{K} \geq \mathcal{K}(\varepsilon, M)$, there is $s(\varepsilon, M, \mathcal{K}) > 0$ such that, for every $s \in (0, s(\varepsilon, M, \mathcal{K}))$ there exists $G(2Ms, Q_0) \subset \Delta(2Ms, Q_0)$ such that $\sigma(\Delta(2Ms, Q_0)) \leq (1+\varepsilon)\sigma(G(2Ms, Q_0))$ and for all $Q \in G(2Ms, Q_0)$,*

$$k_A(Q) \underset{\varepsilon}{\sim} \frac{\omega_n (2Ms)^n}{\sigma(\Delta(2Ms, Q_0))} P_A(\Pi(Q)).$$

The next corollary states that a chord arc domain whose harmonic measure is asymptotically optimally doubling, and whose Poisson kernel has logarithm in VMO has the following property: each surface ball in the boundary can be decomposed into a large very flat piece where both the unit normal and the Poisson kernel have small oscillation, and a very small additional piece.

COROLLARY 5.13. *Let $\Omega \subset \mathbf{R}^{n+1}$ be a $(\delta_1, \infty)$-chord arc domain. There exists $\delta_0 \in (0, \delta_n)$ such that, if $\Omega$ is a $\delta$-Semmes decomposable domain for some $\delta \in (0, \delta_0)$, if $\omega$ is asymptotically optimally doubling, and if $\log h \in \mathrm{VMO}(d\sigma)$, then given $\varepsilon > 0$ and $M > 1$, there exists $\mathcal{K}(\varepsilon, M) > 4$ such that, for $\mathcal{K} \geq \mathcal{K}(\varepsilon, M)$ there is $s(\varepsilon, M, \mathcal{K}) > 0$ so that, for $s \in (0, s(\varepsilon, M, \mathcal{K}))$, there exist an $n$-dimensional plane $L(M\mathcal{K}s, Q_0)$ containing $Q_0$ and a Lipschitz function $\phi : L(M\mathcal{K}s, Q_0) \to \mathbf{R}$ such that $\|\nabla\phi\|_\infty \leq C(n)(\delta + \varepsilon)$. The graph of this function $\mathcal{G} = \{(x,t) \in \mathbf{R}^{n+1} : t = \phi(x)\}$ approximates $\partial\Omega$ in the ball $B(2Ms, Q_0)$, in the sense that*

$$\Delta(2Ms, Q_0) = \mathcal{A}(Ms, Q_0) \cup \mathcal{F}(Ms, Q_0), \tag{5.59}$$

*where*

$$\mathcal{A}(Ms, Q_0) \subset \mathcal{G} \quad \text{and} \quad \sigma(\mathcal{F}(Ms, Q_0)) \leq (C_1 \exp(-\frac{C_2}{\delta}) + \varepsilon)\sigma(\Delta(2Ms, Q_0)), \tag{5.60}$$

*for some $C_1, C_2 > 0$. Moreover if $\Pi : \mathbf{R}^{n+1} \to L(M\mathcal{K}s, Q_0)$ denotes the orthogonal projection*

$$\sup_{\Pi(\mathcal{A}(Ms,Q_0))} |\phi| \leq 2\varepsilon Ms \quad \text{and} \quad \sup_{B(2Ms,Q_0)\cap L(M\mathcal{K}s,Q_0)} |\phi| \leq C(n)(\varepsilon + \delta)Ms. \tag{5.61}$$

*for every $Q \in \mathcal{A}(Ms, Q_0)$,*

$$k_A(Q) \underset{\varepsilon}{\sim} \frac{\omega_n (2Ms)^n}{\sigma(\Delta(2Ms, Q_0))} P_A(\Pi(Q)). \tag{5.62}$$

*Proof of Corollary* 5.13. Let $\delta_0 \in (0, \delta_n)$ be as in Lemma 4.1, and assume that $\Omega$ is a $\delta$-Semmes decomposable domain for $\delta < \delta_0$. Let $\varepsilon > 0$ and $M > 1$. Corollary 5.12 guarantees that there exists $\mathcal{K}(\varepsilon, M) > 4$ such that for $\mathcal{K} > \mathcal{K}(\varepsilon, M)$ there is $s(\varepsilon, M, \mathcal{K}) > 0$ such that, for every $s \in (0, s(\varepsilon, M, \mathcal{K}))$ and $Q_0 \in K \cap \partial\Omega$ there exists $G(2Ms, Q_0) \subset \Delta(2Ms, Q_0)$ such that $\sigma(\Delta(2Ms, Q_0)) \leq (1+\varepsilon)\sigma(G(2Ms, Q_0))$ and for all $Q \in G(2Ms, Q_0)$

$$k_A(Q) \underset{\varepsilon}{\sim} \frac{\omega_n(2Ms)^n}{\sigma(\Delta(2Ms, Q_0))} P_A(\Pi(Q)).$$

Since $\omega$ is asymptotically optimally doubling, we know that $\Omega$ is a Reifenberg flat domain with vanishing constant . There exists $r(\varepsilon, \mathcal{K}) > 0$ so that for $r \in (0, r(\varepsilon, \mathcal{K}))$

$$\frac{1}{r} D[\partial\Omega \cap B(r, Q_0), L(r, Q_0) \cap B(r, Q_0)] < \frac{\varepsilon}{\mathcal{K}}.$$

Thus, if $s \in (0, \frac{r(\varepsilon,\mathcal{K})}{M\mathcal{K}})$, it is easy to check that (see Remark 5.1)

$$\frac{1}{2Ms} D[\partial\Omega \cap B(2Ms, Q_0), L(M\mathcal{K}s, Q_0) \cap B(2Ms, Q_0)] < \varepsilon.$$

Lemma 4.1 guarantees that there exists a Lipschitz function $\phi : L(M\mathcal{K}s, Q_0) \to \mathbf{R}$ such that $\|\nabla\phi\|_\infty \leq C(n)(\delta+\varepsilon)$ whose graph $\mathcal{G} = \{(x,t) \in \mathbf{R}^{n+1} : t = \phi(x)\}$ approximates $\partial\Omega$ in the ball $B(2Ms, Q_0)$, in the sense that

$$\Delta(2Ms, Q_0) = \mathcal{G}(2Ms, Q_0) \cup \mathcal{E}(2Ms, Q_0),$$

where

$$\mathcal{G}(2Ms, Q_0) \subset \mathcal{G} \quad \text{and} \quad \sigma(\mathcal{E}(2Ms, Q_0)) \leq C_1 \exp(-\frac{C_2}{\delta})\sigma(\Delta(2Ms, Q_0)),$$

for some $C_1, C_2 > 0$. Moreover

$$\sup_{\Pi(\mathcal{G}(2Ms, Q_0))} |\phi| \leq 2\varepsilon Ms.$$

Let $0 < s < \min\{\frac{r(\varepsilon,\mathcal{K})}{M\mathcal{K}}, s(\varepsilon, M, \mathcal{K})\}$. Let $\mathcal{A}(Ms, Q_0) = G(2Ms, Q_0) \cap \mathcal{G}(2Ms, Q_0)$ and $\mathcal{F}(Ms, Q_0) = (\Delta(2Ms, Q_0) \backslash \mathcal{A}(Ms, Q_0)$. It is easy to check that (5.59), (5.60), (5.61) and (5.62) are satisfied. □

## 6. Rellich's identity

For $n \geq 2$, let $\Omega \subset \mathbf{R}^{n+1}$ be a $\delta$-chord arc domain, and $A \in \Omega$. Rellich's identity (see [JK3] for the case of a graph, see [KT2] for the general case) asserts that

$$(*) \qquad \int_{\partial\Omega} k_A(Q) \frac{d\sigma(Q)}{|Q-A|^{n-1}} = \frac{1}{\sigma_n} \int_{\partial\Omega} k_A^2(Q) \langle Q - A, \overrightarrow{n}(Q)\rangle d\sigma,$$

where $\sigma_n$ denotes the surface area of the unit sphere in $\mathbf{R}^{n+1}$. Rellich's identity applied to the half space $\{(x,t) \in \mathbf{R}^{n+1} : x \in L, t \geq 0\}$ containing $A$ yields

$$\int_L P_A(x) \frac{dx}{|x-A|^{n-1}} = \frac{1}{\sigma_n} \int_L s P_A^2(x) dx.$$

We first truncate these integrals so that we are only concerned with what is happening in $\Delta(Ms, Q_0)$ and $B(Ms, Q_0) \cap L$. Integration by parts in these regions allows us to capitalize on the fact that Corollary 5.13 holds in $\Delta(Ms, Q_0)$. This provides a good estimate for the oscillation of the unit normal to $\partial\Omega$ on $\Delta(2s, Q_0)$.

Lemma 6.1. *Let $\Omega \subset \mathbf{R}^{n+1}$ be a $(\delta_1, \infty)$-chord arc domain. Assume that $\Omega$ is a Reifenberg flat domain with vanishing constant. Let $K \subset \mathbf{R}^{n+1}$ be a compact set. Given $\eta > 0$ there exists $M(\eta) > 1$ such that for $M \geq M(\eta)$, $s > 0$, $Q_0 \in K \cap \partial\Omega$, if $A \in \Omega$, $\frac{s}{2} \leq d(A, Q_0) \leq s$, and $d(\partial\Omega, A) \geq \frac{s}{16}$ then*

$$\int_{\partial\Omega \cap \{|Q-Q_0| \geq Ms\}} k_A^2(Q) |\langle Q - A, \vec{n}(Q)\rangle| d\sigma < \eta s^{-(n-1)}, \tag{6.1}$$

*and*

$$\int_{\partial\Omega \cap \{|Q-Q_0| \geq Ms\}} \frac{k_A(Q)}{|Q-A|^{n-1}} d\sigma < \eta s^{-(n-1)}. \tag{6.2}$$

Let now $A = Q_0 + s\vec{n}(M\mathcal{K}s, Q_0)$, and $\vec{n}(M\mathcal{K}s, Q_0)$ be the unit normal vector to $L(M\mathcal{K}s, Q_0)$. Note that $A \in \{(x,t) \in \mathbf{R}^{n+1} : x \in L(M\mathcal{K}s, Q_0), t \geq 0\}$, $|A - Q_0| = s$, and $d(A, L(M\mathcal{K}s, Q_0)) = s$. Applying Lemma 6.1 to $\{(x,t) = x + t\vec{n}(M\mathcal{K}s, Q_0) \in \mathbf{R}^{n+1} : x \in L(M\mathcal{K}s, Q_0), t \geq 0\}$, we obtain:

Corollary 6.1. *Let $P_A$ denote the Poisson kernel of the half space $\{(x,t) \in \mathbf{R}^{n+1} : x \in L(M\mathcal{K}s, Q_0), t \geq 0\}$ with pole at $A$. Given $\eta > 0$ there exists $M(\eta) > 1$ such that for $M \geq M(\eta)$*

$$s \int_{L(M\mathcal{K}s, Q_0) \cap \{|x-Q_0| \geq Ms\}} P_A^2(x) dx \leq \eta s^{-(n-1)},$$

*and*

$$\int_{L(M\mathcal{K}s, Q_0) \cap \{|x-Q_0| \geq Ms\}} \frac{P_A(x)}{|x-A|^{n-1}} dx \leq \eta s^{-(n-1)}.$$

*Proof of Lemma* 6.1. Let $M > 4$. We first estimate

$$\begin{aligned} &\int_{\partial\Omega \cap \{|Q-Q_0| \geq Ms\}} k_A^2(Q) |\langle Q - A, \vec{n}(Q)\rangle| d\sigma \\ &\qquad = \sum_{i=0}^{\infty} \int_{2^i Ms \leq |Q-Q_0| \leq 2^{i+1} Ms} k_A^2(Q) |\langle Q - A, \vec{n}(Q)\rangle| d\sigma \end{aligned} \tag{6.3}$$

$$\leq \sum_{i=0}^{\infty} \int_{2^i Ms \leq |Q-Q_0| \leq 2^{i+1}Ms} k_A^2(Q)(|Q-Q_0| + |Q_0 - A|)d\sigma$$

$$\leq 2\sum_{i=0}^{\infty} 2^{i+1}Ms \int_{2^i Ms \leq |Q-Q_0| \leq 2^{i+1}Ms} k_A^2(Q)d\sigma.$$

We look at each term $\int_{2^i Ms \leq |Q-Q_0| \leq 2^{i+1}Ms} k_A^2(Q)d\sigma$ separately. First note that if $Q \in \Delta(2^{i+1}Ms, Q_0)\backslash\Delta(2^i Ms, Q_0)$, $|Q-A| \geq |Q-Q_0| - |Q_0 - A| > 2^{i-1}Ms$. Cover $\Delta(2^{i+1}Ms, Q_0)\backslash\Delta(2^i Ms, Q_0)$ by balls $\Delta(\rho_i, Q_j)$, with

$$Q_j \in \Delta(2^{i+1}Ms, Q_0)\backslash\Delta(2^i Ms, Q_0)$$

and such that the balls $\Delta(\frac{\rho_i}{5}, Q_j)$ are disjoint. Assume that $\rho_i > 0$ is such that $N\rho_i = 2^{i-1}Ms$, where $N = 2N_0 > 2$, and $N_0$ is as in Corollary 5.2. Note that $A \in \Omega\backslash B(N\rho_i, Q_j)$. Corollary 5.2 guarantees that for each $Q_j$

$$\left(\frac{1}{\sigma(\Delta(\rho_i, Q_j))}\int_{\Delta(\rho_i, Q_j)} k_A^2 d\sigma\right)^{\frac{1}{2}} \leq 2\frac{1}{\sigma(\Delta(\rho_i, Q_j))}\int_{\Delta(\rho_i, Q_j)} k_A d\sigma.$$

Recall that, since $\Omega$ is a $(\delta_1, \infty)$-chord arc domain, there exists a constant $C(n) > 1$ depending only on $n$ such that $\sigma(\Delta(\rho_i, Q_j)) \geq C(n)^{-1}\rho_i^n$. Moreover, the fact that $\Omega$ is an unbounded NTA domain, with uniform constants, guarantees that $\omega^A$ is uniformly doubling on $\partial\Omega \cap \{|Q-Q_0| \geq Ms\}$. Therefore the previous inequality implies that

(6.4)
$$\begin{aligned}
\int_{2^i Ms \leq |Q-Q_0| \leq 2^{i+1}Ms} k_A^2(Q)d\sigma &\leq \sum_j \int_{\Delta(\rho_i, Q_j)} k_A^2 d\sigma \\
&\leq 4\sum_j \frac{\omega^A(\Delta(\rho_i, Q_j))}{\sigma(\Delta(\rho_i, Q_j))}\omega^A(\Delta(\rho_i, Q_j)) \\
&\leq C(n)\rho_i^{-n}\sum_j \omega^A(\Delta(\rho_i, Q_j)) \\
&\leq C(n)\rho_i^{-n}\sum_j \omega^A(\Delta(\frac{\rho_i}{5}, Q_j)) \\
&\leq C(n)\rho_i^{-n}\omega^A\left(\Delta(2^{i+1}Ms + \frac{\rho_i}{5}, Q_0)\backslash\Delta(2^i Ms - \frac{\rho_i}{5}, Q_0)\right).
\end{aligned}$$

In particular

(6.5)
$$\begin{aligned}
&\int_{2^i Ms \leq |Q-Q_0| \leq 2^{i+1}Ms} k_A^2(Q)d\sigma \\
&\qquad \leq C(n)\rho_i^{-n}\omega^A(\Delta(2^{i+1}Ms + \frac{2^{i-1}Ms}{N}, Q_0)\backslash\Delta(2^i Ms - \frac{2^{i-1}Ms}{N}, Q_0)).
\end{aligned}$$

Note that $\omega^X(\Delta(2^{i+1}Ms+\frac{2^{i-1}Ms}{N},Q_0)\backslash\Delta(2^iMs-\frac{2^{i-1}Ms}{N},Q_0))$ is a nonnegative harmonic function in $\Omega$, which vanishes on $B(2^iMs-\frac{2^{i-1}Ms}{N},Q_0)\cap\partial\Omega$, and whose supremum is 1. Thus Lemma 3.1 implies that

$$\begin{aligned}(6.6)\qquad \omega^A(\Delta(2^{i+1}Ms+\frac{2^{i-1}Ms}{N},Q_0)\ &\backslash\Delta(2^iMs-\frac{2^{i-1}Ms}{N},Q_0))\\ &\leq C\left(\frac{|A-Q_0|}{2^{i-1}Ms-\frac{2^{i-2}Ms}{N}}\right)^\alpha\\ &\leq C\left(\frac{1}{2^iM}\right)^\alpha.\end{aligned}$$

Combining (6.5) and (6.6) we obtain

$$(6.7)\qquad \int_{2^iMs\leq|Q-Q_0|\leq 2^{i+1}Ms} k_A^2(Q)d\sigma\leq C(n)\rho_i^{-n}\left(\frac{1}{2^iM}\right)^\alpha.$$

Thus (6.3) and (6.7) yield

$$\begin{aligned}&\int_{\partial\Omega\cap\{|Q-Q_0|\geq Ms\}} k_A^2(Q)|\langle Q-A,\overrightarrow{n}(Q)\rangle|d\sigma\\ &\qquad\leq \sum_{i=0}^\infty C(n)2^{i+1}Ms\rho_i^{-n}\left(\frac{1}{2^iM}\right)^\alpha\\ &\qquad\leq \sum_{i=0}^\infty C(n)2^{i+1}Ms\left(\frac{N}{2^{i-1}Ms}\right)^n\left(\frac{1}{2^iM}\right)^\alpha\\ &\qquad\leq C_{n,N_0}\frac{s^{-(n-1)}}{M^{n-1+\alpha}}\sum_{i=0}^\infty\frac{1}{2^{i(n-1+\alpha)}},\end{aligned}$$

and

$$(6.8)\qquad \int_{\partial\Omega\cap\{|Q-Q_0|\geq Ms\}} k_A^2(Q)|\langle Q-A,\overrightarrow{n}(Q)\rangle|d\sigma\leq C(n,N_0)\frac{1}{M^{n-1+\alpha}}s^{-(n-1)}.$$

To estimate the second integral we use the fact that

$$\omega^X(\Delta(2^{i+1}Ms,Q_0)\backslash\Delta(2^iMs,Q_0))$$

is a nonnegative harmonic function in $\Omega$. It is bounded by 1, and vanishes continuously on $B(2^iMs,Q_0)\cap\partial\Omega$. In this case $\omega^A(\Delta(2^{i+1}Ms,Q_0)\backslash\Delta(2^iMs,Q_0))$ $\leq C\frac{1}{(2^iM)^\alpha}$. We then have

$$\begin{aligned}(6.9)\quad &\int_{\partial\Omega\cap\{|Q-Q_0|\geq Ms\}}\frac{k_A(Q)}{|Q-A|^{n-1}}d\sigma\\ &\qquad= \sum_{i=0}^\infty\int_{2^iMs\leq|Q-Q_0|\leq 2^{i+1}Ms}\frac{k_A(Q)}{|Q-A|^{n-1}}d\sigma\\ &\qquad\leq 2^{n-1}\sum_{i=0}^\infty\frac{1}{(2^iMs)^{n-1}}\int_{2^iMs\leq|Q-Q_0|\leq 2^{i+1}Ms}k_A(Q)d\sigma\end{aligned}$$

$$\leq C(n)\sum_{i=0}^{\infty}\frac{1}{(2^iMs)^{n-1}}\cdot\frac{1}{(2^iM)^{\alpha}}$$
$$\leq C(n)\frac{1}{M^{n-1+\alpha}}s^{-(n-1)}.$$

Choosing $M > 4$ so that $C(n)\frac{1}{M^{n-1+\alpha}} < \eta$ and $C(n,N_0)\frac{1}{M^{n-1+\alpha}} < \eta$ we conclude from (6.8) and (6.9) that (6.1) and (6.2) hold. □

Since the proof of Lemma 6.1 relies on the fact that the domain is unbounded we sketch the proof of the corresponding result for bounded domains.

LEMMA 6.2. *There exists $\delta_1 > 0$ so that if $\Omega$ is a bounded $\delta_1$-chord arc domain which is a Reifenberg flat domain vanishing constant, then given $\eta > 0$ there exists $M(\eta) > 1$ such that, for $M \geq M(\eta)$ there is $s(M) > 0$ so that, for $s \in (0, s(M))$, $Q_0 \in \partial\Omega$, if $A \in \Omega$, $\frac{s}{2} \leq d(A, Q_0) \leq s$, and $d(\partial\Omega, A) \geq \frac{s}{16}$ then*

$$\int_{\partial\Omega\cap\{|Q-Q_0|\geq Ms\}} k_A^2(Q)|\langle Q-A, \overrightarrow{n}(Q)\rangle|d\sigma < \eta s^{-(n-1)}, \tag{6.10}$$

*and*

$$\int_{\partial\Omega\cap\{|Q-Q_0|\geq Ms\}} \frac{k_A(Q)}{|Q-A|^{n-1}}d\sigma < \eta s^{-(n-1)}. \tag{6.11}$$

*Proof.* Let $M > 4$ to be chosen. Corollary 5.2 guarantees that there exists $\delta_1$ so that, if $\Omega$ is a bounded $\delta_1$-chord arc domain, there exist $N_0 > 1$ and $s_0 > 0$ such that, for $N = 2N_0$ and $\rho \in (0, s_0)$, if $Q \in \partial\Omega$ and $A \in \Omega\backslash B(N\rho, Q)$ then (5.4) holds with $A$ in place of $X$, and $\rho$ in place of $s$. For $s < s_0/M$, let $i_0 \in \mathbf{N}$ be such that $Ms \leq 2^{i_0}Ms \leq s_0 \leq 2^{i_0+1}Ms$. Since $\Omega$ is a $\delta_1$-chord arc domain we may assume that $s_0$ was chosen so that for $s \leq s_0$, and $Q \in \partial\Omega$, $\sigma(\Delta(s,Q)) \geq C_n s^n$. To prove the first inequality, we estimate the following expression

$$\begin{aligned}\int_{\partial\Omega\cap\{|Q-Q_0|\geq Ms\}} & k_A^2(Q)|\langle Q-A, \overrightarrow{n}(Q)\rangle|d\sigma \\ \leq\ & 2\sum_{i=0}^{i_0-1}2^{i+1}Ms\int_{2^iMs\leq|Q-Q_0|\leq 2^{i+1}Ms}k_A^2(Q)d\sigma \\ & + \int_{|Q-Q_0|\geq 2^{i_0}Ms}|Q-A|k_A^2(Q)d\sigma.\end{aligned} \tag{6.12}$$

The same argument as in the proof of Lemma 6.1 allows us to conclude that

$$\begin{aligned}\sum_{i=0}^{i_0-1}2^{i+1}Ms\int_{2^iMs\leq|Q-Q_0|\leq 2^{i+1}Ms}k_A^2(Q)d\sigma &\leq C_{n,N_0}\frac{s^{-(n-1)}}{M^{n-1+\alpha}}\sum_{i=0}^{i_0}\frac{1}{2^{i(n-1+\alpha)}} \\ &\leq C(n,N_0)\frac{1}{M^{n-1+\alpha}}s^{-(n-1)}.\end{aligned} \tag{6.13}$$

We now estimate the second term in (6.12)

$$\int_{|Q-Q_0|\geq 2^{i_0}Ms} k_A^2(Q)|Q-A|d\sigma \leq \operatorname{diam}\Omega \int_{|Q-Q_0|\geq s_0/2} k_A^2(Q)d\sigma. \tag{6.14}$$

Here $\operatorname{diam}\Omega$ denotes the diameter of $\Omega$. Note that if $|Q-Q_0| \geq \frac{s_0}{2}$, then $|Q-A| \geq \frac{s_0}{4}$. Cover $\partial\Omega\backslash\Delta(\frac{s_0}{2}, Q_0)$ by balls $\Delta(\frac{s_0}{4N}, Q_j)$, where $Q_j \in \partial\Omega\backslash\Delta(\frac{s_0}{2}, Q_0)$ such that the balls $\Delta(\frac{s_0}{20N}, Q_j)$ are disjoint. Corollary 5.2 guarantees that for each $Q_j$

$$\left(\frac{1}{\sigma(\Delta(\frac{s_0}{4N}, Q_j))}\int_{\Delta(\frac{s_0}{4N},Q_j)} k_A^2 d\sigma\right)^{\frac{1}{2}} \leq 2\frac{1}{\sigma(\Delta(\frac{s_0}{4N}, Q_j))}\int_{\Delta(\frac{s_0}{4N},Q_j)} k_A d\sigma.$$

Following the same steps as in inequality (6.4) we have

$$\begin{aligned} \int_{|Q-Q_0|\geq\frac{s_0}{2}} k_A^2(Q)d\sigma &\leq \sum_j \int_{\Delta(\frac{s_0}{4N},Q_j)} k_A^2 d\sigma \\ &\leq C(n)\left(\frac{s_0}{4N}\right)^{-n} \sum_j \omega^A(\Delta(\frac{s_0}{20N}, Q_j)) \\ &\leq C(n)\left(\frac{s_0}{4N}\right)^{-n}. \end{aligned} \tag{6.15}$$

From (6.14) and (6.15) we deduce that

$$\begin{aligned} \int_{|Q-Q_0|\geq 2^{i_0}Ms} k_A^2(Q)|Q-A|d\sigma &\leq C(n)\operatorname{diam}\Omega\left(\frac{s_0}{4N}\right)^{-n} \\ &\leq C(n,N_0)\frac{\operatorname{diam}\Omega}{s_0} s_0^{-(n-1)} \\ &\leq C(n,N_0)\frac{\operatorname{diam}\Omega}{s_0}\cdot\frac{s^{-(n-1)}}{M^{n-1}}. \end{aligned} \tag{6.16}$$

Combining (6.13) and (6.16) we conclude that

$$\int_{\partial\Omega\cap\{|Q-Q_0|\geq Ms\}} k_A^2(Q)|\langle Q-A, \overrightarrow{n}(Q)\rangle|d\sigma \leq C(n,N_0,\Omega)\frac{s^{-(n-1)}}{M^{n-1}}.$$

The proof of the second inequality follows the same pattern as the proof of inequality (6.10). □

Corollary 6.2. *Let $\Omega \subset \mathbf{R}^{n+1}$ be a $(\delta_1, \infty)$-chord arc domain. Assume that $\Omega$ is a Reifenberg flat domain with vanishing constant. Let $K \subset \mathbf{R}^{n+1}$ be a compact set. Given $\eta > 0$, there exist $M(\eta) > 1$ such that, for $M \geq M(\eta)$, $s > 0$, $Q_0 \in K \cap \partial\Omega$, and $A \in \Omega$ such that $\frac{s}{2} \leq d(A, Q_0) \leq s$ and $d(A, \partial\Omega) \geq \frac{s}{16}$, if $\psi \in C_c^\infty(\mathbf{R}^{n+1})$ satisfies $\psi \equiv 1$ on $B(Ms, Q_0)$; $\psi \equiv 0$ on $\mathbf{R}^{n+1}\backslash B(2Ms, Q_0)$, $0 \leq \psi \leq 1$, then*

(6.17)

$$\left|\int_{\partial\Omega}\psi(Q)\frac{k_A(Q)d\sigma(Q)}{|Q-A|^{n-1}}-\frac{1}{\sigma_n}\int_{\partial\Omega}k_A^2(Q)\psi(Q)\langle Q-A;\overrightarrow{n}(Q)\rangle d\sigma(Q)\right|\le \eta s^{-(n-1)}.$$

*Here $\sigma_n$ denotes the surface area of the unit sphere in $\mathbf{R}^{n+1}$ and $\overrightarrow{n}(Q)$ denotes the outward unit normal to $\partial\Omega$.*

*Proof.* Let $\eta' = \eta'(\eta) > 0$, to be chosen. Let $M(\eta') > 1$ be so that Lemma 6.1 holds. Let $M \ge M(\eta')$ and $s > 0$; then by $(*)$

$$\begin{aligned}
&\left|\int_{\partial\Omega}k_A\psi\frac{d\sigma}{|Q-A|^{n-1}}-\frac{1}{\sigma_n}\int_{\partial\Omega}\psi k_A^2\langle Q-A,\overrightarrow{n}(Q)\rangle d\sigma\right|\\
&\quad\le \left|\int_{\partial\Omega}k_A\psi\frac{d\sigma}{|Q-A|^{n-1}}-\int_{\partial\Omega}k_A\frac{d\sigma}{|Q-A|^{n-1}}\right|\\
&\qquad+\frac{1}{\sigma_n}\left|\int_{\partial\Omega}k_A^2\langle Q-A;\overrightarrow{n}(Q)\rangle(1-\psi)d\sigma\right|\\
&\quad\le \int_{\partial\Omega\cap\{|Q-Q_0|\ge Ms\}}k_A\frac{d\sigma}{|Q-A|^{n-1}}\\
&\qquad+\frac{1}{\sigma_n}\int_{\partial\Omega\cap\{|Q-Q_0|\ge Ms\}}k_A^2|\langle Q-A,\overrightarrow{n}(Q)\rangle|d\sigma\\
&\quad< \eta' s^{-(n-1)}+\eta' s^{-(n-1)}=2\eta' s^{-(n-1)}.
\end{aligned}$$

Choosing $\eta' = \eta/4$ we conclude that (6.17) holds for $M \ge M(\eta/4)$. □

Note that if $\Omega$ is a bounded $\delta_1$-chord arc domain which is Reifenberg flat with vanishing constant a similar statement holds for $s > 0$ small enough. Namely, $s < s(M)$, where $s(M)$ is chosen as in Lemma 6.2.

Corollary 6.3. *Let $L(M\mathcal{K}s, Q_0)$, $A$ and $P_A$ be as in Corollary* 6.1. *Given $\eta > 0$ there exists $M(\eta) > 1$ such that, for $M \ge M(\eta)$ if $\psi \in C_c^\infty(\mathbf{R}^{n+1})$ satisfies $\psi \equiv 1$ on $B(Ms, Q_0)$, $\psi \equiv 0$ on $\mathbf{R}^{n+1}\backslash B(2Ms, Q_0)$, $0 \le \psi \le 1$, then*

$$\left|\int_{L(M\mathcal{K}s,Q_0)}\psi P_A\frac{dx}{|x-A|^{n-1}}-\frac{1}{\sigma_n}\int_{L(M\mathcal{K}s,Q_0)}s\psi P_A^2dx\right|\le \eta s^{-(n-1)}.$$

Next, as in Jerison's paper (see [J]) we apply integration by parts in order to combine Corollaries 5.13, 6.2, 6.3 and 5.12, to control the $L^2$ mean oscillation of the unit normal vector (see the proof of (6.25)).

Under the hypothesis of Theorems 5.3 and 5.4, Corollaries 6.2 and 6.3 hold. A careful look at the proof of Corollary 6.5 will reveal that only these two results and the bounded domain versions of Corollaries 5.12 and 5.13 are needed to insure that the conclusions of this corollary are valid in the bounded domain setting.

THEOREM 6.1. *Let $\Omega$ be a $(\delta_1,\infty)$-chord arc domain. There exists $\delta_0 \in (0,\delta_n)$ such that, if $\Omega$ is a $\delta$-Semmes decomposable domain with $\delta \in (0,\delta_0)$, if $\omega$ is asymptotically optimally doubling, and if $\log h \in \mathrm{VMO}(d\sigma)$, for each compact set $K \subset \mathbf{R}^{n+1}$ the following statement is true: Given $\varepsilon > 0$ and $M > 1$, there exists $\mathcal{K}(\varepsilon, M) > 4$ such that for $\mathcal{K} \geq \mathcal{K}(\varepsilon, M)$ there is $s(\varepsilon, M, \mathcal{K}) > 0$ such that, for every $s \in (0, s(\varepsilon, M, \mathcal{K}))$ and $Q_0 \in \partial\Omega \cap K$, there exists a Lipschitz function $\phi : L(M\mathcal{K}s, Q_0) \to \mathbf{R}$ such that $\|\nabla\phi\|_\infty \leq C(n)(\delta+\varepsilon)$, and whose graph $\mathcal{G} = \{(x,t) \in \mathbf{R}^{n+1} : t = \phi(x)\}$ approximates $\partial\Omega$ in the ball $B(2Ms, Q_0)$ in the sense that*

$$\Delta(2Ms, Q_0) = \mathcal{A}(Ms, Q_0) \cup \mathcal{F}(Ms, Q_0),$$

*where*

$$\mathcal{A}(Ms, Q_0) \subset \mathcal{G} \quad \text{and} \quad \sigma(\mathcal{F}(Ms,Q_0)) \leq (C_1 \exp(-\frac{C_2}{\delta}) + \varepsilon)\sigma(\Delta(2Ms, Q_0)),$$

*for some $C_1, C_2 > 0$. Moreover if $\Pi : \mathbf{R}^{n+1} \to L(M\mathcal{K}s, Q_0)$ denotes the orthogonal projection,*

$$\sup_{\Pi(\mathcal{A}(Ms,Q_0))} |\phi| \leq 2\varepsilon Ms, \qquad \sup_{B(2Ms,Q_0)\cap L(M\mathcal{K}s,Q_0)} |\phi| \leq C(n)(\varepsilon+\delta)Ms, \tag{6.18}$$

*and for every $Q \in \mathcal{A}(Ms, Q_0)$*

$$k_A(Q) \underset{\varepsilon}{\sim} \frac{\omega_n (2Ms)^n}{\sigma(\Delta(2Ms,Q_0))} P_A(\Pi(Q)).$$

*Moreover if $\psi \in C_c^\infty(B(2Ms,Q_0))$, $0 \leq \psi \leq 1$, and $\psi \equiv 1$ on $B(Ms,Q_0)$, then*

$$\int_{\mathcal{F}(Ms,Q_0)} k_A(Q)\psi(Q) \frac{d\sigma}{|Q-A|^{n-1}}(Q) \leq C(n) M^n (e^{-C/\delta} + \varepsilon) s^{-(n-1)}, \tag{6.19}$$

$$\int_{\mathcal{F}(Ms,Q_0)} k_A^2(Q)\psi(Q) |\langle Q-A, \overrightarrow{n}(Q)\rangle| d\sigma(Q) \leq C(n) M^{n+1} (e^{-C/\delta} + \varepsilon) s^{-(n-1)}, \tag{6.20}$$

$$s \int_{(\Pi(\mathcal{A}(Ms,Q_0)))^c \cap L(M\mathcal{K}s,Q_0)} \psi(x,\phi(x)) P_A^2(x) dx \leq c(n) M^n (e^{-C/\delta} + \varepsilon) s^{-(n-1)}, \tag{6.21}$$

$$\begin{aligned} &\int_{(\Pi(\mathcal{A}(Ms,Q_0)))^c \cap L(M\mathcal{K}s,Q_0)} \psi(x,\phi(x)) \frac{P_A(x)}{|x-A|^{n-1}} \sqrt{1+|\nabla\phi(x)|^2} dx \\ &\qquad \leq c(n) M^n (e^{-C/\delta} + \varepsilon) s^{-(n-1)}, \end{aligned} \tag{6.22}$$

$$\begin{aligned} &\int_{\mathcal{A}(Ms,Q_0)} \psi(Q) |k_A^2(Q) - \alpha^2 P_A(\Pi(Q))| \, |\langle Q-A, \overrightarrow{n}(Q)\rangle| d\sigma(Q) \\ &\qquad \leq c(n) M^{n+1} \varepsilon s^{-(n-1)}, \end{aligned} \tag{6.23}$$

$$(6.24)\qquad \int_{\mathcal{A}(Ms,Q_0)} |k_A(Q) - \alpha P_A(\Pi(Q))| \frac{\psi(Q)}{|Q-A|^{n-1}} d\sigma(Q) \le c(n) M^n \varepsilon s^{-(n-1)},$$

*where* $\alpha = \frac{\omega_n (2Ms)^n}{\sigma(\Delta(2Ms,Q_0))}$, *and*

$$(6.25)\qquad \left| \int_{\mathcal{A}(Ms,Q_0)} \psi(Q) P_A^2(\Pi(Q)) \langle Q - A, \overrightarrow{n}(Q)\rangle d\sigma(Q) \right.$$
$$\left. - s \int_{L(M\mathcal{K}s,Q_0)} \psi(x,\phi(x)) P_A^2(x) dx \right| \le c(n) M^{n+2} (e^{-C/\delta} + \varepsilon) s^{-(n-1)},$$

$$(6.26)\qquad \left| \int_{\mathcal{A}(Ms,Q_0)} P_A(\Pi(Q)) \frac{\psi(Q)}{|Q-A|^{n-1}} d\sigma(Q) \right.$$
$$\left. - \int_{L(M\mathcal{K}s,Q_0)} \psi(x,\phi(x)) \frac{P_A(x)}{|x-A|^{n-1}} \sqrt{1 + |\nabla\phi(x)|^2} dx \right|$$
$$\le c(n) M^{n+1} (e^{-C/\delta} + \varepsilon) s^{-(n-1)}.$$

*Proof of Theorem* 6.1. Initially assume that $\delta \in (0,\delta_0)$, where $\delta_0$ is as in Corollary 5.13. Given $\varepsilon' = \varepsilon'(\varepsilon) \in (0,1)$ and $M > 1$, Corollary 5.13 guarantees that there exists $\mathcal{K}(\varepsilon', M) > 4$ such that, for $\mathcal{K} \ge \mathcal{K}(\varepsilon', M)$ there is $s(\varepsilon', M, \mathcal{K}) > 0$ such that, for $s \in (0, s(\varepsilon', M, \mathcal{K}))$ and $Q_0 \in K \cap \partial\Omega$ there exist an $n$-dimensional plane $L(M\mathcal{K}s, Q_0)$ and a Lipschitz function $\phi : L(M\mathcal{K}s, Q_0) \to \mathbf{R}$ such that $\|\nabla\phi\|_\infty \le C(n)(\delta + \varepsilon')$, and whose graph $\mathcal{G}$ approximates $\partial\Omega$ in the ball $B(2Ms, Q_0)$, in the sense that

$$\Delta(2Ms, Q_0) = \mathcal{A}(Ms, Q_0) \cup \mathcal{F}(Ms, Q_0)$$

where

$$\mathcal{A}(Ms, Q_0) \subset \mathcal{G} \text{ and } \sigma(\mathcal{F}(Ms, Q_0)) \le (C_1 e^{-C_2/\delta} + \varepsilon')\sigma(\Delta(2Ms, Q_0)).$$

Moreover for every $Q \in \mathcal{A}(Ms, Q_0)$

$$k_A(Q) \underset{\varepsilon'}{\sim} \frac{\omega_n (2Ms)^n}{\sigma(\Delta(2Ms, Q_0))} P_A(\Pi(Q)).$$

In particular, provided that we choose $\delta > 0$ and $\varepsilon' > 0$ small, we can insure that $\frac{1}{2} \le \alpha = \frac{\omega_n (2Ms)^n}{\sigma(\Delta(2Ms,Q_0))} \le 2$ (see Theorem 4.2 and Remark 4.1). Note that

$$(6.27)\qquad \int_{\mathcal{F}(Ms,Q_0)} k_A(Q)\psi(Q) \frac{d\sigma(Q)}{|Q-A|^{n-1}} \le \frac{C(n)}{s^{n-1}} \int_{\mathcal{F}(Ms,Q_0)} K(A,Q) h(Q) d\sigma(Q),$$

where $K(A,Q) = \frac{k_A(Q)}{h(Q)}$. Corollary 5.10 guarantees that there exists $\mathcal{K}(\varepsilon', M) > 4$ such that, for $\mathcal{K} \ge K(\varepsilon', M)$ there exists $s(\varepsilon', M, \mathcal{K}) > 0$ so that for

$s \in (0, s(\varepsilon', M, \mathcal{K}))$ if $Q \in \Delta(2Ms, Q_0)$

$$K(A,Q) \underset{\varepsilon'}{\sim} \frac{\omega_n (2Ms)^n}{\omega(\Delta(2Ms,Q_0))} P_A(\Pi(Q)).$$

Moreover for $Q \in \Delta(2Ms, Q_0)$, $P_A(\Pi(Q)) \leq C(n)s^{-n}$. Thus

$$K(A,Q) \leq C(n) \frac{M^n}{\omega(\Delta(2Ms,Q_0))}. \tag{6.28}$$

Therefore (6.27) becomes

$$\int_{\mathcal{F}(Ms,Q_0)} k_A(Q)\psi(Q)\frac{d\sigma(Q)}{|Q-A|^{n-1}} \leq C(n)\frac{M^n}{\omega(\Delta(2Ms,Q_0))} s^{-(n-1)}\omega(\mathcal{F}(Ms,Q_0)). \tag{6.29}$$

Since $\omega \in A_\infty(d\sigma)$, there exists $\gamma \in (0,1)$ such that

$$\frac{\omega(\mathcal{F}(Ms,Q_0))}{\omega(\Delta(2Ms,Q_0))} \leq C\left(\frac{\sigma(\mathcal{F}(Ms,Q_0))}{\sigma(\Delta(2Ms,Q_0))}\right)^\gamma \leq C(e^{-C/\delta}+\varepsilon')^\gamma. \tag{6.30}$$

Combining (6.29) and (6.30) we obtain

$$\int_{\mathcal{F}(Ms,Q_0)} k_A(Q)\psi(Q)\frac{d\sigma(Q)}{|Q-A|^{n-1}} \leq C(n)M^n(e^{-C/\delta}+\varepsilon')^\gamma s^{-(n-1)}.$$

Choosing $\varepsilon' > 0$ so that $(\varepsilon')^\gamma < \varepsilon$, we have that

$$\int_{\mathcal{F}(Ms,Q_0)} k_A(Q)\psi(Q)\frac{d\sigma(Q)}{|Q-A|^{n-1}} \leq C(n)M^n(e^{-C/\delta}+\varepsilon) s^{-(n-1)}.$$

We now look at the expression in (6.20). Using (6.28) we have

$$\begin{aligned} \int_{\mathcal{F}(Ms,Q_0)} & \psi(Q)k_A^2(Q)|\langle Q-A, \overrightarrow{n}(Q)\rangle| d\sigma(Q) \\ &\leq C(n)Ms\int_{\mathcal{F}(Ms,Q_0)} (K(A,Q)h(Q))^2 d\sigma(Q) \\ &\leq C(n)Ms\left(\frac{M^n}{\omega(\Delta(2Ms,Q_0))}\right)^2 \int_{\mathcal{F}(Ms,Q_0)} h^2 d\sigma, \end{aligned} \tag{6.31}$$

since $\log h \in \mathrm{VMO}(d\sigma)$, $h^2 d\sigma \in A_\infty(d\sigma)$ (see [GCRdF]), and there exists $\gamma' > 0$ such that

$$\begin{aligned} \int_{\mathcal{F}(Ms,Q_0)} h^2 d\sigma &\leq C\left(\int_{\Delta(2Ms,Q_0)} h^2 d\sigma\right)\left(\frac{\sigma(\mathcal{F}(Ms,Q_0))}{\sigma(\Delta(2Ms,Q_0))}\right)^{\gamma'} \\ &\leq C(e^{-C/\delta}+\varepsilon')^{\gamma'}\int_{\Delta(2Ms,Q_0)} h^2 d\sigma. \end{aligned} \tag{6.32}$$

Moreover our choice of $\delta_1 > 0$ insures that that

(6.33)

$$\int_{\Delta(2Ms,Q_0)} h^2(Q)d\sigma(Q) \le 4\sigma(\Delta(2Ms,Q_0))\left(\frac{1}{\sigma(\Delta(2Ms,Q_0))}\int_{\Delta(2Ms,Q_0)} hd\sigma\right)^2$$

(see note after Corollary 5.2). Combining (6.31), (6.32) and (6.33) we have

$$\begin{aligned}&\int_{\mathcal{F}(2Ms,Q_0)} \psi(Q)k_A^2(Q)|\langle Q-A,\overrightarrow{n}(Q)\rangle|d\sigma(Q)\\ &\qquad\le\ C(n)(e^{-C/\delta}+\varepsilon')^{\gamma'}Ms\frac{M^{2n}}{\sigma(\Delta(2Ms,Q_0))}\\ &\qquad\le\ C(n)(e^{-C/\delta}+\varepsilon')^{\gamma'}M^{n+1}s^{-(n-1)}.\end{aligned}$$

Choosing $\varepsilon' > 0$ so that $(\varepsilon')^{\gamma'} < \varepsilon$, we finish the proof of (6.20). In order to prove that (6.21) and (6.22) hold, we first need to estimate $|(\Pi(\mathcal{A}(Ms,Q_0)))^c \cap B(2Ms,Q_0)\cap L(M\mathcal{K}s,Q_0)|$. Note that

$$\begin{aligned}&(\Pi(\mathcal{A}(Ms,Q_0)))^c\cap B(2Ms,Q_0)\cap L(M\mathcal{K}s,Q_0)\\ &\qquad= B(2Ms,Q_0)\cap L(M\mathcal{K}s,Q_0)\backslash\Pi(\mathcal{A}(Ms,Q_0))\\ &\qquad= B(2Ms,Q_0)\cap L(M\mathcal{K}s,Q_0)\backslash\Pi(\Delta(2Ms,Q_0))\\ &\qquad\quad\cup\ \Pi(\Delta(2Ms,Q_0))\backslash\Pi(\mathcal{A}(Ms,Q_0)).\end{aligned}$$

Therefore, since $\Pi(\Delta(2Ms,Q_0)) \subset B(2Ms,Q_0)\cap L(M\mathcal{K}s,Q_0)$,

$$\begin{aligned}&|\,(\Pi(\mathcal{A}(Ms,Q_0)))^c\cap B(2Ms,Q_0)\cap L(M\mathcal{K}s,Q_0)\\ &\qquad\le\omega_n(2Ms)^n-|\Pi(\Delta(2Ms,Q_0))|\\ &\qquad\quad+\ |\Pi(\Delta(2Ms,Q_0))\backslash\Pi(\mathcal{A}(Ms,Q_0))|.\end{aligned}\tag{6.34}$$

Recall that under our hypothesis, $\Omega$ is a Reifenberg flat domain with vanishing constant. Therefore, there exists $s_1(\varepsilon',M,\mathcal{K}) > 0$ such that $s_1(\varepsilon',M,\mathcal{K}) < s(\varepsilon',M,\mathcal{K})$, and for $s\in(0,s_1(\varepsilon',M,\mathcal{K}))$

$$\Gamma(2Ms(1-\varepsilon')^{\frac{1}{n}},Q_0)\subset\Delta(2Ms,Q_0)\subset\Gamma(2Ms,Q_0),$$

where $\Gamma(2Ms,Q_0)=\{(x,t): x\in L(M\mathcal{K}s,Q_0),\ |x-Q_0|\le 2Ms,\ |t|\le 2Ms\}$ $\cap\partial\Omega$ (see Remark 4.1). Thus, as in Remark 4.1, the inclusions above guarantee that

$$\omega_n(2Ms)^n(1-\varepsilon')\le|\Pi(\Delta(2Ms,Q_0))|\le\omega_n(2Ms)^n.\tag{6.35}$$

Combining (6.34) and (6.35) we have

$$\begin{aligned}&|\,(\Pi(\mathcal{A}(Ms,Q_0)))^c\cap B(2Ms,Q_0)\cap L(M\mathcal{K}s,Q_0)|\\ &\qquad\le\varepsilon'\omega_n(2Ms)^n+|\Pi(\mathcal{F}(Ms,Q_0))|\end{aligned}\tag{6.36}$$

$$\begin{aligned}&\leq \varepsilon'\omega_n(2Ms)^n+\sigma(\mathcal{F}(Ms,Q_0))\\&\leq \varepsilon'\omega_n(2Ms)^n+(c_1e^{-c_2/\delta}+\varepsilon')\sigma(\Delta(2Ms,Q_0))\\&\leq C(n)(\varepsilon'+Ce^{-C/\delta})(Ms)^n.\end{aligned}$$

We now look at

$$\begin{aligned}s\int_{(\Pi(\mathcal{A}(Ms,Q_0)))^c\cap L(M\mathcal{K}s,Q_0)}&\psi(x,\phi(x))P_A^2(x)dx\\&\leq C(n)s^{-2n+1}|\,(\Pi(\mathcal{A}(Ms,Q_0)))^c\cap B(2Ms,Q_0)\cap L(M\mathcal{K}s,Q_0)|\\&\leq C(n)s^{-2n+1}(\varepsilon'+Ce^{-C/\delta})M^ns^n\\&\leq C(n)M^n(\varepsilon'+e^{-C/\delta})s^{-(n-1)},\end{aligned}$$

and

$$\begin{aligned}\int_{(\Pi(\mathcal{A}(Ms,Q_0)))^c\cap L(M\mathcal{K}s,Q_0)}&\psi(x,\phi(x))\frac{P_A(x)}{|x-A|^{n-1}}\sqrt{1+|\nabla\phi|^2(x)}dx\\&\leq C(n)s^{-n}\int_{(\Pi(\mathcal{A}(Ms,Q_0)))^c\cap L(M\mathcal{K}s,Q_0)\cap B(2Ms,Q_0)}\frac{dx}{|x-A|^{n-1}}\\&\leq C(n)s^{-n}s^{-n+1}|\,(\Pi(\mathcal{A}(Ms,Q_0)))^c\cap B(2Ms,Q_0)\cap L(M\mathcal{K}s,Q_0)|\\&\leq C(n)M^n(e^{-C/\delta}+\varepsilon')s^{-(n-1)}.\end{aligned}$$

This finishes the proof of inequalities (6.21) and (6.22). We now look at the expressions in inequalities (6.23) and (6.24). Corollary 5.13 guarantees that

$$\begin{aligned}\int_{\mathcal{A}(Ms,Q_0)}\psi(Q)|k_A^2(Q)-\alpha^2P_A^2(\Pi(Q))|\,|\langle Q-A,\overrightarrow{n}(Q)\rangle|d\sigma(Q)&\\\leq 4\varepsilon'\alpha^2\int_{\mathcal{A}(Ms,Q_0)}P_A^2(\Pi(Q))|Q-A|d\sigma(Q)&\\\leq C(n)\varepsilon'Mss^{-2n}\alpha^2|\Pi(\mathcal{A}(Ms,Q_0))|&\\\leq C(n)M^{n+1}\varepsilon's^{-(n-1)},&\end{aligned}$$

and

$$\begin{aligned}\int_{\mathcal{A}(Ms,Q_0)}\psi(Q)|k_A(Q)-\alpha P_A(\Pi(Q))|\frac{d\sigma(Q)}{|Q-A|^{n-1}}&\\\leq\varepsilon'\alpha\int_{\mathcal{A}(Ms,Q_0)}P_A(\Pi(Q))\frac{d\sigma(Q)}{|Q-A|^{n-1}}&\\\leq C(n)\varepsilon'\alpha s^{-n}s^{-(n-1)}(Ms)^n&\\\leq C(n)M^n\varepsilon's^{-(n-1)}.&\end{aligned}$$

In order to prove (6.25) we use an integration by parts introduced by Jerison in [J]. Recall that $A=Q_0+s\overrightarrow{n}(M\mathcal{K}s,Q_0)$, where $\overrightarrow{n}(M\mathcal{K}s,Q_0)$ is the unit normal to $L(M\mathcal{K}s, Q_0)$. If $Q \in \mathcal{A}(Ms, Q_0)$ then $Q = x + \phi(x)\overrightarrow{n}(M\mathcal{K}s, Q_0)$ where $x \in L(M\mathcal{K}s, Q_0)$. Moreover,

$$\overrightarrow{n}(Q) = \left(1 + |\nabla\phi(x)|^2\right)^{-\frac{1}{2}}(\nabla\phi(x), -1) \quad \text{and} \quad d\sigma(Q) = \sqrt{1 + |\nabla\phi(x)|^2}dx.$$

We now look at

(6.37)
$$\begin{aligned}
&\int_{\mathcal{A}(Ms,Q_0)} \psi(Q)P_A^2(\Pi(Q))\langle Q - A, \overrightarrow{n}(Q)\rangle d\sigma(Q) \\
&\quad = -\int_{\Pi(\mathcal{A}(Ms,Q_0))} \psi(x,\phi(x))P_A^2(x)\langle (x,\phi(x)) - (Q_0,s); (-\nabla\phi(x),1)\rangle dx \\
&\quad = -\int_{\Pi(\mathcal{A}(Ms,Q_0)} \psi(x,\phi(x))P_A^2(x)(-\langle x - Q_0, \nabla\phi(x)\rangle + \phi(x) - s)dx \\
&\quad = s\int_{L(M\mathcal{K}s,Q_0)} \psi(x,\phi(x))P_A^2(x)dx \\
&\qquad - s\int_{(\Pi(\mathcal{A}(Ms,Q_0)))^c\cap L(M\mathcal{K}s,Q_0)} \psi(x,\phi(x))P_A^2(x)dx \\
&\qquad - \int_{\Pi(\mathcal{A}(Ms,Q_0))} \psi(x,\phi(x))P_A^2(x)(\phi(x) - \langle x - Q_0, \nabla\phi(x)\rangle)dx.
\end{aligned}$$

Note that (6.21) takes care of the second term. We look carefully at the last term. Recall that $\sup_{\Pi(\mathcal{A}(Ms,Q_0))}|\phi| \leq 2\varepsilon' Ms$, and that $\sup_{B(2Ms,Q_0)\cap L(M\mathcal{K}s,Q_0)}|\phi| \leq C(\varepsilon' + \delta)Ms$. Thus

$$\left|\int_{\Pi(\mathcal{A}(Ms,Q_0))} \psi(x,\phi(x))P_A^2(x)\phi(x)dx\right| \leq C(n)M^{n+1}\varepsilon' s^{-(n-1)}, \tag{6.38}$$

and

(6.39)
$$\begin{aligned}
&\int_{\Pi(\mathcal{A}(Ms,Q_0))} \psi(x,\phi(x))P_A^2(x)\langle x - Q_0, \nabla\phi(x)\rangle dx \\
&\quad = \int_{L(M\mathcal{K}s,Q_0)} \psi(x,\phi(x))P_A^2(x)\langle x - Q_0, \nabla\phi(x)\rangle dx \\
&\qquad - \int_{(\Pi(\mathcal{A}(Ms,Q_0)))^c\cap L(M\mathcal{K}s,Q_0)} \psi(x,\phi(x))P_A^2(x)\langle x - Q_0, \nabla\phi(x)\rangle dx.
\end{aligned}$$

Note that (6.21) guarantees that

(6.40)
$$\begin{aligned}
&\left|\int_{(\Pi(\mathcal{A}(Ms,Q_0)))^c\cap L(M\mathcal{K}s,Q_0)} \psi(x,\phi(x))P_A^2(x)\langle x - Q_0, \nabla\phi(x)\rangle dx\right| \\
&\qquad \leq C(n)Ms\int_{(\Pi(\mathcal{A}(Ms,Q_0)))^c\cap L(M\mathcal{K}s,Q_0)} \psi(x,\phi(x))P_A^2(x)dx \\
&\qquad \leq C(n)M^{n+1}(e^{-C/\delta} + \varepsilon)s^{-(n-1)}.
\end{aligned}$$

We now estimate the remaining term

$$\begin{aligned}(6.41)\qquad &\int_{L(M\mathcal{K}s,Q_0)} \psi(x,\phi(x))P_A^2(x)\langle x-Q_0,\nabla\phi(x)\rangle dx \\ &= \int_{L(M\mathcal{K}s,Q_0)} \operatorname{div}(\phi(x)\psi(x,\phi(x))P_A^2(x)(x-Q_0))dx \\ &\quad -\int_{L(M\mathcal{K}s,Q_0)} \phi(x)\psi(x,\phi(x))P_A^2(x)\operatorname{div}(x-Q_0)dx \\ &\quad -\int_{L(M\mathcal{K}s,Q_0)} \phi(x)\langle x-Q_0,\nabla P_A^2(x)\rangle\psi(x,\phi(x))dx \\ &\quad -\int_{L(M\mathcal{K}s,Q_0)} \phi(x)P_A^2(x)\langle x-Q_0,\nabla\psi(x,\phi(x))\rangle dx.\end{aligned}$$

Since $\psi\in C_c^\infty(B(2Ms,Q_0))$ the first term is 0. Using (6.18) and (6.36) we can control the second term as follows:

$$\begin{aligned}(6.42)\quad &\left|\int_{L(M\mathcal{K}s,Q_0)} \phi(x)\psi(x,\phi(x))P_A^2(x)\operatorname{div}(x-Q_0)dx\right| \\ &\le (n+1)\left|\int_{L(M\mathcal{K}s,Q_0)} \phi(x)P_A^2(x)\psi(x,\phi(x))\right| \\ &\le (n+1)\left|\int_{\Pi(\mathcal{A}(Ms,Q_0))} \phi(x)P_A^2(x)\psi(x,\phi(x))\right| \\ &\quad + C(n)Ms\int_{(\Pi(\mathcal{A}(Ms,Q_0)))^c\cap L(M\mathcal{K}s,Q_0)} P_A^2(x)\psi(x,\phi(x))dx \\ &\le C(n)\varepsilon' Ms\int_{\Pi(\mathcal{A}(Ms,Q_0))} P_A^2(x)\psi(x,\phi(x)) \\ &\quad + C(n)M^{n+1}(e^{-C/\delta}+\varepsilon')s^{-(n-1)} \\ &\le C(n)M^{n+1}(e^{-C/\delta}+\varepsilon')s^{-(n-1)}.\end{aligned}$$

Combining (6.18), (6.36) and the fact that $|\nabla P_A(x)|\le C(n)s^{-n-1}$ we have

$$\begin{aligned}&(6.43)\\ &\left|\int_{L(M\mathcal{K}s,Q_0)} \phi(x)2P_A(x)\langle x-Q_0,\nabla P_A(x)\rangle\psi(x,\phi(x))dx\right| \\ &\quad\le\ 2\left|\int_{\Pi(\mathcal{A}(Ms,Q_0))} \phi(x)P_A(x)\langle x-Q_0,\nabla P_A(x)\rangle\psi(x,\phi(x))dx\right| \\ &\qquad +\ C(n)M^2s^2\int_{(\Pi(\mathcal{A}(Ms,Q_0)))^c\cap L(M\mathcal{K}s,Q_0)} P_A(x)|\nabla P_A(x)|\psi(x,\phi(x))dx\end{aligned}$$

$$\leq C(n)\varepsilon' M^2 s^2 \int_{\Pi(\mathcal{A}(Ms,Q_0))} P_A(x)|\nabla P_A(x)|\psi(x,\phi(x))dx$$
$$+ C(n)M^2 s^2 s^{-n} s^{-n-1}|\,(\Pi(\mathcal{A}(Ms,Q_0)))^c \cap B(2Ms,Q_0)\cap L(M\mathcal{K}s,Q_0)|$$
$$\leq C(n)M^{n+2}(e^{-C/\delta}+\varepsilon')s^{-(n-1)}.$$

Combining (6.18), (6.36) and the fact that $|\nabla\psi| \leq Cs^{-1}$ we obtain

(6.44)
$$\left|\int_{L(M\mathcal{K}s,Q_0)} \phi(x)P_A^2(x)\langle x-Q_0, \nabla\psi(x,\phi(x))\rangle dx\right|$$
$$\leq \left|\int_{\Pi(\mathcal{A}(Ms,Q_0))} \phi(x)P_A^2(x)\langle x-Q_0, \nabla\psi(x,\phi(x))\rangle dx\right|$$
$$+ C(n)M^2 s^2 \int_{\Pi(\mathcal{A}(Ms,Q_0))^c\cap L(M\mathcal{K}s,Q_0)} P_A^2(x)|\nabla\psi(x,\phi(x))|dx$$
$$\leq C(n)\varepsilon' M^2 s^2 \int_{\Pi(\mathcal{A}(Ms,Q_0))} P_A^2(x)|\nabla\psi(x,\phi(x))|dx$$
$$+ C(n)M^2 s^{2-2n-1}|\,(\Pi(\mathcal{A}(Ms,Q_0)))^c \cap B(2Ms,Q_0)\cap L(M\mathcal{K}s,Q_0)|$$
$$\leq C(n)M^{n+2}(e^{-C/\delta}+\varepsilon')s^{-(n-1)}.$$

Combining (6.21), (6.37), (6.38), (6.39), (6.40), (6.41), (6.42), (6.43) and (6.44) we conclude that

$$\left|\int_{\mathcal{A}(Ms,Q_0)} \psi(Q)P_A^2(\Pi(Q))\langle Q-A, \overrightarrow{n}(Q)\rangle d\sigma(Q)\right.$$
$$\left. -s\int_{L(M\mathcal{K}s,Q_0)} \psi(x,\phi(x))P_A^2(x)dx\right|$$
$$\leq C(n)M^{n+2}(e^{-C/\delta}+\varepsilon')s^{-(n-1)}.$$

In order to prove (6.26) we look at

(6.45)
$$\int_{\mathcal{A}(Ms,Q_0)} P_A(\Pi(Q))\frac{\psi(Q)}{|Q-A|^{n-1}}d\sigma(Q)$$
$$= \int_{\Pi(\mathcal{A}(Ms,Q_0))} P_A(x)\frac{\psi(x,\phi(x))}{|(x,\phi(x))-A|^{n-1}}\sqrt{1+|\nabla\phi(x)|^2}dx.$$

For $Q\in\mathcal{A}(Ms,Q_0)$, $Q = x+\phi(x)\overrightarrow{n}(M\mathcal{K}s,Q_0)$, $x\in\Pi(\mathcal{A}(Ms,Q_0))$. Thus

$$\left|\frac{1}{|(x,\phi(x))-A|^{n-1}}-\frac{1}{|x-A|^{n-1}}\right| = \frac{\left||x-A|^{n-1}-|(x,\phi(x))-A|^{n-1}\right|}{|(x,\phi(x))-A|^{n-1}|x-A|^{n-1}}$$
$$\leq C(n)\frac{|\phi(x)|\{|x-A|^{n-2}+|(x,\phi(x))-A|^{n-2}\}}{|(x,\phi(x))-A|^{n-1}|x-A|^{n-1}}$$
$$\leq C(n)|\phi(x)|\left\{\frac{1}{|x-A|\,|(x,\phi(x))-A|^{n-1}}+\frac{1}{|(x,\phi(x))-A|\,|x-A|^{n-1}}\right\}.$$

Note that for $x \in \Pi(\mathcal{A}(Ms,Q_0))$, $|\phi(x)| \leq \varepsilon' Ms$; therefore, provided that $\varepsilon' M < 1/4$, we have that $|(x,\phi(x))-A| \geq \frac{1}{4}s$,

$$\left|\frac{1}{|(x,\phi(x))-A|^{n-1}}-\frac{1}{|x-A|^{n-1}}\right| \leq C(n)M\varepsilon' s^{-(n-1)}. \tag{6.46}$$

Combining (6.45) and (6.46) we have

$$\left|\int_{\mathcal{A}(Ms,Q_0)} P_A(\Pi(Q))\frac{\psi(Q)}{|Q-A|^{n-1}}d\sigma(Q) - \int_{\Pi(\mathcal{A}(Ms,Q_0))} P_A(x)\frac{\psi(x,\phi(x))}{|x-A|^{n-1}}\sqrt{1+|\nabla\phi(x)|^2}dx\right| \tag{6.47}$$
$$\leq C(n)M\varepsilon' s^{-(n-1)}\int_{\Pi(\mathcal{A}(Ms,Q_0))} P_A(x)\psi(x,\phi(x))\sqrt{1+|\nabla\phi(x)|^2}dx$$
$$\leq C(n)M^{n+1}\varepsilon' s^{-(n-1)}.$$

On the other hand

$$\int_{\Pi(\mathcal{A}(Ms,Q_0))} P_A(x)\frac{\psi(x,\phi(x))}{|x-A|^{n-1}}\sqrt{1+|\nabla\phi(x)|^2}dx \tag{6.48}$$
$$= \int_{L(M\mathcal{K}s,Q_0)} P_A(x)\frac{\psi(x,\phi(x))}{|x-A|^{n-1}}\sqrt{1+|\nabla\phi(x)|^2}dx$$
$$- \int_{\Pi(\mathcal{A}(Ms,Q_0))^c\cap L(M\mathcal{K}s,Q_0)} P_A(x)\frac{\psi(x,\phi(x))}{|x-A|^{n-1}}\sqrt{1+|\nabla\phi(x)|^2}dx.$$

Therefore combining (6.47), (6.48) and (6.22) we obtain

$$\left|\int_{\mathcal{A}(Ms,Q_0)} P_A(\Pi(Q))\frac{\psi(Q)}{|Q-A|^{n-1}}d\sigma(Q) - \int_{L(M\mathcal{K}s,Q_0)} P_A(x)\frac{\psi(x,\phi(x))}{|x-A|^{n-1}}\sqrt{1+|\nabla\phi(x)|^2}dx\right|$$
$$\leq C(n)M^{n+1}(\varepsilon' + e^{-C/\delta})s^{-(n-1)}.$$

Choosing $(\varepsilon')^\gamma < \varepsilon$, $\varepsilon' < \varepsilon$, and $\varepsilon' M < 1/4$ we conclude the proof of Theorem 6.1. □

Theorem 6.2. *Let $\Omega$ be a $(\delta_1,\infty)$-chord arc domain. There exists $\delta_0 \in (0,\delta_n)$ such that if $\Omega$ is a $\delta$-Semmes decomposable domain with $\delta \in (0,\delta_0)$, if $\omega$ is asymptotically optimally doubling, and if* $\log h \in \mathrm{VMO}(d\sigma)$, *for each compact set $K \subset \mathbf{R}^{n+1}$ the following statement is true: Given $\varepsilon > 0$ and $\eta > 0$, there exists $M(\eta) > 1$ such that, for $M \geq M(\eta)$ there exists $\mathcal{K}(\varepsilon, M) > 4$ so that, for $\mathcal{K} \geq \mathcal{K}(\varepsilon, M)$ there is $s(\varepsilon,\eta,M,\mathcal{K}) > 0$ such that, for every $s \in (0, s(\varepsilon,\eta,M,\mathcal{K}))$, $Q_0 \in \partial\Omega \cap K$, and every $\psi \in C_c^\infty(B(2Ms,Q_0))$, $0 \leq \psi \leq 1$ and $\psi \equiv 1$ on $B(Ms,Q_0)$, if $\Pi$ denotes the orthogonal projection onto $L(M\mathcal{K}s,Q_0)$ then*

$$\begin{aligned}(6.49)\qquad &\Big|\int_{\partial\Omega} P_A(\Pi(Q))\psi(Q)\frac{d\sigma(Q)}{|Q-A|^{n-1}}\\ &\quad - \frac{\alpha}{\sigma_n}\int_{\partial\Omega}\psi(Q)P_A^2(\Pi(Q))\langle Q-A, \overrightarrow{n}(Q)\rangle d\sigma(Q)\Big|\\ &\quad \leq\ (C(n)M^{n+1}(\varepsilon + e^{-C/\delta}) + \eta)s^{-(n-1)},\end{aligned}$$

*where* $\alpha = \omega_n(2Ms)^n/\sigma(\Delta(2Ms,Q_0))$.

*Proof.* Let $\varepsilon' = \varepsilon'(\varepsilon) > 0$ and $\eta' = \eta'(\eta) > 0$, to be chosen. Choose $M(\eta') > 1$ so that (6.17) holds (for $M \geq M(\eta')$ and $s > 0$), with $\eta'$ instead of $\eta$. For $M \geq M(\eta')$, there exists $\mathcal{K}(\varepsilon', M) > 4$ so that for $\mathcal{K} \geq \mathcal{K}(\varepsilon', M)$ there is $s(\varepsilon',\eta',M,\mathcal{K}) \in (0, s(\eta'))$ so that Theorem 6.1 holds. Then, when $\alpha \geq 1/2$,

$$\begin{aligned}&\Big|\int_{\partial\Omega} P_A(\Pi(Q))\psi(Q)\frac{d\sigma(Q)}{|Q-A|^{n-1}} - \frac{\alpha}{\sigma_n}\int \psi(Q)P_A^2(\Pi(Q))\langle Q-A, \overrightarrow{n}(Q)\rangle d\sigma(Q)\Big|\\ &\quad\leq\ 2\Big|\int_{\partial\Omega}\alpha P_A(\Pi(Q))\psi(Q)\frac{d\sigma(Q)}{|Q-A|^{n-1}}\\ &\quad -\frac{\alpha^2}{\sigma_n}\int\psi(Q)P_A^2(\Pi(Q))\langle Q-A, \overrightarrow{n}(Q)\rangle d\sigma(Q)\Big|\\ &\quad\leq\ 2\Big|\int_{\partial\Omega}(k_A(Q) - \alpha P_A(\Pi(Q)))\psi(Q)\frac{d\sigma(Q)}{|Q-A|^{n-1}}\Big|\\ &\qquad +\ 2\Big|\int_{\partial\Omega}k_A(Q)\psi(Q)\frac{d\sigma(Q)}{|Q-A|^{n-1}} - \frac{1}{\sigma_n}\int_{\partial\Omega}k_A^2\psi(Q)\langle Q-A, \overrightarrow{n}(Q)\rangle d\sigma(Q)\Big|\\ &\qquad +\ \frac{2}{\sigma_n}\Big|\int_{\partial\Omega}(k_A^2 - \alpha^2P_A^2(\Pi(Q)))\langle Q-A, \overrightarrow{n}(Q)\rangle d\sigma(Q)\Big|.\end{aligned}$$

To estimate the first term we combine (6.19), (6.22) and (6.24):

$$\begin{aligned}(6.50)\ &\Big|\int_{\partial\Omega}(k_A(Q) - \alpha P_A(\Pi(Q)))\psi(Q)\frac{d\sigma(Q)}{|Q-A|^{n-1}}\Big|\\ &\qquad\leq\ \int_{\mathcal{A}(Ms,Q_0)}|k_A(Q) - \alpha P_A(\Pi(Q))|\,\psi(Q)\frac{d\sigma(Q)}{|Q-A|^{n-1}}\end{aligned}$$

$$+ \int_{\mathcal{F}(Ms,Q_0)} |k_A(Q) - \alpha P_A(\Pi(Q))|\psi(Q)\frac{d\sigma(Q)}{|Q-A|^{n-1}}$$
$$\leq C(n)M^n(e^{-C/\delta} + \varepsilon')s^{-(n-1)}.$$

To estimate the third term we combine (6.20), (6.21) and (6.23):

$$\left|\int_{\partial\Omega}(k_A^2 - \alpha^2 P_A^2(\Pi(Q)))\langle Q - A, \overrightarrow{n}(Q)\rangle d\sigma(Q)\right| \tag{6.51}$$
$$\leq \int_{\mathcal{A}(Ms,Q_0)} |k_A^2(Q) - \alpha^2 P_A^2(\Pi(Q))|\,|\langle Q - A, \overrightarrow{n}(Q)\rangle d\sigma(Q)$$
$$+ \int_{\mathcal{F}(Ms,Q_0)} |k_A^2(Q) - \alpha^2 P_A^2(\Pi(Q))|\,|Q - A|d\sigma(Q)$$
$$\leq C(n)M^n(e^{-C/\delta} + \varepsilon')s^{-(n-1)}.$$

Combining (6.17), (6.50) and (6.51) we conclude that as long as $\varepsilon' = \varepsilon$ and $\eta' = \eta/4$, (6.49) holds. □

Corollary 6.4. *Let $\Omega$ be a $(\delta_1, \infty)$-chord arc domain. There exists $\delta_0 \in (0, \delta_n)$ such that if $\Omega$ is a $\delta$-Semmes decomposable domain, with $\delta \in (0, \delta_0)$, if $\omega$ is asymptotically optimally doubling, and if $\log h \in \mathrm{VMO}(d\sigma)$, for each compact set $K \subset \mathbf{R}^{n+1}$ the following statement is true: Given $\varepsilon > 0$ and $\eta > 0$, there exists $M(\eta) > 1$ such that, for $M \geq M(\eta)$ there exists $\mathcal{K}(\varepsilon, M) > 4$ so that, for $\mathcal{K} \geq \mathcal{K}(\varepsilon, M)$ there is $s(\varepsilon, \eta, M, \mathcal{K}) > 0$ such that for every $s \in (0, s(\varepsilon, \eta, M, \mathcal{K}))$, and $Q_0 \in \partial\Omega \cap K$, there exist an $n$-dimensional plane $L(M\mathcal{K}s, Q_0)$ containing $Q_0$ and a Lipschitz function $\phi : L(M\mathcal{K}s, Q_0) \to \mathbf{R}$ such that, for every $\psi \in C_c^\infty(B(2Ms, Q_0))$, $0 \leq \psi \leq 1$ and $\psi \equiv 1$ on $B(Ms, Q_0)$, if $\Pi$ denotes the orthogonal projection onto $L(M\mathcal{K}s, Q_0)$ then*

$$\left|\frac{\alpha}{\sigma_n}\int_{L(M\mathcal{K}s,Q_0)} s\psi(x,\phi(x))P_A^2(x)dx \right. \tag{6.52}$$
$$\left. - \int_{L(M\mathcal{K}s,Q_0)} \psi(x,\phi(x))\frac{P_A(x)}{|x-A|^{n-1}}\sqrt{1 + |\nabla\phi(x)|^2}dx\right|$$
$$\leq (C(n)M^{n+2}(\varepsilon + e^{-C/\delta}) + \eta)s^{-(n-1)},$$

*where $\alpha = \omega_n(2Ms)^n/\sigma(\Delta(2Ms, Q_0))$.*

*Proof.* Assume that given $\varepsilon > 0$ and $\eta > 0$, $M(\eta) > 1$, $\mathcal{K}(\varepsilon, M) > 4$ for $M \geq M(\eta)$, and $s(\varepsilon, \eta, M, \mathcal{K})$ for $\mathcal{K} \geq K(\varepsilon, M)$, have been chosen, so that Theorem 6.1 and Theorem 6.2 hold. Then

$$\begin{aligned}
&\left|\frac{\alpha}{\sigma_n}\int_{L(M\mathcal{K}s,Q_0)} s\psi(x,\phi(x))P_A^2(x)dx\right.\\
&\qquad\left.-\int_{L(M\mathcal{K}s,Q_0)}\psi(x,\phi(x))\frac{P_A(x)}{|x-A|^{n-1}}\sqrt{1+|\nabla\phi(x)|^2}dx\right|\\
&\le\ \frac{\alpha}{\sigma_n}\left|\int_{L(M\mathcal{K}s,Q_0)} s\psi(x,\phi(x))P_A^2(x)dx\right.\\
&\qquad\qquad\left.-\int_{\mathcal{A}(Ms,Q_0)}\psi(Q)P_A^2(\Pi(Q))\langle Q-A,\overrightarrow{n}(Q)\rangle d\sigma(Q)\right|\\
&\quad+\ \left|\frac{\alpha}{\sigma_n}\int_{\mathcal{A}(Ms,Q_0)}\psi(Q)P_A^2(\Pi(Q))\langle Q-A,\overrightarrow{n}(Q)\rangle d\sigma(Q)\right.\\
&\qquad\qquad\left.-\int_{\mathcal{A}(Ms,Q_0)}P_A(\Pi(Q))\frac{\psi(Q)}{|Q-A|^{n-1}}d\sigma(Q)\right|\\
&\quad+\ \left|\int_{\mathcal{A}(Ms,Q_0)}P_A(\Pi(Q))\frac{\psi(Q)}{|Q-A|^{n-1}}d\sigma(Q)\right.\\
&\qquad\qquad\left.-\int_{L(M\mathcal{K}s,Q_0)}\psi(x,\phi(x))\frac{P_A(x)}{|x-A|^{n-1}}\sqrt{1+|\nabla\phi(x)|^2}dx\right|.
\end{aligned}$$

The first and the third terms are controlled by (6.25) and (6.26) respectively. We now focus on the second term

$$\begin{aligned}
(6.53)\quad&\left|\frac{\alpha}{\sigma_n}\int_{\mathcal{A}(Ms,Q_0)}\psi(Q)P_A^2(\Pi(Q))\langle Q-A,\overrightarrow{n}(Q)\rangle d\sigma(Q)\right.\\
&\quad\left.-\int_{\mathcal{A}(Ms,Q_0)}P_A(\Pi(Q))\frac{\psi(Q)}{|Q-A|^{n-1}}d\sigma(Q)\right|\\
&\quad\le\ \left|\int_{\partial\Omega}P_A(\Pi(Q))\psi(Q)\frac{d\sigma(Q)}{|Q-A|^{n-1}}\right.\\
&\qquad\quad\left.-\frac{\alpha}{\sigma_n}\int_{\partial\Omega}\psi(Q)P_A^2(\Pi(Q))\langle Q-A,\overrightarrow{n}(Q)\rangle d\sigma(Q)\right|\\
&\qquad+\ \frac{\alpha}{\sigma_n}\int_{\mathcal{F}(Ms,Q_0)}\psi(Q)P_A^2(\Pi(Q))|\langle Q-A,\overrightarrow{n}(Q)\rangle|d\sigma(Q)\\
&\qquad+\ \int_{\mathcal{F}(Ms,Q_0)}P_A(\Pi(Q))\frac{\psi(Q)}{|Q-A|^{n-1}}d\sigma(Q).
\end{aligned}$$

Note that the first term in (6.53) is controlled by (6.49). Moreover, our choice of $A$ guarantees that, for $Q\in\mathcal{F}(Ms,Q_0)$, $\frac{s}{16}\le|Q-A|\le 2Ms$ (see (5.12));

therefore

$$\int_{\mathcal{F}(Ms,Q_0)} \psi(Q) P_A^2(\Pi(Q))|Q-A| d\sigma(Q) + \int_{\mathcal{F}(Ms,Q_0)} P_A(\Pi(Q)) \frac{\psi(Q)}{|Q-A|^{n-1}} d\sigma(Q) \leq c(n) M^{n+1}(e^{-C/\delta}+\varepsilon) s^{-(n-1)}. \tag{6.54}$$

Combining (6.49), (6.53) and (6.54) we conclude that (6.52) holds. □

*Remark* 6.1. Note that $\sup_{B(2Ms,Q_0)\cap L(M\mathcal{K}s,Q_0)} |\phi| \leq \frac{1}{2}s$ (by our choice of $\varepsilon > 0$ and $\delta_0 > 0$). Therefore, an argument similar to the one presented in the proof of Corollaries 6.2 and 6.3 guarantees that, for $\eta > 0$ there exists $M(\eta) > 0$ such that, for $M \geq M(\eta)$, $s > 0$, and $\psi \in C_c^\infty(B(2Ms,Q_0))$, $0 \leq \psi \leq 1$ and $\psi \equiv 1$ on $B(Ms,Q_0)$,

$$\left| \int_{L(M\mathcal{K}s,Q_0)} \psi(x,\phi(x)) \frac{P_A(x)}{|x-A|^{n-1}} dx - \frac{1}{\sigma_n} \int_{L(M\mathcal{K}s,Q_0)} s P_A^2(x) \psi(x,\phi(x)) dx \right| \leq \eta s^{-(n-1)}.$$

Combining Corollary 6.4 and Remark 6.1, we conclude that:

Corollary 6.5. *Let $\Omega$ be a $(\delta_1,\infty)$-chord arc domain. There exists $\delta_0 \in (0,\delta_n)$ such that if $\Omega$ is a $\delta$-Semmes decomposable domain with $\delta \in (0,\delta_0)$, if $\omega$ is asymptotically optimally doubling, and if $\log h \in \mathrm{VMO}(d\sigma)$, for each compact set $K \subset \mathbf{R}^{n+1}$ the following statement is true: Given $\varepsilon > 0$ and $\eta > 0$, there exists $M(\eta) > 1$ such that, for $M \geq M(\eta)$ there exists $\mathcal{K}(\varepsilon, M) > 4$ so that, for $\mathcal{K} \geq \mathcal{K}(\varepsilon,M)$ there is $s(\varepsilon,\eta,M,\mathcal{K}) > 0$ such that, for every $s \in (0, s(\varepsilon,\eta,M,\mathcal{K}))$, and $Q_0 \in \partial\Omega \cap K$, there exist an $n$-dimensional plane $L(M\mathcal{K}s,Q_0)$ containing $Q_0$ and a Lipschitz function $\phi : L(M\mathcal{K}s,Q_0) \to \mathbf{R}$ such that for every $\psi \in C_c^\infty(B(2Ms,Q_0))$, $0 \leq \psi \leq 1$ and $\psi \equiv 1$ on $B(Ms,Q_0)$, if $\Pi$ denotes the orthogonal projection onto $L(M\mathcal{K}s,Q_0)$ then*

$$\left| \int_{L(M\mathcal{K}s,Q_0)} \psi(x,\phi(x)) \frac{P_A(x)}{|x-A|^{n-1}} (\alpha - \sqrt{1+|\nabla\phi(x)|^2}) dx \right| \leq C(n) M^{n+2}(e^{-C/\delta}+\varepsilon) s^{-(n-1)} + \eta s^{-(n-1)},$$

*where*

$$\alpha = \frac{\omega_n (2Ms)^n}{\sigma(\Delta(2Ms,Q_0))}.$$

*Remark* 6.2. Note that since $\omega$ is asymptotically optimally doubling, and $\Omega$ is a Reifenberg flat domain, Theorem 3.4 guarantees that $\Omega$ is a Reifenberg flat domain with vanishing constant. Therefore given $\beta > 0$ and $K \subset \mathbf{R}^{n+1}$ a compact set, there exists $r(\beta) > 0$ such that, for $r \in (0, r(\beta))$ and $Q \in \partial\Omega \cap K$

$$\sigma(B(r,Q)) \geq \frac{1}{1+\beta}\omega_n r^n \tag{6.55}$$

(see Remark 4.1). Hence, if $2Ms \in (0, r(\beta))$, $\alpha \leq 1 + \beta$.

COROLLARY 6.6. *Let $\Omega$ be a $(\delta_1, \infty)$-chord arc domain. There exists $\delta_0 \in (0, \delta_n)$ such that, if $\Omega$ is a $\delta$-Semmes decomposable domain with $\delta \in (0, \delta_0)$, if $\omega$ is asymptotically optimally doubling, and if* $\log h \in \mathrm{VMO}(d\sigma)$, *for each compact set $K \subset \mathbf{R}^{n+1}$ the following statement is true: Given $\varepsilon > 0$ and $\eta > 0$, there exists $M(\eta) > 1$ such that, for $M \geq M(\eta)$ there exists $s(\varepsilon, \eta, M) > 0$ such that, for every $s \in (0, s(\varepsilon, \eta, M))$, and $Q_0 \in \partial\Omega \cap K$, there exists $\overrightarrow{n}(s, Q_0) = \overrightarrow{n}(\varepsilon, \eta, M, s, Q_0) \in \mathbf{S}^n$ such that*

$$\frac{1}{\sigma(\Delta(2s,Q_0))}\int_{\Delta(2s,Q_0)} |\overrightarrow{n} - \overrightarrow{n}(s,Q_0)|^2 d\sigma \leq C(n)M^{n+2}(e^{-C/\delta} + \varepsilon) + \eta.$$

*Proof.* Fix $\varepsilon' = \varepsilon'(\varepsilon) > 0$, and $\eta' = \eta'(\eta) > 0$, to be chosen. Let $M(\eta') > 1$, $\mathcal{K}(\varepsilon', M) > 4$, for $M \geq M(\eta')$ and $s(\varepsilon', \eta', M, \mathcal{K}(\varepsilon', M)) > 0$ be such that Corollary 6.5 and (6.55) are satisfied for $s \in (0, s(\varepsilon', \eta', M, \mathcal{K}(\varepsilon', M))) > 0$ and $\beta = \varepsilon'$. Let $\overrightarrow{n}(s, Q_0) = \overrightarrow{n}(M\mathcal{K}(\varepsilon', M)s, Q_0)$ be the unit normal vector to $L(M\mathcal{K}(\varepsilon', M)s, Q_0)$. Since Corollary 6.5 holds if $\psi \in C_c^\infty(B(2Ms, Q_0))$, $0 \leq \psi \leq 1$ and $\psi \equiv 1$ on $B(Ms, Q_0)$,

$$\left|\int_{L(M\mathcal{K}(\varepsilon',M)s,Q_0)} \psi(x,\phi(x))\frac{P_A(x)}{|x-A|^{n-1}}(\alpha - \sqrt{1+|\nabla\phi(x)|^2})dx\right|$$
$$\leq C(n)M^{n+2}(e^{-C/\delta} + \varepsilon')s^{-(n-1)} + \eta' s^{-(n-1)}$$

Now (6.55) implies that $\alpha \leq 1 + \varepsilon'$. In particular

$$\begin{aligned}\int_{L(M\mathcal{K}(\varepsilon',M)s,Q_0)} &\psi(x,\phi(x))\frac{P_A(x)}{|x-A|^{n-1}}(\sqrt{1+|\nabla\phi(x)|^2} - (1+\varepsilon'))dx \\ &\leq \int_{L(M\mathcal{K}(\varepsilon',M)s,Q_0)} \psi(x,\phi(x))\frac{P_A(x)}{|x-A|^{n-1}}(\sqrt{1+|\nabla\phi(x)|^2} - \alpha)dx \\ &\leq C(n)M^{n+2}(e^{-C/\delta} + \varepsilon')s^{-(n-1)} + \eta' s^{-(n-1)}.\end{aligned} \tag{6.56}$$

Thus

$$\begin{aligned}\int_{L(M\mathcal{K}(\varepsilon',M)s,Q_0)} &\psi(x,\phi(x))\frac{P_A(x)}{|x-A|^{n-1}}(\sqrt{1+|\nabla\phi(x)|^2} - 1)dx \\ &\leq C(n)M^{n+2}(e^{-C/\delta} + \varepsilon')s^{-(n-1)} + \eta' s^{-(n-1)}\end{aligned} \tag{6.57}$$

$$
\begin{aligned}
&+ \ \varepsilon' \int_{L(M\mathcal{K}(\varepsilon',M)s,Q_0)} \psi(x,\phi(x)) \frac{P_A(x)}{|x-A|^{n-1}} dx \\
&\le \ C(n) M^{n+2}(e^{-C/\delta}+\varepsilon') s^{-(n-1)} + \eta' s^{-(n-1)}.
\end{aligned}
$$

Note that

$$
\begin{aligned}
\sqrt{1+|\nabla\phi(x)|^2} - 1 &= \sqrt{1+|\nabla\phi(x)|^2}\left(1 - \frac{1}{\sqrt{1+|\nabla\phi(x)|^2}}\right) \\
&= \sqrt{1+|\nabla\phi(x)|^2}\frac{1}{2}|\overrightarrow{n}(x,\phi(x)) - \overrightarrow{n}(s,Q_0)|^2
\end{aligned}
$$

where $\overrightarrow{n}(x,\phi(x)) = \frac{(-\nabla\phi(x),1)}{\sqrt{1+|\nabla\phi(x)|^2}}$. Moreover,

(6.58)

$$
\begin{aligned}
&\int_{L(M\mathcal{K}(\varepsilon,M)s,Q_0)} \psi(x,\phi(x)) \frac{P_A(x)}{|x-A|^{n-1}} (\sqrt{1+|\nabla\phi(x)|^2} - 1) dx \\
&\quad = \ \frac{1}{2}\int_{L(M\mathcal{K}(\varepsilon,M)s,Q_0)} \psi(x,\phi(x)) \frac{P_A(x)}{|x-A|^{n-1}} \\
&\quad\quad \times \ \sqrt{1+|\nabla\phi(x)|^2}|\overrightarrow{n}(x,\phi(x)) - \overrightarrow{n}(s,Q_0)|^2 dx \\
&\quad \ge \ \frac{1}{2}\int_{\Pi(\Delta(2s,Q_0))} \frac{P_A(x)}{|x-A|^{n-1}}|\overrightarrow{n}(x,\phi(x)) - \overrightarrow{n}(s,Q_0)|^2\sqrt{1+|\nabla\phi(x)|^2} dx \\
&\quad \ge \ \frac{C(n)}{s^{2n-1}}\int_{\Delta(2s,Q_0)\cap \operatorname{graph}\phi} |\overrightarrow{n}(Q) - \overrightarrow{n}(s,Q_0)|^2 d\sigma(Q)
\end{aligned}
$$

where $\overrightarrow{n}(Q)$ denotes the unit normal to graph $\phi$. Combining (6.57) and (6.58) we have that

$$
\frac{1}{s^n}\int_{\Delta(2s,Q_0)\cap \operatorname{graph}\phi} |\overrightarrow{n}(Q) - \overrightarrow{n}(s,Q_0)|^2 d\sigma(Q) \le c(n) M^{n+2}(e^{-C/\delta}+\varepsilon') + C(n)\eta'.
$$

Since $\mathcal{A}(Ms,Q_0) \subset \operatorname{graph}\phi$, our choice of orientation guarantees that the unit normal to the boundary $\partial\Omega$ and the unit normal to graph $\phi$ coincide on $\mathcal{A}(Ms,Q_0)$. Therefore,

$$
\begin{aligned}
&\frac{1}{s^n}\int_{\Delta(2s,Q_0)} |\overrightarrow{n}(Q) - \overrightarrow{n}(s,Q_0)|^2 d\sigma(Q) \\
&\qquad \le \ \frac{1}{s^n}\int_{\Delta(2s,Q_0)\cap \operatorname{graph}\phi} |\overrightarrow{n}(Q) - \overrightarrow{n}(s,Q_0)|^2 d\sigma(Q) \\
&\qquad\quad + \ \frac{1}{s^n}\int_{\mathcal{F}(Ms,Q_0)} |\overrightarrow{n}(Q) - \overrightarrow{n}(s,Q_0)|^2 d\sigma(Q) \\
&\qquad \le \ C(n) M^{n+2}(e^{-C/\delta}+\varepsilon') + C(n)\eta'.
\end{aligned}
$$

Choosing $\varepsilon' = \varepsilon$ and $\eta' = \eta/2C(n)$, we conclude that

$$\frac{1}{\sigma(\Delta(2s,Q_0))}\int_{\Delta(2s,Q_0)} |\overrightarrow{n}(Q) - \overrightarrow{n}(s,Q_0)|^2 d\sigma(Q) \leq C(n)M^{n+2}(e^{-C/\delta} + \varepsilon) + \eta.$$

□

*Proof of the Main Lemma.* Let $\Omega$ be a $(\delta_1,\infty)$-chord arc domain. Assume that $\omega$ is asymptotically optimally doubling and that $\log h \in \mathrm{VMO}(d\sigma)$. Let $\delta_0 > 0$, and $\delta \in (0,\delta_0)$. Assume that for each compact set $K \subset \mathbf{R}^{n+1}$ there exists $R > 0$ such that

$$\sup_{Q\in K\cap\partial\Omega} \|\overrightarrow{n}\|_*(B(R,Q)) \leq \delta.$$

In this case $\Omega$ is $\sqrt{\delta}$-Semmes decomposable (see Theorem 4.1). In particular, Corollary 6.6 holds for $\varepsilon = e^{-c/\sqrt{\delta}}$, $\eta = \delta^4$ and $M > M(\eta) = \left(\frac{C(n)}{\delta^4}\right)^{\frac{1}{n-1+\alpha}}$. Note that this condition on $M$ comes from the proofs of Lemma 6.1 and Lemma 6.2 (see inequalities (6.9) and (6.16)). For $\delta > 0$ small enough we can choose $M = \delta^{-4}$. Then, Corollary 6.6 guarantees that, for each compact set $K \subset \mathbf{R}^{n+1}$, there exists $r(R,\delta) > 0$ such that,

$$\begin{aligned}\sup_{Q\in\partial\Omega\cap K} \|\overrightarrow{n}\|_*(B(r,Q)) &\leq C(n)M^{\frac{n+2}{2}}e^{-c/\sqrt{\delta}} + \delta^2 \\ &\leq C(n)\delta^{-2(n+2)}e^{-c\sqrt{\delta}} + \delta^2.\end{aligned}$$

Choose $\delta_0 > 0$ small enough so that, for $\delta \leq \delta_0$

$$C(n)\delta^{-4(n+2)}e^{-c/\sqrt{\delta}} \leq \delta^2 \leq \frac{1}{4}\delta.$$

We conclude that there exists $\delta_0(n) > 0$ so that, if $\delta \in (0,\delta_0)$, for each compact set $K \subset \mathbf{R}^{n+1}$, there exists $R > 0$ so that if

$$\sup_{Q\in K\cap\partial\Omega} \|\overrightarrow{n}\|_*(B(R,Q)) \leq \delta$$

then there exists $r > 0$ for which

$$\sup_{Q\in\partial\Omega\cap K} \|\overrightarrow{n}\|_*(B(r,Q)) \leq \frac{1}{2}\delta.$$

□

The main theorem as well as Theorems 5.3 and 5.4 have quantitative versions. Their proofs follow the exact same lines and require careful accounting of the constants. In order to be able to state the quantitative results we need one last definition.

*Definition* 6.1. Let $\Omega$ be a $\delta$-chord arc domain. When $f \in L^2_{\mathrm{loc}}(d\sigma)$, we say that $f \in \mathrm{BMO}(\partial\Omega)$ if for each compact set $K \subset \mathbf{R}^{n+1}$, there exists $R > 0$ such that

$$\sup_{Q \in \partial\Omega \cap K} \|f\|_*(B(R,Q)) < \infty$$

where

$$\|f\|_*(B(R,Q)) = \left( \sup_{0<s\le R} \frac{1}{\sigma(\Delta(s,Q))} \int_{\Delta(s,Q)} |f - f_{s,Q}|^2 d\sigma \right)^{\frac{1}{2}}.$$

For $\varepsilon > 0$ and $f \in \mathrm{BMO}(\partial\Omega)$, we say that $\|f\|_* \le \varepsilon$, if for each compact set $K \subset R^{n+1}$, there exists $R > 0$ such that

$$\sup_{Q \in \partial\Omega \cap K} \|f\|_*(B(R,Q)) \le \varepsilon.$$

THEOREM 6.3. *Let $\Omega \subset \mathbf{R}^{n+1}$ be a $(\delta_1, \infty)$-chord arc domain. There exists $\delta = \delta(n) > 0$ so that if $\Omega$ is a $\delta$-chord arc domain, then given $\varepsilon > 0$, there exists $\eta > 0$ such that, if $\omega$ is $\eta$-approximately optimally doubling and $\|\log h\|_* \le \eta$, then $\|\overrightarrow{n}\|_* \le \varepsilon$. Here $\omega$ denotes the harmonic measure either with pole at infinity or with finite pole in $\Omega$, $h = \frac{d\omega}{d\sigma}$ denotes its Poisson kernel, and $\overrightarrow{n}$ denotes the unit normal vector to $\partial\Omega$.*

THEOREM 6.4. *There exists $\delta = \delta(n) > 0$ so that if $\Omega$ is a bounded $\delta$-chord arc domain and $X \in \Omega$, then given $\varepsilon > 0$, there exists $\eta > 0$, such that, if $\omega^X$ is $\eta$-approximately optimally doubling, and $\|\log k_X\|_* \le \eta$, then $\|\overrightarrow{n}\|_* \le \varepsilon$.*

UNIVERSITY OF CHICAGO, CHICAGO, IL
*E-mail address*: cek@math.uchicago.edu

UNIVERSITY OF WASHINGTON, SEATTLE, WA
*E-mail address*: toro@math.washington.edu